\documentclass[12pt,xcolor=pdftex,a4paper,oneside,centertags,reqno]{amsart}

\usepackage{amsmath,amsthm,amssymb,amsfonts,amscd}
\usepackage{latexsym}
\usepackage{graphpap}
\usepackage{epsfig}
\usepackage{amscd}
\usepackage{amsbsy}
\usepackage{enumerate}
\usepackage{eucal}
\usepackage[usenames,dvipsnames]{color}
\usepackage[table, dvipsnames]{xcolor}
\usepackage{tikz}
\usepackage{pict2e,color}
\usepackage{bookmark}
\usepackage{bm}  

\usepackage{needspace}

\usepackage{times}
\usepackage{makecell}

\usepackage{endnotes}
\usepackage{epsfig}
\usepackage{eucal} 
\usepackage{latexsym}
\usepackage{pdfsync}
\usepackage{mathrsfs}
\usepackage[shortlabels]{enumitem} 

\makeatletter
\renewcommand*\env@matrix[1][c]{\hskip -\arraycolsep
  \let\@ifnextchar\new@ifnextchar
  \array{*\c@MaxMatrixCols #1}}
\makeatother

\numberwithin{equation}{section} 

\usepackage[francais]{[babel}
\usepackage[deutsch]{[babel}

\usepackage{graphics}
\usepackage{graphicx}

 \usepackage{hyperref}

\usepackage{pdfsync}
\usepackage{color}
\usepackage{endnotes}

\usepackage{mathrsfs}

\font\gotxii=eufm10 at 12pt

\font\posebni=msam10

\font\gotm=eufm10 at 8pt

\font\posebni=msam10

\newcommand{\3}[0]{\v{s}}

\newcommand{\C}[0]{{\mathbb C}}

\newcommand{\N}[0]{{\mathbb N}}

\newcommand{\R}[0]{\mathbb{R}}

\newcommand{\Z}[0]{{\mathbb Z}}

\newcommand{\bP}{{\mathbf P}}

\newcommand\bS{\mathbf{S}}

\newcommand{\cA}[0]{{\mathcal A}}
\newcommand{\cB}[0]{{\mathcal B}}

\newcommand{\cD}[0]{{\mathcal D}}
\newcommand{\cE}[0]{{\mathcal E}}

\newcommand{\cH}[0]{{\mathcal H}}
\newcommand{\cI}[0]{{\mathcal I}}
\newcommand{\cJ}[0]{{\mathcal J}}
\newcommand{\cK}[0]{{\mathcal K}}
\newcommand{\cL}[0]{{\mathcal L}}
\newcommand{\cM}[0]{{\mathcal M}}

\newcommand{\cR}[0]{{\mathcal R}}

\newcommand{\cT}{{\mathcal T}}
\newcommand{\cU}[0]{{\mathcal U}}
\newcommand{\cV}[0]{{\mathcal V}}

\newcommand{\oV}{{\mathscr V}}

\renewcommand{\a}[0]{\text{\gotxii a}}

\newcommand{\gotb}[0]{\text{\gotxii b}}
\newcommand{\am}[0]{\text{\gotm a}}

\newcommand{\gJ}{\text{\gotxii J}}

\newcommand{\leqsim}[0]{\,\text{\posebni \char46}\,}
\newcommand{\geqsim}[0]{\,\text{\posebni \char38}\,}

\newcommand{\vt}[0]{\vartheta}
\newcommand{\e}[0]{\varepsilon}
\newcommand{\f}[0]{\varphi}

\newcommand{\wt}{\widetilde}
\newcommand{\pd}[0]{\partial}

\newcommand{\nor}[1]{\| #1 \|}

\newcommand{\mn}[2]{\{ #1\, ;\, #2 \}}
\newcommand{\Mn}[2]{\left\{ #1 ; #2 \right\}}
\newcommand{\sk}[2]{\langle #1 , #2\rangle}
\newcommand{\Sk}[2]{\left\langle #1 , #2\right\rangle}

\newtheorem{theorem}{
Theorem}[section]

\newtheorem{proposition}[theorem]{Proposition}
\newtheorem{claim}[theorem]{Claim}
\newtheorem{lemma}[theorem]{Lemma}

\newtheorem{corollary}[theorem]{Corollary}

\renewcommand\leq[0]{\leqslant}
\renewcommand\geq[0]{\geqslant}
\renewcommand\epsilon[0]{\varepsilon}
\renewcommand\theta[0]{\vartheta}

\newcommand\norm[2]{{\left\Vert{#1}\right\Vert_{#2}}}

\renewcommand{\arraystretch}{1.4}

\renewcommand{\div}{\,\text{div}\,}

\newtheorem{predefinition}[theorem]{Definition}
\newenvironment{definition}%
{\begin{predefinition}\rm}{\end{predefinition}}

\newcounter{magnolia blues}

\newtheorem{prequestion}[magnolia blues]{Question}  %
{\begin{prequestion}\rm}{\end{prequestion}}

\newtheorem{preremark}[theorem]{Remark}  \newenvironment{remark}%
{\begin{preremark}\rm}{\end{preremark}}
\newtheorem{prenotation}[theorem]{Notation}  \newenvironment{notation}%
{\begin{prenotation}\rm}{\end{prenotation}}
\newtheorem{preexample}[theorem]{Example}  \newenvironment{example}%
{\begin{preexample}\rm}{\end{preexample}}

\begin{document}

\title[\MakeLowercase{\textit{p}}-ellipticity for complex systems]
{
\MakeLowercase{\textbf{\textit{p}}}-ellipticity for systems of elliptic equations with complex coefficients}


\author[Carbonaro]{Andrea Carbonaro}
\author[Dragi\v{c}evi\'c]{Oliver Dragi\v{c}evi\'c}

\address{Andrea Carbonaro\\ Universit\`a degli Studi di Genova\\ Dipartimento di Matematica\\ Via Dodecaneso\\ 35 16146 Genova\\ Italy}
\email{carbonaro@dima.unige.it}

\address{Oliver Dragičević, Department of Mathematics, Faculty of Mathematics and Physics, University of Ljubljana,  Jadranska ulica 19, SI-1000 Ljubljana, Slovenia, and Institute of Mathematics, Physics and Mechanics, Jadranska ulica 19, SI-1000 Ljubljana, Slovenia}
\email{oliver.dragicevic@fmf.uni-lj.si}


\dedicatory{\it The second-named author dedicates this paper, with gratitude,
\\ to Professor Josip Globevnik, on the occasion of his 80th birthday.}

\begin{abstract}
We extend the concept of $p$-ellipticity for (single) elliptic operators, that we introduced in \cite{CD-DivForm}, to the case of systems of elliptic equations with complex coefficients. 
The $p$-ellipticity condition from this paper turns out to be equivalent, up to a factor depending on $p$, to a condition for systems introduced (and equally named) by Dindoš, Li and Pipher in \cite{DiLiPi20}.
We prove several key properties of $p$-ellipticity akin to those that have been known to hold in the scalar case \cite{CD-DivForm}. In some respects, however, the vector case fundamentally differs from the scalar one. 
As an application we establish the contractivity on $L^p$ of the operator semigroup generated by systems of elliptic equations with complex coefficients
and subject to Dirichlet, Neumann or mixed boundary conditions. 
\end{abstract}

\maketitle

\setcounter{tocdepth}{1}
\tableofcontents
\label{toc}

\setcounter{section}{-1}
\section{Introduction}

The present paper builds upon and extends the notion of {\it $p$-ellipticity} that appeared a few years ago in the context of partial differential equations (PDEs) with complex coefficients. We start by summarizing it. 

\subsection{The concept of $p$-ellipticity for scalar elliptic partial differential equations in divergence form} 
\label{s: scalar p-ell}

Let $\Omega\subseteq\R^{n}$ be an arbitrary open set and $A:\Omega\rightarrow\C^{n\times n}$ a strictly uniformly elliptic (or accretive) matrix function with $L^{\infty}$ coefficients. In \cite{CD-DivForm} we introduced the concept of $p$-ellipticity as follows. Given $A$ as above and $p\in(1,\infty)$, 
let
\begin{equation}
\label{eq: Delta_p}
\Delta_{p}(A):=
\underset{x\in\Omega}{{\rm ess}\inf}
\min_{\substack{\xi\in\C^{n}\\ |\xi|=1}}
\Re\sk{A(x)\xi}{\xi+|1-2/p|\bar\xi}_{\C^{n}}.
\end{equation}
We say that $A$ is {\it $p$-elliptic} if $\Delta_{p}(A)>0$. The special case $p=2$ coincides with the notion of classical ellipticity. We used $\Delta_{p}(A)$ to study the behaviour on $L^p$ of elliptic partial differential operators in divergence form and with complex coefficients, formally $L_Au=-\div(A\nabla u)$, and associated objects.

Around the same time and independently of us, and through their work on problems different from those we considered, M. Dindo\v{s} and J. Pipher 
introduced a condition \cite[(1.3)]{DiPi} which turned out to be exactly equivalent to $\Delta_p(A)>0$.
They did so by strengthening an older condition by Cialdea and Maz'ya \cite{CM}. 
In recent years, through a considerable number of works by various authors, e.g. \cite{CD-DivForm, CD-Mixed, CD-Potentials, Trilinear, Pog24, Pog25, DiPi, DiPi20, DiLiPi20, E20, EHRT21, BER24}, $p$-ellipticity proved to be a fundamental condition that closely relates to two essential concepts in mathematical analysis, namely, elliptic partial differential equations with complex coefficients and $L^p$ spaces. On top of that, we showed it admits strong ties with holomorphic functional calculus on $L^p$ \cite{CD-mult, CD-OU}. 

Much more information about the genesis of $p$-ellipticity and properties of $p$-elliptic operators can be found in \cite{CD-DivForm, CD-Potentials, Trilinear}.

\subsection{The main goal} 

Recently, in a joint work with Li  \cite{DiLiPi20}, Dindoš and Pipher continued their very successful line of exploration of $p$-ellipticity in PDEs by considering systems of PDEs. In order to do so, they introduced a condition that generalizes $p$-ellipticity in the context of elliptic systems. 
See also Cialdea--Maz'ya \cite{CM06}.

In this paper we approach elliptic systems from our own, different point of view, which was in the scalar case already elaborated in our earlier works cited in Section \ref{s: scalar p-ell}. That is, we study {\it heat flows} associated with {\it power functions} (the latter entered the focus of our interest as building blocks of particular {\it Bellman functions}). This leads to our generalizing the numbers $\Delta_p(A)$ and, consequently, to the introduction of an algebraic condition for which we retain the name {\it $p$-ellipticity}, as we find it a natural extension of the $p$-ellipticity for scalar operators. This is the content of Section \ref{s: p-ell}, while the exposition of our path to $p$-ellipticity, and the reasoning that led to it, are deferred to Sections \ref{s: calculation} and \ref{s: rootsyst}. 

Concerning the scalar $p$-ellipticity condition, applications amount to finding different kinds of theorems  which entail $p$-ellipticity as the key assumption. In the vector case, the emphasis lies on extending the results from the scalar case to the vector case. 
Actually, the scalar case \cite{CD-DivForm} often serves as a guideline as to what properties we may aspire to.
Many of these generalizations may have the same formulation as in the scalar case, but proving them in the vector-valued context is often much more involved.

The principal applications of $p$-ellipticity for systems that we pursued whilst developing it here were twofold: the $L^p$ contractivity of the associated semigroup (Theorem \ref{t: contractivity}) and a version of the {\it bilinear embedding theorem} for $p$-elliptic systems (owing to its complexity and consequent length, it will be the subject of another paper). 

In Section \ref{s: comparison} we also elaborate on the Dindoš--Li--Pipher condition from \cite{DiLiPi20} and show it to be essentially equivalent to the one we introduce in this paper; see Proposition \ref{c: DiLiPi} for the precise statement. 
The main differences between the two designs (ours and the one by Dindo\3 et al.) are that:
(i) 
in order to derive $\Delta_p(A)$, we differentiate the flow of the semigroup (a special case of the Bellman function flow), which leads to a dependence on the sphere in $\C^m$ in one factor only; 
(ii) 
we perceive $\Phi_\omega\xi$, which is our counterpart to Dindo\3 et al.'s $\xi(\omega)$, as {\it projections depending on $\omega$ and applied to $\xi$}, which leads to a number of algebraic identities and inequalities and, eventually, properties of the $p$-ellipticity coefficients $\Delta_p(A)$.

\subsection{Organization of the paper}
We split the paper into four main parts, as specified in \hyperref[toc]{Table of Contents}.

\begin{itemize}
\itemsep1em
\item
In Part \ref{part: Prelims} we introduce the notation used throughout the paper, define the main objects and, finally, prove crucial linear algebra results on which later proofs shall rest.

\item
In Part \ref{part: Def Props} we introduce our $p$-ellipticity condition and prove a number of its properties, mostly based on those previously known to be true in the scalar case. We include a couple of examples of tensors when the $p$-ellipticity parameter can be explicitly calculated and use some of them to illustrate distinctions between the scalar and the vector case.

\item
In Part \ref{part: Roots} we present the reasoning that led us to define $p$-ellipticity the way we did in the previous part. Thus the principal purpose of this part is to reveal ideas fueling this paper.

\item
In Part \ref{part: Applications} we prove an extension of our previous result on semigroup contractivity \cite[Theorem 1.3]{CD-DivForm} to the vector case.

\end{itemize}

\part{Preliminaries}
\label{part: Prelims}

\section{Notation}
\label{s: notation}

Setting up consistent yet manageable notation in the context of partial differential equations is known to be a challenging task \cite[p.xx, e.]{Evans}.
This holds even more so when dealing with {\it systems} of differential equations. 
The aim of this section is to bring up reasonable notation and clarify the meaning of symbols used in this paper.
We shall to some extent follow Chen and Wu \cite{CW98} and Giaquinta and Martinazzi \cite{GM12}.

\subsection
{Complex \& real spaces. Identifications}
\label{s: complex vs. real}

For a given $k \in \N$, let $\zeta \in \C^k$ be written as $\zeta = x + i y$, where $x, y \in \R^k$ are the real and imaginary parts, respectively. We define the identification operator $\cV=\cV_k:\C^{k} \rightarrow \R^{2k}$ as follows:
$$
\cV(x+iy) = (x, y).
$$
We shall often write simply $\cV(\zeta) = \widetilde{\zeta}$. Thus 
$$
(x+iy)\wt\ =(x,y)\in\R^k\times\R^k\equiv\R^{2k}.
$$
Observe that the real part of the standard Hermitian product on $\C^k$ corresponds to the Euclidean scalar product on $\R^{2k}$ under the above identification:
\begin{equation}
\label{eq: real vs. complex inner product}
\Re\sk{\zeta}{\eta}_{\C^k} = \sk{\wt\zeta}{\wt\eta}_{\R^{2k}}. 
\end{equation}
Similarly, let $\mathsf{X}$ be any set and $F: \C^k \rightarrow \sf{X}$ a given function. We define its counterpart on $\R^{2k}$, denoted $\widetilde F$, 
via the composition 
$\widetilde F = F \circ \cV^{-1}_k :\R^k\times\R^k\rightarrow\sf{X}$. 
Explicitly,
\begin{equation*}
\widetilde F(\widetilde\zeta) = F(\zeta). 
\end{equation*}
In coordinates, this reads
$$
\widetilde F(x,y) = F(x+iy).
$$

The {\it scalar case} corresponds -- by definition -- to taking $k=1$ everywhere.

\subsection{Wirtinger derivatives}
Let $\Phi:\R^k\times\R^k\rightarrow\C$. 
With $z=(z_1,\hdots,z_k)$, we denote
$$
\aligned
\pd_z\Phi& :=(\pd_{z_1}\Phi,\hdots,\pd_{z_k}\Phi)\\
\pd_{\bar z}\Phi& :=(\pd_{\bar z_1}\Phi,\hdots,\pd_{\bar z_k}\Phi),
\endaligned
$$
where $\pd_{z_j},\pd_{\bar z_j}$ are the usual {\it Wirtinger derivatives}. That is, 
$$
\aligned
\pd_{z_j}&=\frac12(\pd_{x_j}-i\pd_{y_j})\\
\pd_{\bar z_j}&=\frac12(\pd_{x_j}+i\pd_{y_j}),
\endaligned
$$
where 
$\Phi=\Phi(x_1,\hdots,x_k,y_1,\hdots,y_k)$. 
If $F:\C^k\rightarrow\C$, then we frequently (and somewhat imprecisely) write $\pd_{\bar z}F, \pd_zF$ for $\pd_{\bar z}\widetilde F, \pd_z\widetilde F$.

\subsection{Single vs. double indexation of vectors}

Fix $m,n\in\N$. We will deal with complex vectors with $m,n$ or $mn$ components. 
Typically, a position of the index
\begin{itemize}
\item
 in the {\it subscript} \hskip 10.4pt indicates that it lies in the set $\{1,\hdots,n\}$;
\item
 in the {\it superscript} indicates that it lies in the set $\{1,\hdots,m\}$.
\end{itemize}

Accordingly, we will denote by $\C^m$ the usual Cartesian product of $m$ copies of $\C$, while $\C_n$ will be formally exactly the same thing, only with the number of copies ($n$) emphasized by the subscript. Thus we write
$$
\omega=(\omega^\alpha)^{\alpha=1,\hdots,m} \in\C^m.
$$
Along these lines, $\C_n^m$ is the set of all $mn$-tuples of complex numbers with double indices $m$ and $n$, that is
\begin{equation}
\label{{eq: xi components}}
\xi=\big(\xi_h^\alpha)_{h=1,\hdots,n}^{\alpha=1,\hdots,m}\in\C^m_n.
\end{equation}
Here $\omega^\alpha$ and $\xi_h^\alpha$ are complex numbers.

For $\xi\in\C^m_n$ we use the following notation:  
\begin{equation}
\label{eq: xi up-down}
\begin{aligned}
\xi^\alpha &:=\big(\xi_1^\alpha,\hdots,\xi_n^\alpha\big)\in\C_n,\\
\xi_h &:=\big(\xi_h^1,\hdots,\xi_h^m\big)^T\in\C^m. 
\end{aligned}
\end{equation}

\subsection{Inner products and norms}
For $\xi,\eta\in\C^m_n$, denote
$$
\sk{\xi}{\eta}
=\sk{\xi}{\eta}_{\C^m_n}
=\sum_{h=1}^n\sk{\xi_h}{\eta_h}_{\C^m}
=\sum_{\alpha=1}^m\sk{\xi^\alpha}{\eta^\alpha}_{\C_n}
=\sum_{h=1}^n\sum_{\alpha=1}^m\xi_h^\alpha\bar\eta_h^\alpha.
$$
Often we will omit the indices at inner products, assuming that the ambient spaces of the vectors involved indicate which inner product is intended. 

By $|\omega|, |\xi|, |\xi^\alpha|, |\xi_h|$ we understand the Euclidean norms of the vectors in question. For example,
$$
|\xi| 
= \sqrt{\sum_{h=1}^n\sum_{\alpha=1}^m |\xi_h^\alpha|^2}
= \sqrt{\sk{\xi}{\xi}}\,.
$$

For $\omega\in\C^m$ we define the map
$[\,\cdot\, ,\omega]\colon \C^m_n \rightarrow\C_n$ by
\begin{equation}
\label{eq: quasi-inner-product}
[\xi,\omega]=(\sk{\xi_1}{\omega},\hdots,\sk{\xi_n}{\omega}).
\end{equation}
Let us emphasize that $[\cdot,\omega]$ acts componentwise on the lower indices, meaning that
\begin{equation*}
[\xi,\omega]_h=\sk{\xi_h}{\omega}.
\end{equation*}
It is straightforward to verify that
\begin{equation}
\label{eq: quasi inner product - up summation}
[\xi,\omega]
=
\sum_{\alpha=1}^m
\overline{\omega}^\alpha
\xi^\alpha.
\end{equation}

With this notation, one easily confirms that the bracket 
$$
[\cdot,\cdot]: \C^m_n\times\C^m\rightarrow\C_n
$$ 
satisfies the following variant of the Cauchy--Schwarz--Bunyakovsky inequality:
\begin{equation}
\label{eq: CSB}
\left|[\xi,\omega]\right|\leq|\xi||\omega|.
\end{equation}

We draw the reader's attention to the difference between the scalar case ($m=1$) and the vector case ($m>1$). When $m=1$, we have $[\xi,\omega]=\overline{\omega}\xi$. We use square brackets to emphasize that when $m>1$, the map $[\cdot,\cdot]$ maps into vectors rather than scalars; hence, it is not a ``true'' inner product, but rather its $\C_n$-valued counterpart.

\subsection{Spheres in real and complex spaces}
\label{s: spheres}

Choose $m,n\in\N$. We denote by:
\begin{itemize}
\item
$S^{n-1}$ the unit sphere in $\R^n$;
\item
$\Sigma^{m}$ the unit sphere in $\C^m$, that is,
$$
\aligned
\Sigma^{m}
& :=\mn{\omega\in\C^m}{|\omega|=1}\\
& \ = \mn{(r_1e^{i\f_1},\hdots,r_me^{i\f_m})}{r_l,\f_l\in\R\text{ and }r_1^2+\hdots+r_m^2=1}
\,.
\endaligned
$$
\item
$\Sigma_n^{m}$ the unit sphere in $\C_n^m$, that is,
$
\Sigma_n^{m}
=\mn{\xi\in\C_n^m}{|\xi|=1}.
$
\end{itemize}

\subsection{Matrices and tensors}
\label{eq: matrices, tensors}
In the {\it scalar} case, the coefficients of a matrix $A$ as in Section \ref{s: scalar p-ell} are typically denoted by $a_{jk}$ with $j,k\in\{1,\hdots,n\}$.

In the {\it vector} case, we consider functions 
\begin{equation}
\label{eq: ajkab}
a^{\alpha\beta}_{jk} \colon \Omega \rightarrow \C,
\end{equation}
where $j,k \in \{1,\hdots,n\}$ and $\alpha,\beta \in \{1,\hdots,m\}$. For $\xi \in \C^m_n$, we define the action of $A(x)$ by
\begin{equation}
    \label{eq: Axxi}
(A(x)\xi)_j^\alpha := \sum_{\beta=1}^{m}\sum_{k=1}^{n} a_{jk}^{\alpha\beta}(x)\xi_k^\beta.
\end{equation}
The convention for index placement adopted here allows a natural interpretation of $A(x)$ as a matrix whose entries are themselves matrices. This viewpoint gives rise to two distinct realizations: the \emph{up} representation
\begin{equation}
\label{eq: Aup}
A^{\mathrm{up}} \colon \Omega \to (\C^{n\times n})^{m\times m}
\end{equation}
and
\begin{equation}
\label{eq: Adown}
A_{\mathrm{down}} \colon \Omega \to (\C^{m\times m})^{n\times n},
\end{equation}
which we now describe. 

Given $\alpha,\beta\in \{1,\dots,m\}$ and $j,k\in\{1,\dots,n\}$, we define 
$$
A^{\alpha,\beta}:\Omega\rightarrow\C^{n\times n}\quad \textrm{and} \quad
A_{j,k}:\Omega\rightarrow\C^{m\times m}
$$
by 
\begin{equation}
\label{eq: A_up_down_by_entries}
A^{\alpha,\beta}
=
\left[a_{jk}^{\alpha\beta}\right]_{j,k=1,\hdots,n}
\hskip 30pt 
\text{and}
\hskip 30pt 
A_{j,k}
=
\left[a_{jk}^{\alpha\beta}\right]^{\alpha,\beta=1,\hdots,m}.
\end{equation}
We then set
\begin{equation}
\label{eq: A up-down}
A^{\text{up}}=
\begin{bmatrix}
A^{1,1} & \hdots & A^{1,m}\\
\vdots & & \vdots\\
A^{m,1} & \hdots & A^{m,m}
\end{bmatrix}
\hskip 20pt 
\text{and}
\hskip 20pt 
A_{\text{down}}=
\begin{bmatrix}
A_{1,1} & \hdots & A_{1,n}\\
\vdots & & \vdots\\
A_{n,1} & \hdots & A_{n,n}
\end{bmatrix}.
\end{equation}
By using the notation \eqref{eq: xi up-down}, we have
\begin{equation}
\label{eq: Ax}
(A\xi)^\alpha=\sum_{\beta=1}^{m}A^{\alpha,\beta}\xi^\beta
\hskip 30pt 
\text{and}
\hskip 30pt 
(A\xi)_j=\sum_{k=1}^{n}A_{j,k}\xi_k.
\end{equation}
This may be written in a more compact manner as
\begin{equation}
\label{eq: Axi blockwise}
A\xi=
\begin{bmatrix}
A^{1,1} & \hdots & A^{1,m}\\
\vdots & & \vdots\\
A^{m,1} & \hdots & A^{m,m}
\end{bmatrix}
\begin{bmatrix}
\xi^{1}\\
\vdots \\
\xi^{m}
\end{bmatrix}
\hskip 20pt 
\text{or}
\hskip 20pt 
A\xi=
\begin{bmatrix}
A_{1,1} & \hdots & A_{1,n}\\
\vdots & & \vdots\\
A_{n,1} & \hdots & A_{n,n} 
\end{bmatrix}
\begin{bmatrix}
\xi_{1}\\
\vdots \\
\xi_{n}
\end{bmatrix}
.
\end{equation}
An illustration of the interplay between $A^{\text{up}}$ and $A_{\text{down}}$ is displayed in Remark \ref{r: up-to-down colors}.

Such an $A$ is called a {\it tensor} (acting on $\C_n^m$). 
Clearly, when $m=1$, we have just a single matrix (-valued function taking values) in $\C^{n\times n}$.

\begin{remark}
For the sake of consistency with \eqref{eq: A up-down}, we could also write
$$
\aligned
\xi^{\text{up}}&=\big(\xi^1,\hdots,\xi^m\big)^T\\ \xi_{\text{down}}&=\big(\xi_1,\hdots,\xi_n\big)^T.
\endaligned
$$
Then
$$
A\xi\equiv A^{\mathrm{ up}}\,\xi^{\mathrm{up}}\equiv A_{\mathrm{ down}}\,\xi_{\mathrm{down}}.
$$  
Yet in order to foster readability of the paper, in the sense of avoiding too many indices, we will mostly omit writing ``up'' and ``down'', unless the context will require it.
\end{remark}

\begin{remark}
Observe that for $m>1$, the purpose of working with tensors becomes visible when proceeding to the associated {\it operators}, cf. Section \ref{s: operator}. Namely, in the context of systems, we use tensors $A$ to build {\it vector}-valued operators $L_A$.
Thus a tensor $A$ on $\C_n^m$ gives rise to an operator $L_A$ acting on $\C^m$-valued functions.
\end{remark}

\begin{remark}
\label{r: up-to-down colors}
In general, the block-matrix expansions of $A^{\mathrm{up}}$ and $A_{\mathrm{down}}$ yield distinct matrices, since the spot occupied by $a_{jk}^{\alpha\beta}$ in one representation is taken by $a_{\alpha\beta}^{jk}$ in the other. Indeed, consider the (simplest ``nonscalar'') case $m=n=2$. Then
$$
\hskip 9.2pt
A^{\text{up}}
=
\left[
\begin{array}{cc}
A^{1,1} & 
A^{1,2}\\
A^{2,1} & 
A^{2,2}
\end{array}
\right]
=
\left[
\begin{array}{cccc}
\cellcolor{LimeGreen} a_{11}^{11} & \cellcolor{LimeGreen} a_{12}^{11} & 
	\cellcolor{Goldenrod} a_{11}^{12} & \cellcolor{Goldenrod} a_{12}^{12} \\
\cellcolor{LimeGreen} a_{21}^{11} & \cellcolor{LimeGreen} a_{22}^{11} & 
	\cellcolor{Goldenrod} a_{21}^{12} & \cellcolor{Goldenrod} a_{22}^{12} \\
\cellcolor{Lavender} a_{11}^{21} & \cellcolor{Lavender} a_{12}^{21} & 
	\cellcolor{CornflowerBlue} a_{11}^{22} & \cellcolor{CornflowerBlue} a_{12}^{22} \\
\cellcolor{Lavender} a_{21}^{21} & \cellcolor{Lavender} a_{22}^{21} & 
	\cellcolor{CornflowerBlue} a_{21}^{22} & \cellcolor{CornflowerBlue} a_{22}^{22}
\end{array}
\right]\,\ 
$$
and
$$
A_{\text{down}}
=
\left[
\begin{array}{cc}
A_{1,1} & A_{1,2}\\
A_{2,1} & A_{2,2}
\end{array}
\right]
=
\left[
\begin{array}{cccc}
\cellcolor{LimeGreen} a^{11}_{11} & \cellcolor{Goldenrod} a^{12}_{11} & 
	\cellcolor{LimeGreen} a^{11}_{12} & \cellcolor{Goldenrod} a^{12}_{12} \\
\cellcolor{Lavender} a^{21}_{11} & \cellcolor{CornflowerBlue} a^{22}_{11} & 
	\cellcolor{Lavender} a^{21}_{12} & \cellcolor{CornflowerBlue} a^{22}_{12} \\
\cellcolor{LimeGreen} a^{11}_{21} & \cellcolor{Goldenrod} a^{12}_{21} & 
	\cellcolor{LimeGreen} a^{11}_{22} & \cellcolor{Goldenrod} a^{12}_{22} \\
\cellcolor{Lavender} a^{21}_{21} & \cellcolor{CornflowerBlue} a^{22}_{21} & 
	\cellcolor{Lavender} a^{21}_{22} & \cellcolor{CornflowerBlue} a^{22}_{22}
\end{array}
\right].
$$
Clearly, the two matrices on the right are not the same.

Instead, the appropriate way to view $A$ as a block matrix is through its action on vectors, as shown in~\eqref{eq: Axi blockwise}.
\end{remark}

\subsection{More about tensors}
Here we suggest a possible (and slightly different) variant of the constructions above, namely, the one emphasizing connections with standard algebraic objects -- tensor products. 
Note that in this setup, $\xi^\alpha$ are thought of as rows, and not as columns anymore. Thus, an interpretation of $A^{\alpha,\beta}\xi^\beta$ in \eqref{eq: Ax} goes through perceiving this as an action of an endomorphism $A^{\alpha,\beta}$ on an element $\xi^\beta$ from a $n$-dimensional space, rather than an action of a matrix on a vector, cf. Remark \ref{eq: dif_up}.

Fix $m, n \in \N_+$. 
We denote by $\C^m$ the standard Cartesian product of $m$ copies of $\C$ (column vectors), and by $\C_n$ the dual space $(\C^n)^*$ which we identify with row vectors. We consider the space of rank-$2$ tensors $\C^m_n := \C^m \otimes \C_n$. We identify this space with $\C^{m\times n}$, the space of complex matrices with $m$ rows and $n$ columns. 
The identification is canonical: we may view $\xi \in \C^m_n$ as a linear map $\xi \colon \C^n \to \C^m$. By fixing the canonical bases in $\C^n$ and $\C^m$, we identify $\xi$ with its representing $m \times n$ matrix. Regarding the index positioning, 
\begin{itemize}
\item
{\it subscript} denotes a covariant index taking values in $\{1,\hdots,n\}$;
\item
{\it superscript} denotes a contravariant index taking values in $\{1,\hdots,m\}$.
\end{itemize}
As an example, if $\omega\in\C^m$ and 
$\zeta\in\C_n$ then 
$\omega\otimes \zeta\in \C^m_n$ 
and
\begin{equation*}
(\omega\otimes\zeta)_h^\alpha:=\omega^\alpha\zeta_h.
\end{equation*}

We point out a minor deviation from the notation employed in the scalar case (Section~\ref{s: scalar p-ell}); there, matrices $A \in \C^{n\times n}$ were represented with two covariant indices, treating row and column vectors indistinguishably.

\begin{remark}
For $\xi\in\C^m_n\cong \C^{m\times n}$, we define the transpose $\xi^T$ and adjoint $\xi^*$ in $\C^n_m\cong \C^{n\times m}$ by swapping indices as customary:
$$
(\xi^T)^h_\alpha=\xi^\alpha_h, \qquad (\xi^*)^h_\alpha=\overline{(\xi^T)^h_\alpha}=\overline{\xi^\alpha_h}.
$$
In this notation, the bracket operation corresponds to the conjugate transpose of the matrix-vector product:
$$
[\xi,\omega] = (\overline{\xi^*(\omega)})^T.
$$
Furthermore, the standard adjoint identity holds: for any $\zeta \in \C^n$ (column),
$$
\sk{\xi(\zeta)}{\omega}_{\C^m} = \sk{\zeta}{\xi^*(\omega)}_{\C^n} = \sk{\zeta^T}{\overline{[\xi,\omega]}}_{\C_n}.
$$
Note also that $[\cdot,\omega]=\sk{\cdot}{\omega} \otimes \mathrm{Id}_{\C_n}$.
\end{remark}

\subsubsection{Rank-4 tensors}
In the {\it vector} case, the coefficient matrix is replaced by a rank-$4$ tensor-valued function. For a given $x\in\Omega$, $A(x) \colon \C^m_n \rightarrow \C^m_n$ is a linear map; canonically, it would therefore be an element of $\C^m_n \otimes (\C^m_n)^*$. 
However, since we prefer to work with objects whose subscripts and superscripts always run in the same fixed sets, we consider $A(x)$ as an element of $\C^m_n \otimes \C^m_n$ instead. 
That means that we now consider functions \eqref{eq: ajkab}.

The definition \eqref{eq: Axxi} now relies on the identification of $\C^m_n\otimes \C^m_n$ with $\mathrm{End}(\C^m_n)$. For example, for a decomposable rank-$4$ tensor $B=\xi\otimes\eta$, with $\xi,\eta\in\C^m_n$, the corresponding action is given by
$$
B(\gamma)=(\xi\otimes \eta)(\gamma)=\sk{\gamma}{\overline{\eta}}_{\C^m_n}\xi,\quad \forall\gamma\in\C^m_n.
$$

Instead of \eqref{eq: Aup} and \eqref{eq: Adown} we now have 
$$
\aligned
A^{\mathrm{up}} & : \Omega \to \mathrm{End}(\C_n)^{m\times m}\\
A_{\mathrm{down}} & : \Omega \to \mathrm{End}(\C^m)^{n\times n},
\endaligned
$$
as follows.
Given contravariant indices $\alpha,\beta \in \{1,\dots,m\}$ and covariant indices $j,k \in \{1,\dots,n\}$, define the operators
$$
A^{\alpha,\beta} : \Omega \to \mathrm{End}(\C_n)
\qquad \text{and} \qquad
A_{j,k} : \Omega \to \mathrm{End}(\C^m),
$$
by setting
$$
(A^{\alpha,\beta}\zeta)_j := \sum_{k=1}^n a^{\alpha\beta}_{jk}\,\zeta_k,
\qquad
(A_{j,k}u)^\alpha := \sum_{\beta=1}^m a^{\alpha\beta}_{jk}\,u^\beta.
$$
When convenient, we identify $A_{j,k}$ and $A^{\alpha,\beta}$ with their matrix representations with respect to the canonical bases of $\C^m$ and $\C_n$, just as in \eqref{eq: A_up_down_by_entries}.
Likewise, \eqref{eq: A up-down}, \eqref{eq: Ax}, \eqref{eq: Axi blockwise} and Remark \ref{r: up-to-down colors}
continue to apply as before.

\begin{remark}
\label{eq: dif_up}
In the {\it up} representation, the vectors $\xi^\alpha \in \C_n$ are viewed as the
entries of a block column vector, and multiplication is understood as the
natural action of the operators $A^{\alpha,\beta} \in \mathrm{End}(\C_n)$ on
$\xi^\beta$.
\end{remark}

\begin{remark}
\label{r: up/down}
The {\it up} and {\it down} representations correspond to two natural ways of viewing
$A\in\C^m_n\otimes\C^m_n$ as an element of $\mathrm{End}(\C^m_n)=\textrm{End}(\C^m\otimes\C_n)=\mathrm{End}(\C^m)\otimes \mathrm{End}(\C_n)$. Let $\{E^{\alpha,\beta}\}_{\alpha,\beta=1}^m$ and $\{F_{j,k}\}_{j,k=1}^n$ denote the canonical bases of $\mathrm{End}(\C^m)$ and $\mathrm{End}(\C_n)$, respectively, defined by
$$
E^{\alpha,\beta}e_\gamma = \delta_{\beta\gamma}e_\alpha, \qquad F_{j,k}f_h = \delta_{kh}f_j,
$$
where $\{e_\alpha\}_{\alpha=1}^m$ and $\{f_j\}_{j=1}^n$ are the canonical bases of $\C^m$ and $\C_n$. Then we have
\begin{equation*}
\begin{aligned}
A &= \sum_{\alpha,\beta=1}^m E^{\alpha,\beta} \otimes A^{\alpha,\beta} \quad (\text{{\it up} representation}), \\
A &= \sum_{k,j=1}^n A_{k,j} \otimes F_{k,j} \quad (\text{{\it down} representation}).
\end{aligned}
\end{equation*}
\end{remark}

\subsection{Sobolev spaces of vector-valued functions}
Throughout this paper, unless otherwise specified, all integrals are with respect to the Lebesgue measure on $\R^n$.

We define $L^2(\Omega;\C^m)$ to be the space of all functions $u=(u^1,\hdots,u^m)^T:\Omega\rightarrow\C^m$, such that $u^\alpha\in L^2(\Omega)$. 
That is, $L^2(\Omega;\C^m)\cong L^2(\Omega)^m:=L^2(\Omega)\times\hdots\times L^2(\Omega)$, with $m$ copies of $L^2(\Omega)$ on the right.
We define the norm
$$
\nor{u}_{2}:=
\sqrt{
\int_{\Omega}|u|^2
}.
$$
Explicitly,
$$
\nor{u}_{2}=\sqrt{\nor{u^1}_{2}^2+\hdots+\nor{u^m}_{2}^2}.
$$
Similarly, for matrix-valued functions $\xi:\Omega\rightarrow\C_n^m$, we define
$$
\nor{\xi}_{2}
:=
\sqrt{\int_\Omega|\xi|^2}= \sqrt{\int_\Omega\sum^{m}_{\alpha=1}\sum^n_{h=1}|\xi_h^\alpha|^2}\,.
$$

We denote by $H^1(\Omega;\C^{m})$ the subspace of $L^2(\Omega;\C^m)$ consisting of all vector-valued functions whose components belong to the Sobolev space $H^1(\Omega)=W^{1,2}(\Omega;\C)$.
That is, 
$$
H^1(\Omega;\C^m)\cong H^1(\Omega)^m,
$$
endowed with the norm
$$
\norm{u}{H^1(\Omega;\C^m)}:=\sqrt{\sum_{\alpha=1}^m\nor{u^\alpha}^2_{H^1(\Omega)}}=\sqrt{\nor{u}_2^2+\nor{\nabla u}_2^2}.
$$
Here
$$
\nabla u=(\nabla u^1,\dots,\nabla u^m)^T:\Omega\rightarrow \C^m_n
$$
denotes the {\it Jacobian matrix} of $u$. Thus, its rows correspond to the gradients of the components $u^\alpha$, while its columns correspond to the partial derivatives $\partial_h u$ of the vector $u$ with respect to the spatial variables $x_h$. In terms of components,
$$
(\nabla u)^\alpha_h = \partial_h u^\alpha.
$$
On scalar functions we have
$\nabla=(\pd_{x_1},\hdots,\pd_{x_n})$, 
with $(x_1,\hdots,x_n)$ being the Cartesian coordinates in $\Omega\subset\R^n$.

\subsection{Sectors}
\label{sec}
Given $\theta\in (0,\pi)$, define the (open) {\it sector of angle $\theta$} in $\C$ by
\begin{equation*}
\label{autobus karlovac-posedarje}
\bS_{\theta}=\{z\in\C\setminus\{0\}\,;\ |\arg z|<\theta\}.
\end{equation*}
Also set $\bS_{0}=(0,\infty)$.

\section{Operators associated with systems of elliptic PDEs in divergence form}

For a general reference on elliptic systems see Chen and Wu \cite{CW98} or Giaquinta and Martinazzi \cite{GM12}. Let $A$ be a rank-$4$ tensor-valued function on $\Omega$, cf. Section~\ref{eq: matrices, tensors}. We aim to define the associated {\it divergence-form elliptic operator} $L_A$, which constitutes the counterpart of the scalar operator $L_Au=-\div(A\nabla u)$. We present this construction rigorously below.

\subsection{Ellipticity condition (or the Legendre condition) for tensors}
\label{s: Legendre}
Suppose that $A=\mn{A^{\alpha,\beta}}{\alpha,\beta=1,\hdots,m}$ is a family of complex $n\times n$ matrix functions on an open set $\Omega\subset\R^n$. 
We impose on $A=(a^{\alpha\beta}_{kj})^{\alpha,\beta=1,\dots,m}_{k,j=1,\dots,n}$  
the {\it ellipticity} (or {\it accretivity}), also known as the {\it (strong) Legendre condition}.
Specifically, we assume that there exists $\lambda=\lambda(A)>0$ such that for almost every $x\in\Omega$:
$$
\Re\sk{A(x)\xi}{\xi}\geq \lambda|\xi|^2,\quad \forall \xi\in\C^m_n.
$$
Here
$$
\Re\sk{A(x)\xi}{\xi}
=\Re\sum^m_{\alpha,\beta=1}\Sk{A^{\alpha,\beta}(x)\xi^\beta}{\xi^\alpha}
=\Re\sum^m_{\alpha,\beta=1}\sum^n_{k,j=1}a^{\alpha\beta}_{kj}(x)\xi^{\beta}_{j}\,\overline{\xi^{\alpha}_{k}}.
$$
We denote the family of all elliptic tensors on $\Omega$ acting on $\C^m_n$ by $\cA^m(\Omega)$, and we recall that the scalar case corresponds to taking $m=1$.

Also, we require $A$ to be {\it bounded}, meaning that
there exists $\Lambda=\Lambda(A)>0$ such that for almost every $x\in\Omega$,
\begin{equation}
\label{eq: A bddness}
|\sk{A(x)\xi}{\eta}|\leq\Lambda|\xi||\eta|
\hskip 30pt
\forall \xi,\eta\in\C_n^m\,.
\end{equation}

\begin{remark}
There exist other classical ellipticity conditions for systems, such as the {\it integral ellipticity condition} and the {\it (strong) Legendre–Hadamard condition}  (see e.g. \cite[Section 8.1]{CW98}, \cite[p.369]{DiLiPi20}). However, we focus on the Legendre ellipticity condition in this paper, as it most naturally extends the concept of ellipticity for scalar operators.
\end{remark}

\subsection{The sesquilinear form associated with $A$}
\label{s: the form}

We denote by $H^1_{0}(\Omega)$ the closure of $C_c^\infty(\Omega)$ in the Sobolev space $H^1(\Omega)=W^{1,2}(\Omega)$. 
Let $\mathscr{V}$ be a closed linear subspace of $H^1(\Omega)$ satisfying
\begin{equation}
\label{eq: V}
H^1_{0}(\Omega)
\subseteq\mathscr{V}
\subseteq
H^1(\Omega).
\end{equation}
Recall that if $\Omega = \R^d$, then $H^1_{0}(\R^{d})=H^1(\R^{d})$ (see, e.g., \cite[Corollary~3.19]{Adams}).

\subsubsection{Special classes of boundary conditions}
\label{s: bdry conds}

We now describe specific classes of domains $\mathscr{V}$ that, in addition to satisfying \eqref{eq: V}, fulfill certain invariance properties which will be assumed in the sequel.

We say that the subspace $\mathscr{V}\subset H^1(\Omega)$ is {\it invariant} under: 
\begin{itemize}
\item the function $p:\C\rightarrow\C$, if $u\in\mathscr{V} \implies p\circ u\in\mathscr{V}$;
\item the family $\mathscr{P}$ of functions $\C\rightarrow\C$, if it is invariant under every $p\in\mathscr{P}$.
\end{itemize}

Consider the projection function $P:\C\rightarrow\C$ defined by 
$$
P(\zeta)= 
\begin{cases} 
\zeta & \text{if }\  |\zeta| \leq 1, \\ 
\zeta/|\zeta| & \text{if }\  |\zeta| > 1. 
\end{cases}
$$
Throughout this paper, we assume that $\mathscr{V}$ satisfies:
\begin{equation}
\label{eq: invariance}
\mathscr{V} \text{ is invariant under the function } P.
\end{equation}
It is well known (see Ouhabaz \cite[Proposition 4.11]{O}) that condition \eqref{eq: invariance} holds for the following subspaces:
\begin{enumerate}[(a)]

\item 
$\mathscr{V}=H^1_{0}(\Omega)$; 

\item 
$\mathscr{V}=H^1(\Omega)$; 

\item 
$\mathscr{V}$ is the closure in $H^1(\Omega)$ of $\{u|_\Omega : u\in C_c^\infty(\R^{d}\setminus\Gamma)\}$, where $\Gamma$ is a closed subset of $\partial\Omega$. 

\end{enumerate}

In the context of the associated operators (see Section \ref{s: operator}), the above choices of $\oV$ correspond to (a) {\it Dirichlet}, (b) {\it Neumann} or (c) {\it mixed boundary conditions}.

\subsubsection{The form}
\label{s: form def}
Let $\oV^m$ denote the space of vector-valued functions whose components belong to $\mathscr{V}$:
$$
\oV^m = \left\{ u=(u^1,\dots,u^m)^T\in H^1(\Omega;\C^m):\  u^\alpha \in \mathscr{V} \text{ for all } \alpha=1,\dots,m \right\}.
$$
We define the sesquilinear form $\mathfrak{a} = \mathfrak{a}_{A,\oV^m}: \oV^m \times \oV^m \rightarrow \C$ by
\begin{equation*}
\label{eq: form_def}
\mathfrak{a}(u,v):=\int_\Omega \sk{A\nabla u}{\nabla v}_{\C^m_n}.
\end{equation*}

By adopting the terminology from Ouhabaz \cite[Definitions 1.4 and 1.7]{O} and adapting the proofs from Section 4.1 therein, we obtain that if $A$ satisfies the Legendre ellipticity condition (see Section~\ref{s: Legendre}), then the form $\a$ is
densely defined, sectorial
and closed. The norms $\nor{\cdot}_{\am}$ and $\nor{\cdot}_{H^1(\Omega;\C^m)}$ are equivalent on $D(\a)=\oV^m$.

\subsection{The operator associated with $A$}
\label{s: operator}

We define $L=L_{A}=L_{A,\oV^m}$ to be the operator on $L^2(\Omega;\C^m)$, 
associated with $\a=\a_{A,\oV^m}$. See \cite[Section 1.2.3]{O} for this construction. 
So we have
\begin{equation}
\label{eq: L def}
\int_\Omega
\sk{Lu}{v}_{\C^m}
=
\int_\Omega
\sk{A\nabla u}{\nabla v}_{\C^m_n}
\hskip 20pt
\forall u\in D(L),\ \forall
v\in D(\a).
\end{equation}

That is,
\begin{equation*}
\aligned
\sk{Lu}{v}_{L^2(\Omega;\C^m)}
=
\sum^m_{\alpha,\beta=1} 
\sk{A^{\alpha,\beta}\nabla u^\beta}{\nabla v^\alpha}_{L^2(\Omega;\C_n)}
=\sum^m_{\alpha,\beta=1}\sum^n_{k,j=1}\sk{a^{\alpha\beta}_{kj}\partial_j u^\beta}{\partial_kv^\alpha}_{L^2(\Omega)}.
\endaligned
\end{equation*}
Thus, formally (since the individual tensors $A^{\alpha\beta}$ need not be elliptic), one may write $Lu:\Omega\rightarrow\C^m$ with
$$
Lu=
\begin{bmatrix}
L_{A^{1,1}} & \hdots & L_{A^{1,m}}\\
\vdots & & \vdots\\
L_{A^{m,1}} & \hdots & L_{A^{m,m}}
\end{bmatrix}
\begin{bmatrix}
u^1\\
\vdots\\
u^m
\end{bmatrix}.
$$

When $\oV$ corresponds to one of the cases (a)–(c) in Section~\ref{s: bdry conds}, we say that $L$ is subject to {\it Dirichlet}, {\it Neumann}, or {\it mixed boundary conditions}, respectively.

Sectorial forms are automatically accretive and continuous; see e.g. \cite[Proposition 1.8]{O}.
Therefore, by \cite[Proposition~1.27 and Theorem~1.54]{O}, the operator $-L$ generates on $L^2(\Omega;\C^m)$ a strongly continuous semigroup of operators
$$
T_{t}=T^{A}_{t}=T^{A,\oV^m}_{t}=
e^{-tL},\quad t>0,
$$
which is holomorphic and contractive in a sector of positive angle. Hence, for every $t>0$, $T_{t}$ maps  $L^{2}(\Omega;\C^m)$ into $D(L)\subseteq\oV^m$ 
\cite[Theorem II.4.6]{EN}, 
and since $\oV^m\subseteq H^1(\Omega;\C^m)$, 
the spatial gradient $\nabla  T_{t}f$ is
well defined. 

By \cite[p.~72]{Stein}, given $u\in L^{2}(\Omega;\C^m)$, we may redefine each $T_tu$ on a set of measure zero so that, for almost every $x\in\Omega$, the function $t \mapsto T_tu(x)$ is real-analytic on $(0,\infty)$.


\section{Auxiliary Euclidean operators}
\label{s: auxops}

We introduce a couple of operators on $\C^m_n$ that will be essential as we proceed.

\subsection{Linear algebra preliminaries}
First we formulate a few elementary results.
Recall that we identify elements in $\R^\ell$ with column vectors.

\begin{proposition}
\label{p: real symmetric}
Let $\ell\geq 2$ and $\mathsf{x}_1,\hdots,\mathsf{x}_n,\mathsf{y}_1,\hdots,\mathsf{y}_n\in\R^{\ell}$.
Consider the associated $\ell\times\ell$ real symmetric matrix $M$ defined by:
\begin{equation}
\label{eq: M}
M=\frac12\sum_{j=1}^n\left(\mathsf{x}_j\mathsf{y}_j^T+\mathsf{y}_j\mathsf{x}_j^T\right)
=\frac12\sum_{j=1}^n\left(\mathsf{x}_j\otimes\mathsf{y}^T_j+\mathsf{y}_j\otimes\mathsf{x}^T_j\right).
\end{equation}
Then the following hold:
\begin{enumerate}[(i)]

\item For every $\sigma\in S^{\ell-1}$, 
\begin{equation}
\label{eq: M quad form}
\sum_{j=1}^n
\sk{{\mathsf{x}_j}}{\sigma}\sk{{\mathsf{y}_j}}{\sigma}
=\sk{M\sigma}{\sigma}_{\R^{\ell}}.
\end{equation}

\item 
If $\mu_1\leq\cdots\leq\mu_{\ell}$ are the eigenvalues of $M$, then
\begin{equation}
\label{eq: mu min/max}
\aligned
\min_{\sigma\in S^{\ell-1}} \sum_{j=1}^n
\sk{{\mathsf{x}_j}}{\sigma}\sk{{\mathsf{y}_j}}{\sigma}
& = \mu_1\\
\max_{\sigma\in S^{\ell-1}} \sum_{j=1}^n
\sk{{\mathsf{x}_j}}{\sigma}\sk{{\mathsf{y}_j}}{\sigma}
& = \mu_{\ell}.
\endaligned
\end{equation}

\item 
When $n=1$, we have
\begin{equation}
\label{eq: mu min/max n=1}
\aligned
\mu_1 & =\frac12(\sk{\mathsf{x}}{\mathsf{y}}-|\mathsf{x}||\mathsf{y}|),\\
\mu_{\ell} & =\frac12(\sk{\mathsf{x}}{\mathsf{y}}+|\mathsf{x}||\mathsf{y}|).
\endaligned
\end{equation}
Consequently, $\mu_1\leq0\leq\mu_{\ell}$.
If, in addition, $\ell>2$, then $\mu_2=\cdots=\mu_{\ell-1}=0$.

\item
Furthermore,
\begin{equation}
\label{trace formula}
\sum_{j=1}^n
\sk{{\mathsf{x}_j}}{{\mathsf{y}_j}}
=\mu_1+\hdots+\mu_\ell.
\end{equation}

\end{enumerate}

\end{proposition}

\begin{proof} 
By definition, for any ${\mathsf x},{\mathsf y},\sigma\in\R^\ell$ we have 
\begin{equation}
\label{eq: Kronecker2}
({\mathsf x}{\mathsf y}^T)\sigma=({\mathsf x}\otimes{\mathsf y}^T)\sigma=\sk{\sigma}{{\mathsf y}}{\mathsf x}.
\end{equation}
Therefore, 
$$
\sk{({\mathsf x}{\mathsf y}^T)\sigma}{\sigma}
=\sk{({\mathsf x}\otimes{\mathsf y}^T)\sigma}{\sigma}
=\sk{{\mathsf y}}{\sigma}\sk{{\mathsf x}}{\sigma},
$$
from which we immediately deduce \eqref{eq: M quad form}. 

By the spectral theorem for real symmetric matrices, all the eigenvalues of $M$ are real. 
Moreover, there exists a real orthogonal matrix $Q$ such that $M=QDQ^T$, where $D$ is a diagonal matrix with eigenvalues of $M$ on the diagonal. Together with \eqref{eq: M quad form}, this yields \eqref{eq: mu min/max}.
Alternatively, we could have recalled the variational characterization of eigenvalues (Rayleigh quotient).

From \eqref{eq: Kronecker2} we infer that ${\mathsf x}{\mathsf y}^T$ is a rank 1 operator (or matrix) with ${\mathsf x}$ being the eigenvector for the only possibly nonzero eigenvalue $\sk{{\mathsf x}}{{\mathsf y}}$.
Likewise, ${\mathsf x}{\mathsf y}^T+{\mathsf y}{\mathsf x}^T$ is an operator of rank 2 (if ${\mathsf x},{\mathsf y}$ are linearly independent) and one checks that ${\mathsf x}/|{\mathsf x}|\pm{\mathsf y}/|{\mathsf y}|$ are the eigenvectors for the eigenvalues $\sk{{\mathsf x}}{{\mathsf y}}\pm|{\mathsf x}||{\mathsf y}|$.  
This gives \eqref{eq: mu min/max n=1}.

Let us present a few more details in the last part. 
If $n=1$ then 
$$
M=\frac12({\mathsf x}\otimes{\mathsf y}^T+{\mathsf y}\otimes{\mathsf x}^T ),
$$
for some ${\mathsf x},{\mathsf y}\in\R^\ell$. 
The spaces $\textrm{span}\{{\mathsf x},{\mathsf y}\}$ and $\textrm{span}\{{\mathsf x},{\mathsf y}\}^\perp$ are invariant under $M$. 
Clearly, $M$ vanishes on $\textrm{span}\{{\mathsf x},{\mathsf y}\}^\perp$. 
Restricted to the (at most) two-dimensional subspace $\textrm{span}\{{\mathsf x},{\mathsf y}\}$, the operator $M$ acts on a generic vector $v = a{\mathsf x}+b{\mathsf y}$, $a,b\in\R$, as:
$$
M(a{\mathsf x}+b{\mathsf y}) = \frac{1}{2} \left( a\sk{\mathsf x}{\mathsf y}{\mathsf x} + a|\mathsf x|^2{\mathsf y} + b|\mathsf y|^2{\mathsf x} + b\sk{\mathsf x}{\mathsf y}{\mathsf y} \right).
$$
Suppose that ${\mathsf x},{\mathsf y}$ are linearly independent. 
In the basis $\{{\mathsf x},{\mathsf y}\}$, this action is represented by the matrix multiplication:
$$
\begin{bmatrix}
a\\
b
\end{bmatrix}
\mapsto
 \frac12
\begin{bmatrix}
  \sk{\mathsf{x}}{\mathsf{y}} & |\mathsf{y}|^2 \\
  |\mathsf{x}|^2 & \sk{\mathsf{x}}{\mathsf{y}}
\end{bmatrix}.\begin{bmatrix}
a\\
b
\end{bmatrix}.
$$
The eigenvalues of this $2\times 2$ matrix are $(1/2)(\sk{{\mathsf x}}{{\mathsf y}}\pm |{\mathsf x}||{\mathsf y}|)$, which yields \eqref{eq: mu min/max n=1}.

At last we prove \eqref{trace formula}. 
By a standard result from linear algebra, 
$
\mu_1+\hdots+\mu_{\ell}
=\text{tr}\,M
$.
Since for $M$ as in \eqref{eq: M} one has
$\text{tr}\,M=\sum_{j=1}^n
\sk{\mathsf{x}_j}{\mathsf{y}_j}$, 
this finishes the proof.
\end{proof}

\begin{remark}
\label{r: range of M}
The definition \eqref{eq: M} implies,
for $\sigma\in\R^{\ell}$,
\begin{equation*}
M\sigma=\frac12\sum_{j=1}^n
\left(
\sk{\sigma}{\mathsf{y}_j}\mathsf{x}_j+\sk{\sigma}{\mathsf{x}_j}\mathsf{y}_j
\right).
\end{equation*}
From 
here 
we see that 
the range of $M$ is at most $2n$-dimensional.  

Thus \framebox{if $\ell>2n$}, then the kernel of $M$ is nontrivial, in other words, $0$ is among the eigenvalues of $M$,
regardless of the concrete choice of $\mathsf{x}_j$'s and $\mathsf{y}_j$'s.
\end{remark}

\begin{corollary}
\label{c: bilateral est}
For any $m,n\in\N$ and $u_1,\hdots,u_n,v_1,\hdots,v_n\in\C^m$ we have
\begin{equation}
\label{eq: (u,v) bilateral est complex}
\min_{\omega\in\Sigma^m}
\sum_{j=1}^n
\Re\sk{u_j}{\omega}\Re\sk{v_j}{\omega}
\leq
\frac1{2m}
\sum_{j=1}^n\Re\sk{u_j}{v_j}
\leq
\max_{\omega\in\Sigma^{m}}
\sum_{j=1}^n
\Re\sk{u_j}{\omega}\Re\sk{v_j}{\omega},
\end{equation}
where all the inner products above are the standard Hermitian inner products in $\C^{m}$.
\end{corollary}

\begin{proof}
We use the identification operator $\cV:\C^m\rightarrow \R^{2m}$ introduced in Section~\ref{s: notation} and write $\mathsf{x}_j=\widetilde{u}_j$, $\mathsf{y}_j=\widetilde{v}_j$, $\sigma=\widetilde{\omega}$. Consequently, establishing \eqref{eq: (u,v) bilateral est complex} is equivalent to showing that for all $\mathsf{x}_1,\dots,\mathsf{x}_n$, $\mathsf{y}_1,\dots,\mathsf{y}_n\in\R^{2m}$ the following holds:
\begin{equation}
\label{eq: (u,v) bilateral est real}
\min_{\sigma\in S^{2m-1}}
\sum_{j=1}^n
\sk{\mathsf{x}_j}{\sigma}\sk{\mathsf{y}_j}{\sigma}
\leq
\frac1{2m}\sum_{j=1}^n
\sk{\mathsf{x}_j}{\mathsf{y}_j}
\leq
\max_{\sigma\in S^{2m-1}}
\sum_{j=1}^n
\sk{\mathsf{x}_j}{\sigma}\sk{\mathsf{y}_j}{\sigma},
\end{equation}
where all inner products above are the standard Euclidean inner products in $\R^{2m}$.

Now simply apply identities \eqref{eq: mu min/max} and \eqref{trace formula} from Proposition \ref{p: real symmetric}
with $\ell=2m$.
\end{proof}

\begin{remark}
\label{c: bilateral est complex}
The inequalities \eqref{eq: (u,v) bilateral est complex} and \eqref{eq: (u,v) bilateral est real} can be expressed in terms of the bracket  \eqref{eq: quasi-inner-product} as follows: 
\begin{enumerate}
\item For any $\xi,\eta\in\C_n^m$ we have
\begin{equation*}
\min_{\omega\in\Sigma^{m}}
\Sk{\Re[\xi,\omega]}{\Re[\eta,\omega]}_{\R_n}
\leq
\frac1{2m}
\,\Re\sk{\xi}{\eta}_{\C^m_n}
\leq
\max_{\omega\in\Sigma^{m}}
\Sk{\Re[\xi,\omega]}{\Re[\eta,\omega]}_{\R_n}.
\end{equation*}
\item For any $\mathsf{X},\mathsf{Y}\in\R^{2m}_n$ we have
\begin{equation*}
\min_{\sigma\in S^{2m-1}}
\Sk{[\mathsf{X},\sigma]}{[\mathsf{Y},\sigma]}_{\R_{n}}
\leq
\frac1{2m}
\,\sk{\mathsf{X}}{\mathsf{Y}}_{\R_n^{2m}}
\leq
\max_{\sigma\in S^{2m-1}}
\Sk{[\mathsf{X},\sigma]}{[\mathsf{Y},\sigma]}_{\R_{n}}.
\end{equation*}
\end{enumerate}
Thus we may call \eqref{eq: (u,v) bilateral est complex} and \eqref{eq: (u,v) bilateral est real} {\it double inner product estimates}.
\end{remark}


\subsection{Projections $\Phi_\omega$}

Given $\omega\in\Sigma^m$ we denote by $\Pi_{\widetilde{\omega}}=\widetilde{\omega}\otimes\widetilde{\omega}^T:\R^{2m}\rightarrow\R^{2m}$ the orthogonal projection onto the (real) linear span of $\wt\omega\in S^{2m-1}$, namely 
$$
\Pi_{\widetilde{\omega}}\,\mathsf{x}=\sk{\mathsf{x}}{\widetilde{\omega}}_{\R^{2m}}\, \widetilde{\omega},\quad \forall\mathsf{x}\in\R^{2m}.
$$
We denote by $\pi_\omega:\C^m\rightarrow\C^m$ the operator on $\C^m$ induced by $\Pi_{\widetilde{\omega}}$ via the isomorphism $\cV_m$. 
That is,
$$
\pi_\omega=\cV_m^{-1}\circ \Pi_{\wt\omega}\circ\cV_m.
$$
We further define 
$$
\Phi_\omega:\C^m_n\to\C^m_n
$$ 
as the direct sum of $n$ copies of $\pi_\omega$:
$$
\Phi_\omega
=\pi_\omega^{\oplus n}
:=\underbrace{\pi_\omega\oplus\hdots\oplus\pi_\omega}_{n}:\big(\C^m\big)^n\rightarrow\big(\C^m\big)^n,
$$
in the sense that
\begin{equation}
\label{eq: Phi_omega}
\Phi_\omega(\xi_1,\hdots,\xi_n)
=
(\pi_\omega \xi_1,\hdots,\pi_\omega \xi_n),
\end{equation}
where $\xi_h\in\C^m$ for all $h\in\{1,\hdots,n\}$.

Explicitly, for $u\in\C^m$, 
\begin{equation}
\label{eq: pi_omega}
\pi_\omega u
=
\Re\sk{u}{\omega}_{\C^m}\omega,
\end{equation}
while for $\xi\in\C^m_n$, 
$$
\left(\Phi_\omega \xi\right)_h^\alpha
=\omega^\alpha\Re\sk{\xi_h}{\omega}_{\C^m}
=\omega^\alpha\sum_{\beta=1}^m\Re\big(\xi_h^\beta\bar\omega^\beta\big).
$$
In particular,
\begin{equation}
\label{eq: Phi up/down}
\aligned
\left(\Phi_\omega \xi\right)^\alpha
& =
\omega^\alpha
\,\Re[\xi,\omega]\\
\left(\Phi_\omega \xi\right)_h
& =
\Re\sk{\xi_h}{\omega}_{\C^m}
\omega.
\endaligned
\end{equation}
Recalling the bracket notation from \eqref{eq: quasi-inner-product}, we obtain the compact tensor identity
\begin{equation*}
\label{Lacrimosa 5.16}
\Phi_\omega \xi=\omega\otimes\Re[\xi,\omega].
\end{equation*}
Notice also that $\Re\sk{\pi_\omega u}{v}=\Re\sk{u}{\pi_\omega v}$.

The background behind the definition of operators $\Phi_\omega$ is explained in Part \ref{part: Roots}; see e.g. Theorem \ref{t: nabla |f|^{p-2}f}.

\begin{remark}
\label{r: Phi as an operator}
We note that our $\Phi_\omega \xi$ corresponds to the object denoted $\xi(\omega)$ in the Dindo\v{s}--Li--Pipher notation \cite[p.373, (18)]{DiLiPi20}; see Section \ref{s: comparison} for a more detailed comparison.
As indicated by our construction, we prefer to interpret this object in terms of a family of $\R$-linear operators $\Phi_\omega$, parametrized by $\omega\in\Sigma^m$ and acting on $\xi\in\C^m_n$.
\end{remark}

This view (of $\Phi_\omega$ as projections) is reflected in the following statement, which follows directly from \eqref{eq: Phi_omega}.

\begin{proposition}
\label{p: Phi properties}
Let $\omega\in\Sigma^m$. Then
\begin{enumerate}[i)]
\item
\label{p: Phi prop 1}
$\Phi_\omega$ is $\R$-linear;
\item
\label{p: Phi prop 2}
$\Phi_\omega$ is a contraction, i.e., $|\Phi_\omega\xi|\leq|\xi|$ for all $\xi\in\C_{n}^{m}$;
\item
\label{p: Phi prop 3}
$\Phi_\omega$ is an idempotent, i.e., $\Phi_\omega\circ\Phi_\omega=\Phi_\omega$.
\end{enumerate}
\end{proposition}

In the sequel, we will frequently encounter expressions of the type $\Re\sk{A\xi}{\Phi_\omega\xi}_{\C_n^m}$. 
They can be interpreted in several ways, as we show next.

\begin{proposition}
\label{p: Phi bilinear form}
For $X,Y\in\C_n^m$ and $\omega\in\Sigma^m$ we have
$$
\aligned
\Re\sk{X}{\Phi_\omega Y}_{\C_n^m}
&=
\Re\sk{\Phi_\omega X}{Y}_{\C_n^m}\\
&=
\sk{\Phi_\omega X}{\Phi_\omega Y}_{\C_n^m}\\
&=
\sk{\Re[X,\omega]}{\Re[Y,\omega]}_{\R^n}\\
&=
\sum_{j=1}^n
\Re\sk{X_j}{\omega}_{\C^m}\Re\sk{Y_j}{\omega}_{\C^m}\\
&=
\sum_{j=1}^n
\sk{\wt{X_j}}{\wt\omega}_{\R^{2m}}\sk{\wt{Y_j}}{\wt\omega}_{\R^{2m}}\\
&=
\sk{S\wt\omega}{\wt\omega}_{\R^{2m}}\,,
\endaligned
$$
where
$$
S=S_{X,Y}
=
\frac12\sum_{j=1}^n \big(\wt{X_j}\cdot\wt{Y_j}^T+\wt{Y_j}\cdot\wt{X_j}^T\big).
$$
\end{proposition}

\begin{proof}
We have, by the second line in \eqref{eq: Phi up/down},
\begin{equation}
\label{p: Phi bilinear - expansion}
\aligned
\Re\sk{X}{\Phi_\omega Y}_{\C^m_n}
&=
\sum_{j=1}^n
\Re\sk{X_j}{(\Phi_\omega Y)_j}_{\C^m}\\
&=
\sum_{j=1}^n
\Re\sk{X_j}{\Re\sk{Y_j}{\omega}_{\C^m}\omega}_{\C^m}\\
&=
\sum_{j=1}^n
\Re\sk{X_j}{\omega}_{\C^m}\Re\sk{Y_j}{\omega}_{\C^m}.
\endaligned
\end{equation}
By symmetry, this establishes the first identity.

Now we jump to the third identity: applying the first line in \eqref{eq: Phi up/down} gives
$$
\aligned
\sk{\Phi_\omega X}{\Phi_\omega Y}_{\C_n^m}
&=
\sum_{\alpha=1}^m
\sk{(\Phi_\omega X)^\alpha}{(\Phi_\omega Y)^\alpha}_{\C_n}\\
&=
\sum_{\alpha=1}^m
\sk{\omega^\alpha\Re[X,\omega]}{\omega^\alpha\Re[Y,\omega]}_{\C_n}\\
&=
\sum_{\alpha=1}^m
\omega^\alpha\overline{\omega}^\alpha
\sk{\Re[X,\omega]}{\Re[Y,\omega]}_{\C_n}
\\
&=
\sk{\Re[X,\omega]}{\Re[Y,\omega]}_{\R_n}.
\endaligned
$$

The second equality follows from the fact that, by Proposition \ref{p: Phi properties} \ref{p: Phi prop 3} and the first identity,
$$
\Re\sk{X}{\Phi_\omega Y}_{\C_n^m}
=
\Re\sk{X}{\Phi_\omega\Phi_\omega Y}_{\C_n^m}
=
\Re\sk{\Phi_\omega X}{\Phi_\omega Y}_{\C_n^m}.
$$
By the third identity that we just proved, $\sk{\Phi_\omega X}{\Phi_\omega Y}_{\C_n^m}$ is already a real number, thus taking the real part in the last term is superfluous.

The fourth equality follows from the previous ones and \eqref{p: Phi bilinear - expansion}.

The fifth representation follows by recalling \eqref{eq: real vs. complex inner product}.

For the last one apply Proposition \ref{p: real symmetric}.
\end{proof}

\begin{proposition}
\label{p: pi - unitary similarity}
Assume the notation introduced in this section.
\begin{enumerate}[(i)]
\item
If the linear map $R:\C^m\rightarrow\C^m$ is unitary, then 
$$
\pi_{R\omega}=R\circ\pi_\omega\circ R^*.
$$

\item
If $\omega_1,\hdots,\omega_k\in\Sigma^m$ then for $u\in\C^m$ we have
$$
\left(
	\pi_{\omega_1}\circ\hdots\circ\pi_{\omega_k}
\right)
(u)
=\Re\sk{u}{\omega_k}\Re\sk{\omega_k}{\omega_{k-1}}\hdots\Re\sk{\omega_2}{\omega_{1}}\omega_1.
$$
\end{enumerate}

\end{proposition}

The proof is straightforward, hence we omit it.

\begin{proposition}
We have, as a consequence of \eqref{eq: Phi up/down},
$$
\left[(\Phi_\omega\circ\Phi_\sigma)(\xi)\right]_h
=
\Re\sk{\xi_h}{\sigma}_{\C^m}
\Re\sk{\sigma}{\omega}_{\C^m}
\omega.
$$
\end{proposition}

\subsection{The operators $\cJ_{p,\omega}$}
\label{s: J_{p,omega}}

For $r\in\R$ and $\omega\in\Sigma^m$ define, following \cite{CD-DivForm}, the operators 
$\cJ_{r,\omega}:
\C_n^m
\rightarrow 
\C_n^m$ 
by
\begin{equation}
\label{eq: J_pvect}
\cJ_{r,\omega}X
=
X+(r-2)
\Phi_\omega X.
\end{equation}
Of course, $\cJ_{r,\omega}$ also admit the ``tensor structure'' akin to \eqref{eq: Phi_omega}.
Indeed, define the maps $\iota_{r,\omega}:\C^m\rightarrow\C^m$ by
\begin{equation}
\label{eq: iota}
\iota_{r,\omega}={\rm id}+(r-2)\pi_\omega.
\end{equation}
Then 
\begin{equation}
\label{eq: J vs. iota}
\cJ_{r,\omega}=\iota_{r,\omega}^{\oplus n}
=\underbrace{\iota_{r,\omega}\oplus\hdots\oplus\iota_{r,\omega}}_{n}.
\end{equation}

Of special interest for us are the cases $r>1$ and the case $r=0$, that is, 
\begin{equation}
\label{eq: J_{0,omega}}
\cJ_{0,\omega}X=X-2\Phi_\omega X\,.
\end{equation}

By the Weinstein–Aronszajn identity, the $(2m)\times(2m)$ matrix, corresponding to $\wt\iota_{r,\omega}$,
namely, $I+(r-2)\wt\omega\cdot\wt\omega^T$, has for any $\wt\omega\in S^{2m-1}$ determinant $r-1$. Thus $\iota_{r,\omega}$ is invertible, unless $r=1$. Hence, the same holds for $\cJ_{r,\omega}$. In fact, we can say more.

\begin{corollary}
\label{c: JoJ}
For any $r,s\in\R$ and $\omega\in\Sigma^m$ we have
\begin{alignat}{2}
\label{eq: JoJ}
\cJ_{r,\omega}\circ\cJ_{s,\omega}
=\cJ_{1+(r-1)(s-1),\omega}.
\end{alignat}
In particular, if $pq=p+q$ (which is $1/p+1/q=1$ if $p,q\ne0$)  
then 
$$
\cJ_{p,\omega}\circ\cJ_{q,\omega}=\rm{id}.
$$
In particular, $\cJ_{0,\omega}\circ\cJ_{0,\omega}=\rm{id}$. 
\end{corollary}

\begin{proof}
Indeed, for arbitrary $\xi\in\C_n^m$ we get from  Proposition \ref{p: Phi properties} {\it \ref{p: Phi prop 3}} 
\begin{alignat*}{2}
\cJ_{r,\omega}(\cJ_{s,\omega}\xi)
&=\cJ_{r,\omega}(\xi+(s-2)\Phi_\omega \xi)\\
&=\left(\xi+(s-2)\Phi_\omega \xi\right) + (r-2)\Phi_\omega\left(\xi+(s-2)\Phi_\omega \xi\right)\\
&=\xi+\left((r-2)+(s-2)+(r-2)(s-2)\right)\Phi_\omega \xi\\
&=\xi+\left((r-1)(s-1)-1\right)\Phi_\omega \xi\\
&=\cJ_{1+(r-1)(s-1),\omega}\xi.
\tag*{\qedhere}
\end{alignat*}
\end{proof}

\begin{remark}
Note that for $\gJ_{p,\omega}:=\cJ_{p+1,\omega}$, the identity \eqref{eq: JoJ} reads 
$\gJ_{r,\omega}\circ\gJ_{s,\omega}=\gJ_{rs,\omega}$.
In particular, for all $r\in\R\backslash\{0\}$ we have $\gJ_{r,\omega}\circ\gJ_{1/r,\omega}=\rm{id}$. That is, 
$$
\gJ_{r,\omega}^{-1} = \gJ_{r^{-1},\omega}
\hskip 30pt
\forall r\in\R\backslash\{0\}, \omega\in\Sigma^m.
$$
\end{remark}

\begin{corollary}
\label{c: norm of J}
For $p\ne1$ we have 
\begin{equation*}
|\cJ_{p,\omega}X|\sim|X|. 
\end{equation*}
More precisely, for every $p\in\R$, $X\in\C_n^m\backslash\{0\}$ and $\omega\in\Sigma^m$ we have
$$
\min\{1,|p-1|\}
\leq
\frac{|\cJ_{p,\omega}X|}{|X|}
\leq
\max\{1,|p-1|\}
\,.
$$
Both estimates are optimal.

In particular, for every $X\in\C_n^m$ and $\omega\in\Sigma^m$ we have
\begin{equation}
\label{eq: J_0 isometry}
|\cJ_{0,\omega}X|=|X|.
\end{equation}
\end{corollary}

\begin{proof}
From \eqref{eq: J vs. iota} it immediately follows that it suffices to prove the statements with $\cJ_{r,\omega}$ replaced by $\iota_{r,\omega}$ and $X\in\C_n^m\backslash\{0\}$ by $\zeta\in\C^m\backslash\{0\}$.

For the \underline{upper estimate} use \eqref{eq: iota} and \eqref{eq: pi_omega} to calculate
\begin{equation}
\label{eq: norm of iota}
\aligned
|\iota_{p,\omega}\zeta|^2
& = \sk{\zeta+(p-2)\Re\sk{\zeta}{\omega}_{\C^m}
\omega}{\zeta+(p-2)\Re\sk{\zeta}{\omega}_{\C^m}
\omega}_{\C^m}\\
&=|\zeta|^2+p(p-2)\Re\sk{\zeta}{\omega}^2\,.
\endaligned
\end{equation}
If $p\in[0,2]$ (that is, $p(p-2)\leq0$ or $|p-1|\leq1$), we have
$$
|\iota_{p,\omega}\zeta|^2\leq |\zeta|^2\,.
$$
If $p\in(-\infty,0]\cup[2,\infty)$ (that is, $p(p-2)\geq0$ or $|p-1|\geq1$), we have
$$
|\iota_{p,\omega}\zeta|^2\leq |\zeta|^2+p(p-2)|\zeta|^2=(p-1)^2|\zeta|^2\,.
$$
Equalities are obtained above by taking 
\begin{itemize}
\item
$\wt\zeta\in\wt\omega^\perp$ if $|p-1|\leq1$;
\item
$\zeta=\omega$ \hskip 6pt if $|p-1|\geq1$.
\end{itemize}

For the \underline{lower estimate} we use Corollary \ref{c: JoJ} and the part proven so far: if $p\ne1$ and $q=p/(p-1)$ then
$$
\aligned
|\zeta|
=|\iota_{q,\omega}
\big(\iota_{p,\omega}\zeta\big)|
\leq
\max\{1,|q-1|\}|\iota_{p,\omega}\zeta|.
\endaligned
$$
Hence
\begin{alignat*}{2}
|\iota_{p,\omega}\zeta|&\geq\frac{|\zeta|}{\max\{1,|q-1|\}}=\min\{1,|p-1|\}|\zeta|\,.
\end{alignat*}
This estimate extends to the case $p=1$ through \eqref{eq: norm of iota}, which gives
$$
|\iota_{1,\omega}\zeta|^2
=|\zeta|^2-\Re\sk{\zeta}{\omega}^2\,.
$$

Equalities are obtained above by taking 
\begin{itemize}
\item
$\zeta=\omega$ \hskip 6pt if $|p-1|\leq1$;
\item
$\wt\zeta\in\wt\omega^\perp$ if $|p-1|\geq1$.\hfill\qedhere
\end{itemize}
\end{proof}


\subsection{The operators $\cI_p,\cJ_p,\cI_{p,\omega},\cJ_{p,\omega}$}
\label{s: aux ops}

Let us start with the following identity, valid for $\xi\in\C^n$ and $p\in\R\backslash\{0\}$:
\begin{equation}
\label{eq: bar vs. Re}
\xi+(1-2/p)\bar\xi=\frac2p\left(\xi+(p-2)\Re\xi\right).
\end{equation}
For any $p\in(1,\infty)$ define operators $\cI_p,\cJ_p:\C^k\rightarrow\C^k$ by
$$
\aligned
\cI_p\eta & =\eta+(1-2/p)\bar \eta\\
\cJ_p\eta & =\eta+(p-2)\Re\eta.
\endaligned
$$
The operator $\cI_p$ appeared explicitly, without changes, in \cite{CD-DivForm}. 
Now \eqref{eq: bar vs. Re} reads
\begin{equation}
\label{eq: I vs. J}
\cI_p=\frac2p\,\cJ_p.
\end{equation} 
This said, one thinks of defining, for $\omega\in\Sigma^m$, the operators $\cI_{p,\omega}:\C_n^m\rightarrow\C_n^m$ by
\begin{equation}
\label{eq: I vs. J vect}
\cI_{p,\omega}:=\frac2p\,\cJ_{p,\omega}\,.
\end{equation}

Direct calculations involving \eqref{eq: I vs. J} and \eqref{eq: I vs. J vect} yield the following identities.

\begin{lemma}
\label{l: I,J vect vs. I,J scal}
For any $p>1$ and $X\in\C_n^m$, $\omega\in\Sigma^m$ we have
$$
\aligned
\cI_{p,\omega} X
&=X-(1-2/p)\cJ_{0,\omega}X
&&\hskip -11pt =\cI_{p}([X,\omega])\omega+\frac2p(X-[X,\omega]\omega)\\
\cJ_{p,\omega}X
&=X+ (p-2)\Phi_\omega X
&&\hskip -11pt =\cJ_{p}([X,\omega])\omega+\hskip7.5pt (X-[X,\omega]\omega).
\endaligned
$$
\end{lemma}

\begin{remark}
\label{r: I,J at m=1}
When \framebox{$m=1$}, the above operators simplify significantly. Namely, when $m=1$, the term $X-[X,\omega]\omega$ vanishes for all $X\in\C_n$ and $\omega\in\Sigma^1$, and thus we get 
$$
\aligned
\cI_{p,\omega} X
&=\cI_{p}([X,\omega])\omega
&&\hskip -10pt=\omega\cI_{p}(\bar\omega X)\\
\cJ_{p,\omega} X
&=\cJ_{p}([X,\omega])\omega
&&\hskip -10pt=\omega\cJ_{p}(\bar\omega X).
\endaligned
$$
That is,
$$
\aligned
\cJ_{p,\omega} & = \cM_\omega \cJ_p \cM_{\omega}^{-1}\\
\cI_{p,\omega} & = \cM_\omega \cI_p \cM_{\omega}^{-1},
\endaligned
$$
where $\cM_\omega:\C^n\rightarrow\C^n$ is the multiplication operator $\xi\mapsto\omega\xi$.
\end{remark}


\subsection{The operators $I_{p,\omega}$}

For $p>0$ and $\omega\in \Sigma^m$, we define the $\R$-linear operator $I_{p,\omega}:\C^m\rightarrow\C^m$ by
$$
I_{p,\omega}:=\frac{2}{p}(\mathrm{Id}_{\C^m}+(p-2)\pi_\omega ).
$$
The operators $I_{p,\omega}$ admit a precise geometric interpretation in terms of orthogonal reflections. To see this, we first define the complementary projections
$$
\pi_{\omega^\perp}=\mathrm{Id}_{\C^m}-\pi_{\omega}\quad \textrm{and}\quad  \mathrm{\Phi_{\omega^\perp}}=\mathrm{Id}_{\C^m_n}-\Phi_{\omega}=\pi_{\omega^\perp}\otimes\mathrm{Id}_{\C_n}.
$$
\begin{remark}
The operator $\pi_{\omega^\perp}$ is the complex counterpart, via conjugation with the identification operator $\cV_m$, of the orthogonal projection $\Pi_{\widetilde{\omega}^\perp}$ onto the subspace $\textrm{span}\{\widetilde{\omega}\}^\perp$ of $\R^{2m}$. 
\end{remark}

We now consider the $\R$-linear operator $R_\omega: \C^m \rightarrow \C^m$, which corresponds to the orthogonal reflection in $\R^{2m}$ across the line spanned by $\widetilde{\omega}$. We  further denote by $\cR_\omega$ its extension to the tensor space $\C^m_n$. More precisely, we set
\begin{equation*}
\label{eq: reflections}
R_\omega:=\pi_\omega-\pi_{\omega^\perp}\quad \textrm{and}\quad \cR_\omega:=R_\omega\otimes\mathrm{Id}_{\C_n}=\Phi_\omega-\Phi_{\omega^\perp}.
\end{equation*}
In short, 
$$
\cR_\omega=-\cJ_{0,\omega};
$$
recall that, up to the sign, the latter operators were already introduced in \eqref{eq: J_{0,omega}}. 

By utilizing the identity $2\pi_\omega-\mathrm{Id}_{\C^m}=R_\omega$, a direct calculation shows that for any $p>0$ and $\omega \in \Sigma^m$, we have
\begin{equation}
\label{eq: J and R}
I_{p,\omega}=\mathrm{Id}_{\C^m}+\left(1-2/p\right)R_\omega\quad \textrm{and}\quad \cI_{p,\omega}=\mathrm{Id}_{\C^m_n}+\left(1-2/p\right)\cR_\omega.
\end{equation}

Alternatively, letting $q = p/(p-1)$ be the conjugate exponent of $p$, the operator $\cI_{p,\omega}$ can be expressed as a weighted sum of the complementary projections:
\begin{equation*}
\cI_{p,\omega} = \frac{2}{q} \Phi_\omega + \frac{2}{p} \Phi_{\omega^\perp}.
\end{equation*}
From this representation, the symmetry between $p$ and $q$ is immediate, yielding the identity
\begin{equation*}
\cI_{p,\omega}+\cI_{q,\omega}=2\mathrm{Id}_{\C^m_n}.
\end{equation*}

\subsection{The operator $R_\omega$ as a unitary conjugation}

It is natural to regard the real orthogonal group $O(2m)$ as the group of
$\R$-linear isometries of $\C^m$. Accordingly, we define the set
of real-linear orthogonal operators on $\C^m$ by
\begin{equation*}
O_{\R}(\C^m) := \mathcal{V}_m^{-1} O(2m) \mathcal{V}_m.
\end{equation*}
An operator $O \in O_{\R}(\C^m)$ is characterized by the property
that it preserves the real part of the Hermitian inner product:
\begin{equation*}
\label{eq: real_isometry}
\Re\sk{Ou}{Ov}_{\C^m}
= \Re\sk{u}{v}_{\C^m},
\quad \forall\, u,v \in \C^m.
\end{equation*}
Equivalently, $O$ is an isometry of $\C^m$ endowed with its underlying
real Euclidean structure. The unitary group $U(m)$, consisting of all
$\C$-linear isometries, is a proper subgroup of
$O_{\R}(\C^m)$; indeed,
$$
U(m) = O_{\R}(\C^m) \cap GL(m,\C).
$$

By basically repeating Proposition \ref{p: pi - unitary similarity}, we see that any $O \in O_{\R}(\C^m)$ and $\omega \in \Sigma^m$, one has
$$
\pi_{O\omega} = O \circ \pi_\omega \circ O^{-1}.
$$
In particular, if $U \in U(m)$, then $\pi_{U\omega} = U \circ \pi_\omega \circ U^*$.

Owing to the linear relation $R_\omega = 2\pi_\omega - \mathrm{Id}$, it follows that the reflection operators satisfy the same conjugation law: 
$$
R_{O\omega} = O \circ R_\omega \circ O^{-1},
$$ 
for all $O \in O_{\R}(\C^m)$ and $\omega\in\Sigma^m$. Since $U(m)$ acts transitively on the unit sphere $\Sigma^m$, for every $\omega\in\Sigma^m$ there exists a unitary operator $U_\omega$ such that $U_\omega e_1 = \omega$, where 
$$
e_1 = (1, 0, \dots, 0)^T\in\C^m.
$$
Under this unitary change of basis, $R_\omega$ is conjugate to the reflection $R_{e_1}$. In the $e_1$ basis, the operator simplifies significantly:
\begin{equation*}
R_{e_1}(u^1,\dots,u^m)^T=\left(\overline{u^1},-u^2,\dots,-u^m\right)^T,
\end{equation*}
for all $u=(u^1,\dots,u^m)^T\in\C^m$. This establishes the following result.

\begin{proposition}
\label{p: conj}
Let $\kappa:\C\rightarrow\C$ denote the conjugation $\kappa (z)=\overline{z}$. Consider the block diagonal operator $B\in O_\R(\C^m)$ defined by 
$$
B:=\kappa\oplus (-\mathrm{Id}_{\C^{m-1}}).
$$
Then, $B=R_{e_1}$ and, for every $\omega\in \Sigma^m$, the reflection $R_\omega$ is unitarily conjugate to $B$. Namely,
$$
R_\omega=U_\omega B U^{*}_\omega,
$$
where $U_\omega\in U(m)$ is any unitary operator such that $U_\omega e_1=\omega$.
\end{proposition}

\subsection{The operator $\cI_{p,\omega}$ as a unitary conjugation}
Maintaining the notation of the previous sections, we denote by $\cB$ and $\cU$ the extensions to $\C^m_n$ of the operators $B$ and $U\in U(m)$, respectively:
$$
\cB=B\otimes\mathrm{Id}_{\C_n}=R_{e_1}\otimes\mathrm{Id}_{\C_n}=-\cJ_{0,e_1};\quad \cU=U\otimes\mathrm{Id}_{\C_n}.
$$
It follows from \eqref{eq: J and R} and Proposition~\ref{p: conj} that
\begin{equation}
\label{eq: cJ and unitary}
\aligned
&I_{p,\omega}=U_\omega\left(\mathrm{Id}_{\C^m}+\left(1-2/p\right)B\right)U^*_\omega\\ &\cI_{p,\omega}=\cU_\omega\left(\mathrm{Id}_{\C^m_n}+\left(1-2/p\right)\cB\right)\cU^*_\omega.
\endaligned
\end{equation}
We emphasize that $\cB$, $\cU_\omega$, and $\cR_\omega$ preserve the real part of the scalar product in $\C^m_n$. Moreover, $\cU_\omega$ is an element of the unitary group on $\C^m_n$. These facts follow immediately from observing that $\cB$ and $\cU_\omega$ are obtained by tensorizing with $\mathrm{Id}_{\C_n}$ the operators $B\in O_\R(\C^m)$ and $U_\omega\in U(m)$, respectively, and that 
$$
\cR_\omega=\cU_\omega\cB\cU^*_\omega.
$$
Note that, in the {\it up} representation, $\cB$ acts as
\begin{equation}
\label{eq: K in up}
\cB(\xi^1,\xi^2,\dots,\xi^m)^T=(\overline{\xi^1},-\xi^2,\dots,-\xi^m)^T.
\end{equation}
Equivalently, in components,

$$
(\cB\xi)^\alpha=
\begin{cases}
\overline{\xi^\alpha},&\alpha=1,\\
-\xi^\alpha,&\alpha\neq 1.
\end{cases}
$$

\part{Definition and properties of $p$-ellipticity}

\label{part: Def Props}

\section{The concept of $p$-ellipticity for systems of elliptic differential equations in divergence form} 
\label{s: p-ell}

In this section we revisit the concept of $p$-ellipticity, encountered in Section \ref{s: scalar p-ell} in the realm of elliptic equations, and extend it to elliptic systems.
We first present our construction; later, in Part \ref{part: Roots} we will explain how it was conceived and hence why we find it natural. 
Briefly, we wanted to find a sufficient algebraic condition that allows the semigroup contractivity, which is known to hold on $L^2$, to extrapolate to $L^p$, along the lines of the corresponding theorem in the scalar case, see \cite[Theorem 1.3]{CD-DivForm}. We will show that this naturally leads to the definition which we describe next.

\begin{definition}
Recall the notation from Section \ref{s: spheres}.
Suppose that $A$ is a bounded (cf. Section \ref{s: Legendre}) complex tensor on $\Omega$ and $p\in(1,\infty)$.
For $\Phi_{\omega}$ as in \eqref{eq: Phi up/down} define
\begin{equation}
\label{eq: lambda_p}
\lambda_p(A):=
\underset{x\in\Omega}{{\rm ess}\inf}
\min_{\xi\in\Sigma_n^m}
\min_{\omega\in\Sigma^m}
\Re\sk{A(x)\xi}{\xi+(p-2)\Phi_\omega\xi}_{\C^m_n}.
\end{equation}
That is,
\begin{equation}
\label{eq: lambda_p via J_p}
\lambda_p(A)=
\underset{x\in\Omega}{{\rm ess}\inf}
\min_{\xi\in\Sigma_n^m}
\min_{\omega\in\Sigma^m}
\Re\sk{A(x)\xi}{\cJ_{p,\omega}\xi}_{\C^m_n}.
\end{equation}
We say that $A$ is {\it $p$-elliptic} if $\lambda_p(A)>0$. This is equivalent to requiring that
\begin{equation}
\label{eq: p-ell}
\aligned
&\exists\, \lambda_p(A)>0\ \ni: \ {\rm a.e.}\ x\in\Omega:\\
&\hskip 20pt
\forall\xi\in\C_n^m,\ \forall\omega\in\Sigma^m:\hskip 30pt
\Re\sk{A(x)\xi}{\xi+(p-2)\Phi_\omega\xi}\geq\lambda_p(A)|\xi|^2.
\endaligned
\end{equation}
We denote the family of all $p$-elliptic tensors on $\Omega$ acting on $\C^m_n$ by $\cA_{p}^m(\Omega)$.

\medskip
In accordance with \cite{CD-DivForm} we also define
\begin{equation}
\label{eq: Delta vs. lambda}
\Delta_p(A)=\frac 2p\lambda_p(A).
\end{equation}
See Section \ref{s: aux ops} for more insight into this relation.
\end{definition}

By the unitary conjugation property \eqref{eq: cJ and unitary}, we may equivalently write
\begin{equation}
\label{eq: Deltap bis}
\Delta_p(A)=\underset{x\in\Omega}{{\rm ess}\inf}
\min_{\xi\in\Sigma_n^m 
}
\min_{\omega\in\Sigma^m
}
\Re\Sk{\cU^*_\omega A(x)\cU_\omega\xi}{\xi+(1-2/p)\cB\xi}_{\C^m_n}.
\end{equation}

\subsection{Comparison with the scalar $p$-ellipticity}

By Proposition \ref{p: Phi bilinear form} we may write, for any tensor $A$,
\begin{equation}
\label{eq: lambda_p scalar}
\lambda_p(A)=
\min_{\xi\in\Sigma_n^m}
\min_{\omega\in\Sigma^m}
\left(
\Re\sk{A\xi}{\xi}_{\C_n^m}+(p-2)\sk{\Re[A\xi,\omega]}{\Re[\xi,\omega]}_{\R_n}
\right).
\end{equation}
We note that when \framebox{$m=1$}, by taking $\eta=[\xi,\omega]=\overline\omega\xi$, the expression above becomes
$$
\Re\sk{A\eta}{\eta}_{\C_n}+(p-2)\sk{\Re(A\eta)}{\Re\eta}_{\R_n},
$$
that is, it depends on $\overline\omega\xi$ only. Thus since $|\eta|=|\xi|$, after taking the minimum, the action of the sphere $\Sigma^m=\Sigma^1$ disappears.  
A trivial reason behind this is that when $m=1$, the projections $[X,\omega]\omega$ recover $X$, meaning simply that
$$
X=[X,\omega]\omega
$$
regardless of the choice of $X,\omega$, which allows for the above simplification. See Lemma \ref{l: I,J vect vs. I,J scal} and Remark \ref{r: I,J at m=1}.

The situation for $m>1$ is markedly different.

\medskip
In different words, when $m=1$, the quantity in \eqref{eq: Deltap bis} coincides with the number $\Delta_p(A)$ in our definition of $p$-ellipticity for complex matrix-valued functions in \cite{CD-DivForm} because in the scalar case, the unitary map $U_\omega$ acts simply as multiplication by the phase factor:
$$
U_\omega u=e^{i\mathrm{arg}(\omega)}u,\quad \forall \omega\in\Sigma^1,\ \forall u\in\C^1.
$$ 
Consequently, its tensor extension acts as
$$
\cU_\omega\xi=e^{i\mathrm{arg}(\omega)}\xi,\quad \forall \omega\in\Sigma^1,\ \forall\xi\in\C^1_n.
$$
Crucially, for $m=1$, the tensor $A(x)$ is a complex $n\times n$ matrix ($A(x)\in\C^1_n\otimes\C^1_n$) and thus trivially commutes with the multiplication by $e^{i\theta}$. This implies that in \eqref{eq: Deltap bis} the action of the sphere $\Sigma^1$ disappears. Furthermore, in the scalar case, the operator $B$ is the standard complex conjugation (denoted by $\kappa$), so 
$$
\cB\xi=\overline{\xi},\quad \forall \xi\in\C^1_n.
$$

Hence, for $m=1$, the expression \eqref{eq: Deltap bis} reduces to
$$
\aligned
\Delta_p(A)&=\underset{x\in\Omega}{{\rm ess}\inf}
\min_{\xi\in\Sigma_n^1 
}
\Re\sk{A(x)\xi}{(\mathrm{Id}+(1-2/p)\cB)\xi}_{\C^1_n}\\
&=\underset{x\in\Omega}{{\rm ess}\inf}
\min_{\xi\in\Sigma_n^1 
}
\Re\sk{A(x)\xi}{\xi+(1-2/p)\overline{\xi}}_{\C^1_n},
\endaligned
$$
which coincides with our definition of scalar $p$-ellipticity in \cite{CD-DivForm}.

The situation for $m>1$ is different because, in general, rank-$4$ tensors $A(x)$ need not commute with the unitary transformations $\cU_\omega$. Therefore, in the definition of $\Delta_p(A)$, it is necessary to minimize over the unitary orbit of the tensor $A(x)$ generated by the operators $\cU_\omega$, $\omega\in\Sigma^m$.

\begin{remark}
\label{r: commute}    
The condition that a rank-4 tensor $A$ commutes with $\mathcal{U}_\omega$ for all $\omega \in \Sigma^m$ is highly restrictive. Indeed, by Schur's Lemma, such commutation holds if and only if $A = \mathrm{Id}_{\C^m} \otimes C$ for some $C\in\mathrm{End}(\C_n)$. Let us show this.

Writing $A$ in the down representation, we see that $\cU_\omega A=A\cU_\omega$ is equivalent to 
$$
U_\omega A_{j,k}=A_{j,k}U_\omega,\quad \forall j,k=1,\dots,n.
$$
Since the action of $\{U_\omega:\ \omega\in\Sigma^m\}$ on $\C^m$ is irreducible, Schur's Lemma implies that each block $A_{j,k}$ is a multiple of the identity:
$$
A_{j,k}=c_{jk}\mathrm{Id}_{\C^m}.
$$
Equivalently, if we set $C=(c_{jk})$, we have 
$$
A^{\alpha\beta}_{jk}=\delta^{\alpha\beta}c_{jk}.
$$
Therefore, $A^{\textrm{up}}=\textrm{diag}(C,\dots,C)$. This (see Remark~\ref{r: up/down}) can be restated as
$$
A=\mathrm{Id}_{\C^m}\otimes C,
$$
where $C\in\mathrm{End}(\C_n)$ is represented by the matrix $(c_{jk})$ in the canonical basis of $\C_n$.

In this specific case, the minimization over $\Sigma^m$ in the definition of $\Delta_p(A)$ becomes redundant, and as we will prove in Corollary~\ref{c: A^m vs. A p-ell},
we recover the scalar theory applied to $C$. 
\end{remark}

\subsection{Operators $\Theta_A^\pm$}
Define
\begin{equation}
\label{eq: Theta^pm}
\aligned
\Theta^-(\xi)&=\Theta^-_A(x,\xi):=
\min_{\omega\in\Sigma^m}
\Re\sk{A(x)\xi}{\Phi_\omega\xi}_{\C^m_n}\\
\Theta^+(\xi)&=\Theta^+_A(x,\xi):=
\max_{\omega\in\Sigma^m}
\Re\sk{A(x)\xi}{\Phi_\omega\xi}_{\C^m_n}.
\endaligned
\end{equation}
Then we may rewrite \eqref{eq: lambda_p} as 
\begin{equation}
\label{eq: lambda via Theta}
\lambda_p(A)=
\left\{
	\begin{array}{rcl}
		{\displaystyle
			\underset{x\in\Omega}{{\rm ess}\inf}
				\min_{\xi\in\Sigma_n^m
				}
					\left(
						\Re\sk{A(x)\xi}{\xi}+(p-2)\Theta_A^+(x,\xi)
					\right)
		} 
& ; & \text{ if }p\leq2\\
		{\displaystyle
			\underset{x\in\Omega}{{\rm ess}\inf}
				\min_{\xi\in\Sigma_n^m}
					\left(
						\Re\sk{A(x)\xi}{\xi}+(p-2)\Theta_A^-(x,\xi)
					\right)
		} 
& ; & \text{ if }p\geq2.
	\end{array}
\right.
\end{equation}

We can get around the action of the unit sphere $\Sigma^m\subset\C^m$ in \eqref{eq: Theta^pm}, as follows.

\begin{theorem}
\label{t: S}
Given a tensor $A$ and an element $\xi\in\C_n^m$, introduce the real $(2m)\times(2m)$ matrix
$$
S=S_{A,\xi}
=
\frac12\sum_{j=1}^n \big(\wt{(A\xi)}_j\cdot\wt\xi_j^T+\wt\xi_j\cdot\wt{(A\xi)}_j^T\big).
$$
Then, for every $\omega\in\Sigma^m$,
\begin{equation*}
\Re\sk{A\xi}{\Phi_\omega\xi}_{\C^m_n}
=
\sk{S\wt\omega}{\wt\omega}_{\R^{2m}}\,.
\end{equation*}
Consequently,
\begin{equation*}
\aligned
\Theta^-_A(x,\xi)
=
\min_{\omega\in\Sigma^m}
\Re\sk{A\xi}{\Phi_\omega\xi}_{\C^m_n}
& =
\text{the smallest eigenvalue of }
S\\
\Theta^+_A(x,\xi)
=
\max_{\omega\in\Sigma^m}
\Re\sk{A\xi}{\Phi_\omega\xi}_{\C^m_n}
& =
\text{the largest eigenvalue of }
S.
\endaligned
\end{equation*}
\end{theorem}

\begin{proof}
The first statement is a special case of Proposition \ref{p: Phi bilinear form}, applied with $X=A\xi$ and $Y=\xi$. 
The second statement follows from Proposition \ref{p: real symmetric}.
\end{proof}

Proposition \ref{p: Phi bilinear form} and the displays \eqref{eq: Theta^pm}, \eqref{eq: lambda via Theta} were the reason for studying inequalities from Corollary \ref{c: bilateral est}. Indeed, in the special case of  
$$
\aligned
(u_1,\hdots,u_n)&:=A\xi\\
(v_1,\hdots,v_n)&:=\xi,
\endaligned
$$
we obtain the following result:

\begin{corollary}
\label{c: A form vs. Theta}
For any bounded tensor $A$ and any $\xi\in\C_n^m$,
$$
\Theta^-_A(\xi)
\leq
\frac{1}{2m}\,
\Re\sk{A\xi}{\xi}_{\C^m_n}
\leq
\Theta^+_A(\xi).
$$
\end{corollary}
\begin{proof}
We get from Proposition \ref{p: Phi bilinear form} that
$$
\Re\sk{A\xi}{\Phi_\omega\xi}_{\C^m_n}
=
\sum_{h=1}^n
\Re\sk{(A\xi)_h}{\omega}_{\C^m}\Re\sk{\xi_h}{\omega}_{\C^m}.
$$
Now recall the definition \eqref{eq: Theta^pm} and apply Corollary \ref{c: bilateral est}.
\end{proof}
Slightly more generally, for $u,v\in\C_n^m$ introduce
\begin{equation*}
\aligned
T^-(u,v) &:= \min_{\omega\in\Sigma^m} \Re\sk{u}{\Phi_\omega v}_{\C^m_n}\\
T^+(u,v) &:= \max_{\omega\in\Sigma^m} \Re\sk{u}{\Phi_\omega v}_{\C^m_n}.
\endaligned
\end{equation*}
Thus
$$
\Theta_A^\pm(\xi)=T^\pm(A\xi,\xi).
$$
The above can be restated as 
\begin{equation}
\label{eq: Re(uv) vs. T}
T^-(u,v) \leq \frac1{2m}\,\Re \sk{u}{v} \leq T^+(u,v)
\end{equation}
for $u,v\in\C_n^m$.

\medskip
Observe that when $m=1$, we can actually {\it calculate} $T^\pm(u,v)$, as specified in \eqref{eq: Theta & m=1}. 
(This was in fact important, because we used this to prove the monotony of the (scalar) $p\mapsto\Delta_p(A)$, see \cite[Proof of Corollary 5.16]{CD-DivForm}.)

Indeed, then we have $\omega\in\Sigma^1$, $\xi\in\C^n$ and
$$
\Phi_\omega\xi 
 = \omega\,\Re(\overline\omega\xi)
 = \frac{\omega(\overline\omega\xi+\omega\overline\xi)}{2}
 = \frac{\xi+\omega^2\overline\xi}{2}\,.
$$
Consequently,
$$
\Re\sk{u}{\Phi_\omega v}=
\frac{\Re\sk{u}{v}+\Re(\overline\omega^2\sk{u}{\overline v})}2,
$$
thus optimization over $\omega\in\Sigma^1$ clearly gives, for $u,v,\in\C^n$,
\begin{equation}
\label{eq: T for m=1}
T^\pm(u,v)=\frac{\Re\sk{u}{v}\pm|\sk{u}{\overline v}|}2\,.
\end{equation}
This implies that already when \framebox{$m=1$}, none of the estimates in \eqref{eq: Re(uv) vs. T} is sharp.

\subsection{Improving \eqref{eq: Re(uv) vs. T}}
We can improve \eqref{eq: Re(uv) vs. T} as follows.

\begin{proposition}
\label{eq: Phi_omega quad est}
For any $m\in\N$ and $u,v\in\C_n^m$ we have 
\begin{equation*}
\aligned
\min_{\omega\in\Sigma^m} \Re\sk{u}{\Phi_\omega v}_{\C^m_n}
\leq & \ \frac1m\sum_{k=1}^m\frac{\Re\sk{u^k}{v^k}-\left|\sk{u^k}{\overline v^k}\right|}2
& \\
&\ \frac1m\sum_{k=1}^m\frac{\Re\sk{u^k}{v^k}+\left|\sk{u^k}{\overline v^k}\right|}2\!\! 
 &\!\!\!   \leq \max_{\omega\in\Sigma^m} \Re\sk{u}{\Phi_\omega v}_{\C^m_n}.
\endaligned
\end{equation*}
In other words,
\begin{equation*}
\aligned
T^-(u,v) & \leq\frac1m\sum_{k=1}^m T^-\left(u^k,v^k\right)\\
T^+(u,v) & \geq\frac1m\sum_{k=1}^m T^+\left(u^k,v^k\right).
\endaligned
\end{equation*}

\end{proposition}

This estimate not only improves \eqref{eq: Re(uv) vs. T}, but also turns into an equality when $m=1$, as we have seen in \eqref{eq: T for m=1}. 

\begin{proof}
We shall use the following identity which was a part of Proposition \ref{p: Phi bilinear form}:
$$
\Re\sk{u}{\Phi_\omega v}_{\C_n^m}
=
\sum_{j=1}^n
\Re\sk{u_j}{\omega}_{\C^m}\Re\sk{v_j}{\omega}_{\C^m}.
$$

In order to shed more light on the proof of the general case, we first consider the case $m=2$ on which the general case is built.

\subsubsection*{The case $m=2$}
We introduce the notation $\omega=(r_1\sigma^1,\hdots,r_m\sigma^m)$ with $r_j\in[-1,1]$ and 
$$
\sigma^j\in
\Sigma^1_+
:=
\mn{\zeta\in\C}{|\zeta|=1\text{ and }\Re\zeta\geq0}\,.
$$
Note that we allow $r_j$ to lie in $[-1,1]$ and not only $[0,1]$, while $\sigma$'s are restricted to the right half of the unit circle in $\C$.

Write
$$
\aligned
a_j & =\Re\big(u_j^1\overline{\sigma^1}\big) & \ \hskip 20pt \ 
c_j & =\Re\big(v_j^1\overline{\sigma^1}\big)\\
b_j & =\Re\big(u_j^2\overline{\sigma^2}\big) & \ \hskip 20pt \  
d_j & =\Re\big(v_j^2\overline{\sigma^2}\big).
\endaligned 
$$
Then
$$
\sum_{j=1}^n
\Re\sk{u_j}{\omega}_{\C^2}
\Re\sk{v_j}{\omega}_{\C^2}
=
\Big(\sum_{j=1}^na_jc_j\Big) r_1^2
+
\Big(\sum_{j=1}^n(a_jd_j+b_jc_j)\Big) r_1r_2
+
\Big(\sum_{j=1}^nb_jd_j\Big) r_2^2.
$$
We quickly check that, for any $A,B,C\in\R$,
$$
\aligned
\min_{r_1,r_2\in[-1,1] \atop r_1^2+r_2^2=1}
(Ar_1^2+Br_1r_2+Cr_2^2)
& =
\frac{A+C-\sqrt{(A-C)^2+B^2}}{2}\\
& \leq 
\frac{A+C}{2}\,.
\endaligned
$$
Note that here we discarded the entire square root.

This gives
$$
\aligned
\min_{r_1,r_2\in[-1,1] \atop r_1^2+r_2^2=1}
\sum_{j=1}^n
& \Re\sk{u_j}{\omega}_{\C^2}
\Re\sk{v_j}{\omega}_{\C^2}\\
& \leq
\sum_{j=1}^n\frac{a_jc_j+b_jd_j}{2}\\
& =
\frac14
\sum_{j=1}^n
\left(
\Re\big(u_j^1\overline{v_j^1}\big)
+
\Re\big(u_j^2\overline{v_j^2}\big)
+
\Re\big(\overline{\sigma^1}^2u_j^1v_j^1\big)
+
\Re\big(\overline{\sigma^2}^2u_j^2v_j^2\big)
\right).
\endaligned
$$
Now, by noting that $\overline{\Sigma_+^1}^2=\Sigma^1$, we continue as
$$
\aligned
\min_{\sigma^1,\sigma^2\in\Sigma_+^1}
& \frac14
\bigg(
\sum_{j=1}^n
\Re\big(u_j^1\overline{v_j^1}\big)
+
\sum_{j=1}^n
\Re\big(u_j^2\overline{v_j^2}\big)
+
\Re\Big(\overline{\sigma^1}^2
\sum_{j=1}^n
u_j^1v_j^1\Big)
+
\Re\Big(\overline{\sigma^2}^2
\sum_{j=1}^n
u_j^2v_j^2\Big)
\bigg)\\
& = 
\frac14
\bigg(
\sum_{j=1}^n
\Re\sk{u_j}{v_j}
-
\Big|
\sum_{j=1}^n
u_j^1v_j^1
\Big|
-
\Big|
\sum_{j=1}^n
u_j^2v_j^2
\Big|
\bigg)\\
& =
\frac1{2}
\sum_{k=1}^2
\frac1{2}
\bigg(
\Re\sk{u^k}{v^k}
-
\Big|
\sum_{j=1}^n
u_j^kv_j^k
\Big|
\bigg).
\endaligned
$$
This finishes the first inequality of Proposition \ref{eq: Phi_omega quad est} when $m=2$. In a similar vein we obtain the second estimate, as well.

\subsubsection*{The case of general $m\in\N$}

Denote $u_j=(u_j^1,\hdots,u_j^m)\in\C^m$ and  $v_j=(v_j^1,\hdots,v_j^m)\in\C^m$.
Write also $\omega=(\omega^1,\hdots,\omega^m)\in\C^m$ with 
$$
\omega^k=r_k\sigma^k,
$$
where $\sigma^k\in\Sigma_+^1$ and $r_k\in[-1,1]$ are such that 
$$
r_1^2+\hdots+r_m^2=1. 
$$
With this notation, $\omega=\big(r_1\sigma^1,\hdots,r_m\sigma^m\big)$ and
\begin{equation}
\label{eq: double min}
\min_{\omega\in\Sigma^m}
\leq
\min_{\sigma^1,\hdots,\sigma^m\in\Sigma^1}
\min_{r_1,\hdots,r_m\in[-1,1] \atop r_1^2+\hdots+r_m^2=1}.
\end{equation}
We then have
$$
\aligned
\sum_{j=1}^n
\Re\sk{u_j}{\omega}\Re\sk{v_j}{\omega}
&=
\sum_{j=1}^n
\Big(
\Re\big(u_j^1\overline{\omega^1}\big)+\hdots+\Re\big(u_j^m\overline{\omega^m}\big)
\Big)
\Big(
\Re\big(v_j^1\overline{\omega^1}\big)+\hdots+\Re\big(v_j^m\overline{\omega^m}\big)
\Big)\\
&=
\sum_{j=1}^n
\big(
a_j^1r_1+\hdots+a_j^mr_m
\big)
\big(
b_j^1r_1+\hdots+b_j^mr_m
\big),
\endaligned
$$
where
$$
\aligned
a_j^k&=\Re\big(u_j^k\overline{\sigma^k}\big)\\
b_j^k&=\Re\big(v_j^k\overline{\sigma^k}\big).
\endaligned
$$
We continue with
$$
\left(
a_j^1r_1+\hdots+a_j^mr_m
\right)
\left(
b_j^1r_1+\hdots+b_j^mr_m
\right)
=
\sum_{k=1}^ma_j^kb_j^kr_k^2+
\sum_{k\ne l}
a_j^kb_j^l
r_kr_l.
$$
Hence
$$
\sum_{j=1}^n
\Re\sk{u_j}{\omega}\Re\sk{v_j}{\omega}
=
\sum_{k=1}^m\alpha_k r_k^2
+
\sum_{k\ne l}\beta_{k,l}r_kr_l,
$$
where
\begin{equation}
\label{eq: alpha,beta}
\aligned
\alpha_k&=\sum_{j=1}^na_j^kb_j^k\\
\beta_{k,l}&=\sum_{j=1}^n
a_j^kb_j^l.
\endaligned
\end{equation}
This leads to
$$
\aligned
\min_{r_1,\hdots,r_m\in[-1,1] \atop r_1^2+\hdots+r_m^2=1}
\sum_{j=1}^n
\Re\sk{u_j}{\omega}\Re\sk{v_j}{\omega}
& =
\min_{r_1,\hdots,r_m\in[-1,1] \atop r_1^2+\hdots+r_m^2=1}
\Big(
\sum_{k=1}^m\alpha_k r_k^2
+
\sum_{k\ne l}\beta_{k,l}r_kr_l
\Big).
\endaligned
$$

At this point we propose a symmetrization argument. 
We may rename parameters $r_k$, say by shifting $r_k\rightarrow r_{k+d}$, with addition of the indices intended in $\Z_m$, meaning that if $k+d>m$, then we define $r_{k+d}:=r_{k+d-m}$. That is, we may consider
$$
N(d)=N(d,\vec r,u,v,\sigma)
=
\sum_{k=1}^m\alpha_k r_{k+d}^2
+
\sum_{k\ne l}\beta_{k,l}r_{k+d}r_{l+d}.
$$
Since the set of $r_j$'s over which we are taking the minimum is invariant under permutation of indices, we still get the same value. This is to say, 
\begin{equation*}
\min_{\vec r\in[-1,1]^m\atop |\vec r|=1}N(d)
=
\min_{\vec r\in[-1,1]^m\atop |\vec r|=1}N(0)
\hskip 30pt
\text{ for any }d=1,2,\hdots,m-1.
\end{equation*}
Therefore
$$
\min_{\vec r\in[-1,1]^m\atop |\vec r|=1}N(0)
=
\frac{\min_{\vec r} N(0)+\hdots+\min_{\vec r} N(m-1)}m
\leq
\min_{\vec r}\frac{N(0)+\hdots+N(m-1)}{m}.
$$
The point is that in the last average the first sum becomes ``$\vec r$-free'':
$$
\aligned
\frac{N(0)+\hdots+N(m-1)}{m}
&=
\frac1m
	\left(
\sum_{k=1}^m\alpha_k
+
\sum_{k\ne l}\beta_{k,l}\sum_{d=1}^mr_{k+d}r_{l+d}.
	\right)\\
&=
\frac1m
	\bigg(
\sum_{k=1}^m\alpha_k
+
\sum_{k\ne l}r_{k}r_{l}
\underbrace{\sum_{d=1}^m\beta_{k+d,l+d}}_{=\gamma_{k,l}}
	\bigg).
\endaligned
$$
As a consequence we derive
\begin{equation}
\label{eq: N avg est}
\aligned
\min_{r_1,\hdots,r_m\in[-1,1] \atop r_1^2+\hdots+r_m^2=1}
\sum_{j=1}^n
\Re\sk{u_j}{\omega}\Re\sk{v_j}{\omega}
&\leq\min_{\vec r}
\frac{N(0)+\hdots+N(m-1)}{m}\\
& =
\frac1m
	\bigg(
\sum_{k=1}^m\alpha_k
+
\min_{\vec r}
\sum_{k\ne l}\gamma_{k,l}r_{k}r_{l}
	\bigg)\\
&\leq 
\frac1m
\sum_{k=1}^m\alpha_k.
\endaligned
\end{equation}
The last inequality stems from choosing, say, $\vec r=(1,0,\hdots,0)$ which turns all the terms in the double sum inside zero, as $r_1=1$ is in any of those terms multiplied by $r_l$ for $l\ne1$, which is zero.

Now we calculate $\alpha_k$. Recalling \eqref{eq: alpha,beta}, we get to calculate
$$
\aligned
a_j^kb_j^k
&=
\Re\big(u_j^k\overline{\sigma^k}\big)\Re\big(v_j^k\overline{\sigma^k}\big)\\
&=
\frac{
\big(u_j^k\overline{\sigma^k}+\overline{u_j^k}{\sigma^k}\big)
\big(v_j^k\overline{\sigma^k}+\overline{v_j^k}{\sigma^k}\big)
}{4}\\
&=
\frac{
u_j^kv_j^k\big(\overline{\sigma^k}\big)^2
+
\overline{u_j^k}v_j^k
+
u_j^k\overline{v_j^k}
+
\overline{u_j^k}\overline{v_j^k}\big(\sigma^k\big)^2
}{4}\\
&=
\frac{
\Re\big(u_j^k\overline{v_j^k}\big)
+
\Re\big(
\big(\overline{\sigma^k}\big)^2
u_j^kv_j^k
\big)}{2}\,.
\endaligned
$$
This gives
$$
\aligned
\frac1m
\sum_{k=1}^m\alpha_k
&=
\frac1m
\sum_{k=1}^m
\sum_{j=1}^na_j^kb_j^k\\
&=
\frac1m
\sum_{k=1}^m
\sum_{j=1}^n
\frac{
\Re\left(u_j^k\overline{v_j^k}\right)
+
\Re\left(
\big(\overline{\sigma^k}\big)^2
u_j^kv_j^k
\right)}{2}\\
&=
\frac1{2m}
\bigg[
\sum_{j=1}^n
\Re\sk{u_j}{v_j}
+
\sum_{k=1}^m\Re\Big(
\big(\overline{\sigma^k}\big)^2
\sum_{j=1}^n
u_j^kv_j^k
\Big)
\bigg]\,.
\endaligned
$$
Consequently,
\begin{alignat*}{2}
\min_{\sigma^1,\hdots,\sigma^m\in\Sigma^1}
\frac1m
\sum_{k=1}^m\alpha_k
& =
\frac1{2m}
\bigg(
\sum_{j=1}^n
\Re\sk{u_j}{v_j}
-
\sum_{k=1}^m\Big|
\sum_{j=1}^n
u_j^kv_j^k
\Big|
\bigg)\\
& =
\frac1{m}
\sum_{k=1}^m
\frac1{2}
\bigg(
\Re\sk{u^k}{v^k}
-
\Big|
\sum_{j=1}^n
u_j^kv_j^k
\Big|
\bigg)
\,.
\end{alignat*}
In view of \eqref{eq: N avg est} and
\eqref{eq: double min},
this completes the proof.
\end{proof}

\begin{remark}
By recalling \eqref{eq: real vs. complex inner product} again, we obtain that $S=\cV T \cV^{-1}$, where the operator $T=T_{A,\xi}:\C^m\rightarrow\C^m$ is $\R$-linear and maps
$$
Tz=\frac12\sum_{j=1}^n
\left[
\Re\sk{z}{{(A\xi)_j}} \cdot {\xi_j}+\Re\sk{z}{{\xi_j}}\cdot {(A\xi)_j} 
\right]\,.
$$
Let us, for the purpose of this remark, extend the definition of the projections $\pi_\omega$, that appeared in Remark \ref{r: Phi as an operator}, to the case when the symbol is an arbitrary (that is, not necessarily unimodular) element of $\R^{2m}$. By repeating the defining formula verbatim, we get
$
\pi_\zeta X:=\Re\sk{X}{\zeta}_{\C^{m}}\cdot\zeta
$
for $\zeta,X\in\C^{m}$. Then we quickly see that
$$
T=\frac12\sum_{j=1}^n
\left(
\pi_{(A\xi)_j}+\pi_{\xi_j}-\pi_{(A\xi-\xi)_j}
\right).
$$
\end{remark}


\subsection{Comparison with the Dindoš--Li--Pipher condition}
\label{s: comparison}

Here we show that the two notions of $p$-ellipticity for systems, the one from \cite{DiLiPi20} and the one from this paper, coincide up to some factors depending on $p$ only. 

In \cite[(17)]{DiLiPi20}, Dindo\v{s}, Li and Pipher introduced the following condition they called $p$-ellipticity; their terminology followed the one we introduced in \cite{CD-DivForm}, where we treated the scalar case.
\medskip

{\it For all $\xi=(\xi_h^\alpha)\in\C^{n\times m}$, $\omega\in\C^m$ with $|\omega|=1$ and a.e. $x\in\Omega$,
\begin{equation}
\label{eq: p-ell DiLiPi}
\Re\sk{A(x)\left(\xi-\left(1-2/p\right)\xi(\omega)\right)}{\xi+\left(1-2/p\right)\xi(\omega)}\geq C|\xi|^2\,.
\end{equation}}
\noindent
Here $\xi(\omega)$ is $\Phi_\omega\xi$ in our notation.

\begin{proposition}
\label{c: DiLiPi}
The optimal $C$ in \eqref{eq: p-ell DiLiPi} is comparable to $\lambda_p(A)$, that appears in \eqref{eq: lambda_p}, by multiplicative constants depending on $p$ only. 

Consequently, the $p$-ellipticity condition \cite[(17)]{DiLiPi20} by Dindoš, Li and Pipher, reproduced in \eqref{eq: p-ell DiLiPi}, is equivalent to the $p$-ellipticity introduced in \eqref{eq: p-ell}.
\end{proposition}

\begin{proof}
Let $\xi\in\C_n^m$, $\omega\in\Sigma^m$ and $p>1$. Denote
$$
\aligned
X & =\xi-(1-2/p)\Phi_\omega\xi\\
Y & =\xi+(1-2/p)\Phi_\omega\xi.
\endaligned
$$
We can rewrite this as
$$
\aligned
X&=\cJ_{1+2/p,\omega}\xi\\
Y&=\cJ_{3-2/p,\omega}\xi.
\endaligned
$$
Now Corollary \ref{c: JoJ} implies
$$
\aligned
Y&=\cJ_{p,\omega}X\\
\xi&=\cJ_{1+p/2,\omega}X.
\endaligned
$$
Consequently, the condition \eqref{eq: p-ell DiLiPi} reads
\begin{equation}
\label{eq: p-ell DiLiPi via J_p}
\Re\sk{A(x)X}{\cJ_{p,\omega}X}\geq C|\cJ_{1+p/2,\omega}X|^2\,.
\end{equation}
By Corollary \ref{c: norm of J}, 
$|\cJ_{1+p/2,\omega}X|\sim_p|X|$.
Hence the inequality \eqref{eq: p-ell DiLiPi via J_p} is equivalent (up to some factor depending on $p$ only) to the one obtained by replacing $\cJ_{1+p/2,\omega}X$ on the right-hand side with $X$, which is exactly our $p$-ellipticity condition $\lambda_p(A)>0$.
\end{proof}

\begin{remark}
Apart from what was highlighted in Remark \ref{r: Phi as an operator},
perhaps the main difference between our formulation of the $p$-ellipticity for systems \eqref{eq: p-ell} and the one by Dindo\v{s}--Li--Pipher \eqref{eq: p-ell DiLiPi} is that in their case, dependence on $p$ and $\omega$ appears in both factors of the inner product, while in our case it only appears in the right one, as in the scalar case. 
(Or else, cf. \eqref{eq: Deltap bis}, $\omega$ appears only on the left and $p$ only on the right.)
In fact, when formally differentiating the $L^p$ norm of the semigroup, which was the basis for our definition (cf. Section \ref{s: rootsyst}), we necessarily end up with dependence 
on $p$ 
in one factor only, just due to the chain rule.
\end{remark}

\section{Monotonicity of $p$-ellipticity}

\label{s: monotonellipticity}

Recall the notation \eqref{eq: Delta vs. lambda}.
In \cite[Corollary 5.16]{CD-DivForm} we proved that $p\mapsto\Delta_p(A)$ was a decreasing function on $[2,\infty)$. In particular, that implied an important ``chain-type'' property of  the $p$-ellipticity classes:
\begin{equation*}
\aligned
\cA_{p_1}(\Omega)
& \subset \cA_{p_2}(\Omega) \hskip 20pt \text{ for } 1<p_1<p_2\leq 2\\
\cA_{p_1}(\Omega)
& \supset \cA_{p_2}(\Omega) \hskip 20pt \text{ for } 2\leq p_1<p_2<\infty.
\endaligned
\end{equation*}
We want to investigate whether the same property holds for systems. That is, whether we have, for any $m\in\N$,
\begin{equation}
\label{eq: A_p inclusions}
\aligned
\cA_{p_1}^m(\Omega)
& \subset \cA_{p_2}^m(\Omega) \hskip 20pt \text{ for } 1<p_1<p_2\leq 2\\
\cA_{p_1}^m(\Omega)
& \supset \cA_{p_2}^m(\Omega) \hskip 20pt \text{ for } 2\leq p_1<p_2<\infty.
\endaligned
\end{equation}

In \cite[Theorem 1.1 {\it (v)}]{DiLiPi20}, a weaker substitute for \eqref{eq: A_p inclusions} was proved. Namely, the authors showed that if $1<p_1<p_2<\infty$ then 
$
r\in(p_1,p_2) 
\Rightarrow
\cA^m_{p_1} 
\cap
\cA^m_{p_2} 
\subset 
\cA^m_{r} 
$.
Here we fully recover \cite[Corollary 5.16]{CD-DivForm} in the ``vector'' context:

\begin{theorem}
\label{t: A_p inclusions}
For every $m\in\N$, the inclusions \eqref{eq: A_p inclusions} hold.
In particular, for any $p\in[1,\infty)$ we have
$\cA_p^m\subset \cA_2^m$.
\end{theorem}

The theorem will immediately follow from the proposition below, once we 
recall that, by definition,
$
\cA_p^m(\Omega)
=
\mn{A\text{ bounded tensor}}{\lambda_p(A)>0}.
$

\begin{proposition}
Let $p\in(1,\infty)$. For any bounded tensor $A$ the following holds:
\begin{enumerate}[(i)]
\item
If $1< p<p+\e\leq2$ then
$$
\lambda_{p+\e}(A)
\geq
\min
	\left\{
\lambda_p(A),\frac{p+\e-1}{p-1}\,\lambda_p(A)
	\right\}.
$$
\item
If $2\leq p-\e<p<\infty$ then
$$
\lambda_{p-\e}(A)
\geq
\min
	\left\{
\lambda_p(A),\frac{p-\e-1}{p-1}\,\lambda_p(A)
	\right\}.
$$
\end{enumerate}

\end{proposition}

\begin{proof}
Let $A$ be a bounded tensor. We may actually assume that $A$ is constant, i.e., that it does not depend on $x$. Let  $S=S_{A,\xi}$ be the matrix that appears in Theorem \ref{t: S}. 
Denote by 
$
\mu_1(\xi)\leq\hdots\leq\mu_{2m}(\xi)
$
its eigenvalues. 
In the continuation we will for simplicity often write just $\mu_j$, but will keep in mind that in fact $\mu_j=\mu_j(\xi)$.
It follows from \eqref{trace formula} and \eqref{eq: real vs. complex inner product} that
\begin{equation*}
\Re\sk{A\xi}{\xi}=
\mu_1(\xi)+\hdots+\mu_{2m}(\xi).
\end{equation*}
Also, by Theorem \ref{t: S} we have
\begin{equation}
\label{eq: Theta vs. mu}
\aligned
\Theta_A^-(\xi) & = \mu_{1}(\xi)\\
\Theta_A^+(\xi) & = \mu_{2m}(\xi).
\endaligned
\end{equation}
Now \eqref{eq: lambda via Theta} implies the following expression for $\lambda_s(A)$, on which our proof rests:
\begin{equation}
\label{eq: lambda via mu}
\lambda_s(A)
=
\inf_{\xi\in\Sigma_n^m}
\left\{
\begin{array}{lcl}
(\mu_1+\hdots+\mu_{2m-1})+(s-1)\mu_{2m} & ; & s\leq 2\\
(s-1)\mu_{1}+(\mu_2+\hdots+\mu_{2m}) & ; & s\geq 2.
\end{array}
\right.
\end{equation}

{\it (i)} First suppose that \framebox{$1< p<p+\e\leq2$}. Take any $\xi\in\Sigma_n^m$.

\begin{itemize}
\medskip
\item
If \underline{$\mu_{2m}(\xi)\geq0$}, then 
$$
\aligned
(\mu_1+\hdots+\mu_{2m-1})+(p+\e-1)\mu_{2m}
&\geq
(\mu_1+\hdots+\mu_{2m-1})+(p-1)\mu_{2m}\\
&\geq\lambda_p(A).
\endaligned
$$

\medskip
\item
If \underline{$\mu_{2m}(\xi)\leq0$}, then $\mu_j\leq0$ for $j\in\{1,\hdots,2m\}$ and thus $\mu_1+\hdots+\mu_{2m-1}\leq0$. Consequently,
$$
\aligned
(\mu_1+\hdots+& \mu_{2m-1}) + (p+\e-1)\mu_{2m}\\
& \geq 
\frac{p+\e-1}{p-1}\,
(\mu_1+\hdots+\mu_{2m-1})+(p+\e-1)\mu_{2m}\\
& = 
\frac{p+\e-1}{p-1}\,
\left[
(\mu_1+\hdots+\mu_{2m-1})+(p-1)\mu_{2m}
\right]
\\
& \geq 
\frac{p+\e-1}{p-1}\,
\lambda_p(A).
\endaligned
$$ 
\end{itemize}

We proved that for any $\xi\in\Sigma_n^m$ we have
$$
(\mu_1+\hdots+ \mu_{2m-1}) + (p+\e-1)\mu_{2m}
\geq
\min
	\left\{
\lambda_p(A),\frac{p+\e-1}{p-1}\,\lambda_p(A)
	\right\}.
$$
Using 
\eqref{eq: lambda via mu} settles the case {\it (i)}.

\bigskip
{\it (ii)} Now suppose that \framebox{$2\leq p-\e<p<\infty$}. Again take any $\xi\in\Sigma_n^m$.

\begin{itemize}
\medskip
\item
If \underline{$\mu_{1}(\xi)\geq0$}, then $\mu_j\geq0$ for $j\in\{1,\hdots,2m\}$ and thus $\mu_2+\hdots+\mu_{2m}\geq0$. Consequently,
$$
\aligned
(p-1-\e)\mu_{1}+(\mu_2+\hdots+\mu_{2m})
& \geq 
\frac{p-\e-1}{p-1}\,
\left[
(p-1)\mu_{1}+(\mu_2+\hdots+\mu_{2m})
\right]
\\
& \geq 
\frac{p-\e-1}{p-1}\,
\lambda_p(A).
\endaligned
$$ 

\medskip
\item
If \underline{$\mu_{1}(\xi)\leq0$}, then 
$$
\aligned
(p-1-\e)\mu_{1}+(\mu_2+\hdots+\mu_{2m})
& \geq 
(p-1)\mu_{1}+(\mu_2+\hdots+\mu_{2m})
\\
& \geq 
\lambda_p(A).
\endaligned
$$ 
\end{itemize}
We proved that for any $\xi\in\Sigma_n^m$ we have
$$
(p-1-\e)\mu_{1}+(\mu_2+\hdots+\mu_{2m})
\geq
\min
	\left\{
\lambda_p(A),\frac{p-\e-1}{p-1}\,\lambda_p(A)
	\right\}.
$$
Using \eqref{eq: lambda via mu} proves the proposition.
\end{proof}

\begin{theorem}
\label{t: lambda is Lipschitz}
Given any bounded tensor $A$, the function $p\mapsto \lambda_p(A)$ is Lipschitz continuous on $(1,\infty)$.
\end{theorem}

\begin{proof}
From Theorem \ref{t: S} (last part) and Proposition \ref{p: Phi properties} {\it \ref{p: Phi prop 2}} we infer
$$
\left|\Theta^\pm(x,\xi)\right|
\leq\Lambda(A)
|\xi|^2
\hskip 30pt
\forall\xi\in\C_n^m.
$$
Recall \eqref{eq: Theta vs. mu} and \eqref{eq: lambda via mu}. 
If $p,r\in[1,2]$ then, for $\xi\in\Sigma_n^m$,
$$
\aligned
(\mu_1+\hdots+\mu_{2m-1})+(p-1)\mu_{2m}
&=
(\mu_1+\hdots+\mu_{2m-1})+(r-1)\mu_{2m}+(p-r)\mu_{2m}\\
&\geq\lambda_r(A)-|p-r|\Lambda(A).
\endaligned
$$
Taking infimum over $\xi\in\Sigma_n^m$ gives
$
\lambda_p(A)\geq\lambda_r(A)-|p-r|\Lambda(A).
$
By reversing the order of $p,r$ we prove 
$$
|\lambda_p(A)-\lambda_r(A)|\leq\Lambda(A)|p-r|.
$$

Similarly if $p,r\in[2,\infty)$. So $p\mapsto\lambda_p(A)$ is Lipschitz separately on $[1,2]$ and $[2,\infty)$, hence it is Lipschitz on $(1,\infty)$.
\end{proof}

\begin{corollary}
\label{eq: Delta is Lipschitz}
Given any bounded tensor $A$, 
the function $p\mapsto \Delta_p(A)$ is Lipschitz continuous on $(1,\infty)$.
\end{corollary}

\begin{proof}
Denote $f(p)=\lambda_p(A)$ and 
$$
g(p)=\frac{\Delta_p(A)}{2}=\frac{f(p)}{p}\,.
$$
we used \eqref{eq: Delta vs. lambda} for the last identity. 
Then we have, for $p,s\geq1$,
$$
g(p)-g(s)=
\frac1p
\left[
f(p)-f(s)-\frac{f(s)}{s}(p-s)
\right].
$$
Consequently,
$$
|g(p)-g(s)|
\leq
|f(p)-f(s)|+\frac{|f(s)|}{s}|p-s|.
$$
By \eqref{eq: lambda_p via J_p} and Corollary  \ref{c: norm of J}
$$
\frac{|f(s)|}{s}\leq \Lambda(A)\max\left\{\frac1s,1-\frac1s\right\}\leq \Lambda(A).
$$
Finally, use that $f$ is Lipschitz continuous, which was proven in Theorem \ref{t: lambda is Lipschitz}.
\end{proof}

\subsection{Operators $\Theta_A^\pm$ revisited}

According to \eqref{eq: lambda via Theta}, {\it if} for a.e. $x\in\Omega$ we have
\begin{equation}
\label{eq: Theta vs. 0}
\Theta^-(x,\xi)
\leq
0 
\leq
\Theta^+(x,\xi)
\hskip 30pt
\forall\xi\in\C_n^m,
\end{equation}
then $p\mapsto\lambda_p(A)$ is increasing on $(1,2)$ and decreasing on $(2,\infty)$, which implies \eqref{eq: A_p inclusions}. 
When $m>n$, by Theorem \ref{t: S} and Remark \ref{r: range of M}, the condition \eqref{eq: Theta vs. 0} holds. 
However, in general the first inequality in \eqref{eq: Theta vs. 0} does {\it not} hold, while the last does hold for any elliptic $A$, by Corollary \ref{c: A form vs. Theta}.
 
\begin{example}
\label{ex: Theta & m=1}
Consider \framebox{$m=1<n$}. 
For every $\xi\in\C_n$, \eqref{eq: T for m=1} gives
\begin{equation}
\label{eq: Theta & m=1}
\Theta^\pm(\xi)
=\frac{\Re\sk{A\xi}{\xi}\pm|\sk{A\xi}{\bar\xi}|}2\,.
\end{equation}
However, $\Theta_A^-$ can be positive even if $A$ is elliptic rather than merely bounded: take, for example, rotation matrices $A=e^{i\f}I_n$ with $-\pi/2<\f<\pi/2$ (which ensures ellipticity), where we get
$$
\Theta^-(\xi)
=\frac{\cos\f|\xi|^2-|\sk{\xi}{\bar\xi}|}2\,.
$$
Thus for choosing $\xi\ne0$ such that, say, $\sk{\xi}{\bar\xi}=0$ (which is possible since $n>1$), we get $\Theta^-(\xi)=\cos\f|\xi|^2/2>0$.
\end{example}


\section{Relation with real elliptic matrices}
\label{s: real}

The starting point of this section is the following important property of $p$-ellipticity, which was in the scalar case established in \cite{CD-DivForm}:
\begin{theorem}\cite[p.3179]{CD-DivForm}
\label{t: scalar real via p-ell}
An elliptic matrix $A$ is
$$
\text{real } \Longleftrightarrow \text{ $p$-elliptic for all }p>1. 
$$
\end{theorem}
Afterwards, 
the ``$\Leftarrow$'' part was 
enhanced by
an elegant quantitative estimate due to ter Elst, Haller-Dintelmann, Rehberg and Tolksdorf \cite{EHRT21}:

\begin{theorem}\cite[Lemma 2.15]{EHRT21}
Take $p\in(1,\infty)\backslash\{2\}$.
Let $A=U+iV$ be $p$-elliptic. Then for almost every $x\in\Omega$ we have
$$
\nor{V(x)}_{\cL(\R^n)}\leq 2\frac{\sqrt{p-1}}{|p-2|}\nor{U(x)}_{\cL(\R^n)},
$$
where $\nor{\Lambda}_{\cL(\R^n)}$ denotes the operator norm of the linear map $\Lambda:\R^n\rightarrow\R^n$.
\end{theorem}

\begin{remark}
Let us make a few observations regarding the above theorem. 
\begin{enumerate}
\item
Write $1-2/p=\cos\f$ for some $\f\in\R$. Then $2\sqrt{p-1}/|p-2|=|\tan\f|$.

\item
The theorem was in \cite{EHRT21} formulated for $p>2$ only. The estimate for $1<p<2$ follows from the one for $p>2$ by using that, in the scalar case at least, $\cA^m_p=\cA^m_q$ with $1/p+1/q=1$ and noticing that
$$
2\frac{\sqrt{p-1}}{|p-2|}
=
2\frac{\sqrt{q-1}}{|q-2|}.
$$

\item
The case $p=2$ is different. That is, by assuming (2-)ellipticity only, one cannot prove a general estimate of the type $\nor{V}\leqsim\nor{U}$. In order to see it, consider the rotation matrices $e^{i\phi}I$ with $\phi\in(-\pi/2,\pi/2)$.
\end{enumerate}
\end{remark}

We prove a variant of this theorem for systems, as summarized in Theorem \ref{t: EHRT vect}. 
As a consequence, we establish a vector version of Theorem \ref{t: scalar real via p-ell}; see Corollary \ref{c: realness}.

\medskip
First we need an auxiliary technical result. 

\begin{lemma}
\label{l: R=C}
Let $x,y\in\R_n^m$. Suppose that 
$$
R:=\max_{\sigma\in S^{m-1}}\Sk{[x,\sigma]}{[y,\sigma]}\geq0.
$$
Then 
\begin{equation}
\label{eq: R=C}
R=\max_{\omega\in \Sigma^{m}}\Sk{\Re [x,\omega]}{\Re [y,\omega]}.
\end{equation}
\end{lemma}

\begin{proof}
Denote the quantity on the right-hand side of \eqref{eq: R=C} by $C$. 

If $\omega\in\Sigma^m$ then for every {\it real} $x,y$ we have
\begin{equation}
\label{eq: Re[] ids}
\Sk{\Re [x,\omega]}{\Re [y,\omega]}
=\Sk{[x,\Re \omega]}{[y,\Re \omega]}
=|\Re\omega|^2\Sk{[x,\sigma]}{[y,\sigma]},
\end{equation}
where $\sigma:=\Re\omega/|\Re\omega|\in S^{m-1}$ if $\Re\omega\ne0$, while if $\Re\omega=0$ then any $\sigma\in S^{m-1}$ fits.
Hence for
$$
\aligned
{\mathcal R}={\mathcal R} (x,y) &:= \Mn{\Sk{[x,\sigma]}{[y,\sigma]}}{\sigma\in S^{m-1}}\\
{\mathcal C}={\mathcal C} (x,y) &:= \Mn{\Sk{\Re [x,\omega]}{\Re [y,\omega]}}{\omega\in \Sigma^{m}}
\endaligned
$$
we get
$
{\mathcal C}=[0,1]\cdot {\mathcal R}.
$
Indeed, ``$\subseteq$'' follows directly from \eqref{eq: Re[] ids}, while in order to confirm ``$\supseteq$'' take any $\tau\in[0,1]$, $\sigma\in S^{m-1}$, write $\tau=\cos^2\f$ with $\f\in[0,\pi/2]$ and define $\omega=e^{i\f}\sigma$. Then 
$\sigma\in\Sigma^m$ and 
$\tau\Sk{[x,\sigma]}{[y,\sigma]}=\Sk{\Re [x,\omega]}{\Re [y,\omega]}$, again by \eqref{eq: Re[] ids}.

This implies 
$$
C=\max{\mathcal C} = \max\left([0,1]\cdot {\mathcal R}\right) = \max{\mathcal R}=R,
$$
where we used the assumption $R\geq0$ for the penultimate identity.
\end{proof}

\begin{theorem}
\label{t: EHRT vect}
Take $m,n\in\N$ and $p>2$. 
Let $A=U+iV$ be a $p$-elliptic tensor acting on $\C_n^m$. 
Then for a.e. $x\in\Omega$ we have
$$
\nor{V(x)}_{\cL(\R_n^m)}\leq 4m\frac{\sqrt{p-1}}{p-2}\nor{U(x)}_{\cL(\R_n^m)}.
$$
\end{theorem}

\begin{proof}
We adapt the approach from \cite{EHRT21} and combine it with our linear algebra estimate \eqref{eq: (u,v) bilateral est complex}. 
We shall use the $p$-ellipticity as expressed in the form of \eqref{eq: lambda_p scalar}.
Write $A=A(x)$. 
For a.e. $x\in\Omega$, the $p$-ellipticity implies, for any $\xi\in\C^m_{n}$ and $\omega\in\Sigma^m$,
\begin{equation}
\label{eq: Delta geq0}
0\leq
\Re\sk{A\xi}{\xi}_{\C_n^m}+(p-2)\sk{\Re[A\xi,\omega]}{\Re[\xi,\omega]}_{\R_n}.
\end{equation}
Write $\xi=u+iv$ and assume that $\omega$ is \underline{real}. Then 
$$
\aligned
\Re\sk{A\xi}{\xi}_{\C_n^m} 
& = \Re\sk{(U+iV)(u+iv)}{u+iv}_{\C_n^m}\\
& = 
\sk{Uu-Vv}{u}_{\R_n^m}
+
\sk{Vu+Uv}{v}_{\R_n^m};\\
\sk{\Re[A\xi,\omega]}{\Re[\xi,\omega]}_{\R_n}
&= 
\sk{[Uu-Vv,\omega]}{[u,\omega]}_{\R_n}.
\endaligned
$$
Thus we transformed \eqref{eq: Delta geq0} into
$$
\aligned
\sk{Vv}{u}_{\R_n^m}
+&(p-2)\sk{[Vv,\omega]}{[u,\omega]}_{\R_n}\\
&\hskip 10pt \leq
\sk{Uu}{u}_{\R_n^m}
+
\sk{Vu+Uv}{v}_{\R_n^m}
+(p-2)\sk{[Uu,\omega]}{[u,\omega]}_{\R_n}.
\endaligned
$$
At this point we take the maximum over all $\omega\in S^{m-1}$. 

We estimate the left-hand side from below by using Corollary \ref{c: bilateral est}, that is, its variant formulated in the first inequality of Remark \ref{c: bilateral est complex}.
In order to do so, we have to make sure that the maximum over $\omega\in S^{m-1}$ is the same as the maximum over $\omega\in \Sigma^{m}$. This is provided by Lemma \ref{l: R=C}; its assumption is fulfilled due to the very Remark \ref{c: bilateral est complex}, for if needed, we may replace $u$ by $-u$.

As for the right-hand side, we estimate it by using the Cauchy--Schwarz--Bunyakovsky inequality and its version \eqref{eq: CSB}; in both cases we also use the assumption that $p\geq2$.
This brings us to
$$
\left(1+
\frac{p-2}{2m}
\right)
\sk{Vv}{u}
\leq
\nor{U}
\left[
\left(
p-1 
\right)
|u|^2
+
|v|^2
\right]
+\nor{V}|u||v|.
$$
By choosing $u=Vv/|Vv|$ (that is, maximizing over $u\in\R_n^m$, $|u|=1$) we get 
$$
\left(1+\frac{p-2}{2m}\right)
|Vv|
\leq
\nor{U}
\left[
(p-1)
+
|v|^2
\right]
+\nor{V}|v|.
$$
Replacing $v$ by $tv$, dividing by $t$ and minimizing the right-hand side over $t>0$ gives
$$
\aligned
\left(1+\frac{p-2}{2m}\right)
|Vv|
&\leq
\big(
2\sqrt{p-1}\nor{U}+\nor{V}
\big)
|v|.
\endaligned
$$
This means
$$
\left(1+\frac{p-2}{2m}\right)
\nor{V}
\leq
2\sqrt{p-1}\nor{U}+\nor{V},
$$
which is equivalent to the desired inequality. 
\end{proof}

\begin{corollary}
\label{c: realness}
Any elliptic tensor $A$ is
\begin{equation}
\label{eq: p-ell vs. real}
\text{ $p$-elliptic for all }
p\in(1,2]
\text{ or all }
p\in[2,\infty) 
\hskip 10pt
\Rightarrow 
\hskip 10pt
\text{real}. 
\end{equation}
As for the implication in the other direction:
\begin{enumerate}[(i)]
\item
if $m=1$, it holds;
\item
if $m>1$, then for every $p\ne2$ there exists a real elliptic tensor acting on $\C_n^m$ that is not $p$-elliptic.
\end{enumerate}
Thus in general, ``$\Leftarrow$'' holds in \eqref{eq: p-ell vs. real} if and only if $m=1$.
\end{corollary}

\begin{proof}
($\Rightarrow$):
Assume that $A$ is $p$-elliptic for all $p\geq2$. Now apply Theorem \ref{t: EHRT vect} and send $p\rightarrow\infty$.
If, however, $A$ is $p$-elliptic for all $1<p\leq2$, then $A^*$ is $q$-elliptic for all $q\geq2$ (Proposition \ref{p: p-ell of A vs. A*}), thus by the part just proven we have that $A^*$ is real, i.e., $A$ is real.

\medskip
{\it (i)}
We already know this from \cite{CD-DivForm}, cf. Theorem \ref{t: scalar real via p-ell}; 

\medskip
{\it (ii)}
If $m=2$, adequate examples are given in Theorem \ref{t: Delta of V_theta} and Theorem \ref{t: aI+bI}.
For $m>2$ use these examples and add $m-2$ identity blocks by direct sum.
\end{proof}

\subsection{Summary}
Here we summarize the findings from the Sections \ref{s: monotonellipticity} and \ref{s: real}.

Recall from Section \ref{s: p-ell} the notation
$$
\cA_p^m(\Omega):=\{\text{all }p\text{-elliptic tensors on }\Omega \text{ acting on }\C_n^m\}.
$$
We have, for $2\leq p\leq r$, 
\renewcommand{\arraystretch}{1.6}
$$
\begin{array}{rcc}
\{\text{elliptic tensors on }\Omega\} & = & \cA_2^m(\Omega)\\
&& \text{\rotatebox{90}{$\subset$}}\\ 
&& \cA_p^m(\Omega)\\
&& \text{\rotatebox{90}{$\subset$}}\\
&& \cA_r^m(\Omega)\\
&& \text{\rotatebox{90}{$\subset$}}\\
\{ \text{real elliptic tensors on }\Omega\} & \supset & 
{\displaystyle\bigcap_{s\in[2,\infty)}\cA_s^m(\Omega)}
\end{array}
$$
The inclusion in the bottom row is an equality if and only if $m=1$.
This, however, means that in the vector case ($m>1$), assuming $p$-ellipticity for $p\ne2$ may bring something new even for real tensors. In contrast, in the scalar case ($m=1$), assuming $p$-ellipticity for $p\ne2$ is for real tensors (in case of $m=1$, matrices) equivalent to merely assuming ($2$-)ellipticity.

Moreover, $\cA_1^m(\Omega)=\cA_\infty^m(\Omega)=\emptyset$.

\section{The decoupling property}

Recall the notation \eqref{eq: Delta vs. lambda}. 
It was proven in \cite[Proposition 5.15]{CD-DivForm} that when $m=1$ and $A\in\cA(\Omega)$,
\begin{equation}
\label{eq: scalar decoupling}
\Delta_p(A)
\sim\mu(A)-|1-2/p|,
\end{equation}
with the implied constants depending on $A$ only, where $\mu(A)$ was given by
\begin{equation}
\label{eq: scalar mu}
\mu(A):={\rm ess}\inf\Re\frac{\sk{A(x)\xi}{\xi}}{\left|\sk{A(x)\xi}{\bar\xi}\right|}\,.
\end{equation}
Following Egert \cite[p.707]{E20}, we call \eqref{eq: scalar decoupling} the {\it decoupling property} (or {\it decoupling formula}, for it gives the endpoints of the $p$-ellipticity interval pertaining to $A$, simply by solving for $p$ the equation $|1-2/p|=\mu(A)$). 

In \cite[Theorem 1.1 {\it (v)}]{DiLiPi20}, a weaker substitute for \cite[Proposition 5.15]{CD-DivForm} was proved. Namely, the authors showed that the set of all $p$'s such that $A$ is $p$-elliptic forms an open interval. 
This also follows from Theorem \ref{t: lambda is Lipschitz}.
Here we fully recover \cite[Proposition 5.15]{CD-DivForm} in the ``vector'' context, by providing for any $m\in\N$ the formula for the endpoints of that interval.


We invoke another set of auxiliary operators, closely related to $\cJ_{p,\omega}$, that have their origin in \cite{CD-DivForm}. 

\subsubsection{Operators $\Xi^\pm$}

\ \medskip

Assume the notation from Section \ref{s: aux ops}. 
By \eqref{eq: lambda_p via J_p}, 
\eqref{eq: Delta vs. lambda}, \eqref{eq: I vs. J vect}
and Lemma \ref{l: I,J vect vs. I,J scal},
\begin{equation}
\label{eq: Delta_p expanded}
\aligned
\Delta_p(A)
&=
\underset{x\in\Omega}{{\rm ess}\inf}
\min_{\xi\in\Sigma^m_n}
\min_{\omega\in\Sigma^m}
\Re\sk{A(x)\xi}{\cI_{p,\omega}\xi}_{\C^m_n}\\
&=
\underset{x\in\Omega}{{\rm ess}\inf}
\min_{\xi\in\Sigma^m_n}
\min_{\omega\in\Sigma^m}
\Re\sk{A(x)\xi}{\xi-(1-2/p)\cJ_{0,\omega}\xi}_{\C^m_n}\\
&=
\underset{x\in\Omega}{{\rm ess}\inf}
\min_{\xi\in\Sigma^m_n}
\big(
\Re\sk{A(x)\xi}{\xi}-
\max_{\omega\in\Sigma^m}
(1-2/p)\Re\sk{A(x)\xi}{\cJ_{0,\omega}\xi}_{\C^m_n}
\big).
\endaligned
\end{equation}
This makes us introduce
\begin{equation}
\label{eq: Xi^pm}
\aligned
	\Xi_A^-(x,\xi)
		&=\min_{\omega\in\Sigma^m}\Re\sk{A(x)\xi}{\cJ_{0,\omega}\xi}_{\C^m_n}\\
	\Xi_A^+(x,\xi)
		&=\max_{\omega\in\Sigma^m}\Re\sk{A(x)\xi}{\cJ_{0,\omega}\xi}_{\C^m_n},
\endaligned
\end{equation}
for then we get the following counterpart to \eqref{eq: lambda via Theta}:
\begin{equation}
\label{eq: Delta via Xi}
\Delta_p(A)=
\left\{
	\begin{array}{rcl}
		{\displaystyle
			\underset{x\in\Omega}{{\rm ess}\inf}
				\min_{\xi\in\Sigma^m_n}
					\left(
						\Re\sk{A(x)\xi}{\xi}-(1-2/p)\Xi^-_A(x,\xi)
					\right)
		} 
& ; & \text{ if }p\leq2\\
		{\displaystyle
			\underset{x\in\Omega}{{\rm ess}\inf}
				\min_{\xi\in\Sigma^m_n}
					\left(
						\Re\sk{A(x)\xi}{\xi}-(1-2/p)\Xi^+_A(x,\xi)
					\right)
		} 
& ; & \text{ if }p\geq2.
	\end{array}
\right.
\end{equation}

\noindent
Note that, by \eqref{eq: J_{0,omega}}, 
\begin{equation}
\label{eq: Xi via Theta}
\Xi_A^\pm(\xi)=\Re\sk{A\xi}{\xi}_{\C^m_n}-2\Theta_A^\mp(\xi).
\end{equation}
From \eqref{eq: Theta & m=1} and \eqref{eq: Xi via Theta} we get
\begin{equation}
\label{eq: Xi at m=1}
\Xi_A^\pm(x,\xi)=\pm|\sk{A(x)\xi}{\bar\xi}|
\hskip 30pt
\text{ when }m=1,
\end{equation}
The identity \eqref{eq: Xi at m=1} shows that the formula \eqref{eq: Delta via Xi} is a generalization of \cite[(5.16)]{CD-DivForm}.


Indeed, with $A\in\cA^m(\Omega)$ and $\Xi^\pm$ as in \eqref{eq: Xi^pm}, we define 
\begin{equation}
\label{eq: mu vector}
\mu^\pm(A):={\rm ess}\inf\Re\frac{\sk{A(x)\xi}{\xi}}{|\Xi^\pm_A(x,\xi)|}\,.
\end{equation}
That is,
\begin{equation*}
\mu^+(A)
={\rm ess}\inf\frac{\hskip -16.4pt \Re\sk{A(x)\xi}{\xi}}{\displaystyle{\max_{\omega\in\Sigma^m}}\Re\sk{A(x)\xi}{\xi-2\Phi_{\omega}\xi}}
={\rm ess}\inf\frac{\hskip -51.5pt \Re\sk{A(x)\xi}{\xi}}{\Re\sk{A(x)\xi}{\xi}-2\Theta^-_A(\xi)}
\end{equation*}
and similarly for $\mu^-(A)$. 
Like in the scalar case \eqref{eq: scalar mu}, we intend by \eqref{eq: mu vector} that $\mu^\pm(A)={\rm ess}\inf\f$, where the function $\f:\Omega\rightarrow\R\cup\{+\infty\}$ is defined by
$$
\f(x)=\inf\Re\frac{\sk{A(x)\xi}{\xi}}{|\Xi^\pm_A(x,\xi)|}\,
$$
and the above infimum runs over all $\xi\in\C^m_n$ for which $\Xi^\pm_A(x,\xi)\ne0$.
From \eqref{eq: Xi at m=1} we see that \eqref{eq: mu vector} is indeed a generalization of \eqref{eq: scalar mu}.

\medskip

We continue by an explicit comparison between $\lambda_p(A)$ and $\mu^\pm(A)$. 
This extends, with the same constants, the result proven for $m=1$ in \cite[Proposition 5.15]{CD-DivForm}.

\begin{theorem}
\label{t: decoupling}
Suppose that $A\in\cA(\Omega)$ and $p\in[2,\infty)$. Then $\Delta_p(A)=\lambda_p(A)=0$ if and only if $1-2/p=\mu^+(A)$. 
If these equalities are not satisfied then
\begin{equation}
\label{eq: vector decoupling}
\frac{\lambda(A)}{\mu^+(A)}
\leq
\frac{\Delta_p(A)}{\mu^+(A)-|1-2/p|}
\leq
\Lambda(A).
\end{equation}
Therefore
$$
\Delta_p(A)\sim\mu^+(A)-|1-2/p|.
$$ 
If $p\in(1,2]$, the same holds with $\mu^+(A)$ replaced with $\mu^-(A)$.
\end{theorem}

\begin{proof}
The proof is based on the formula \eqref{eq: Delta via Xi}. There, one factors out $\Re\!\sk{A(x)\xi}{\xi}$ or $\Xi_A^\pm(x,\xi)$, depending on the part of \eqref{eq: vector decoupling} which is being proven, recalls ellipticity \& boundedness of $A$ and \eqref{eq: mu vector}, and applies properties of essential infimum. 
The complete proof is rather elementary yet tedious, therefore we leave it out.
\end{proof}

\noindent
In the scalar case ($m=1$), calculating $\lambda_p(A)$, or verifying $p$-ellipticity, 
is for any $p>1$ in terms of difficulty comparable to calculating $\lambda_2(A)$, or verifying ($2$-)ellipticity. In the vector case ($m>1$) we saw that this is not the case. Therefore it would be of use to come up with a simpler condition sufficient for $p$-ellipticity. Theorem \ref{t: decoupling} enables that, as we show next.

\begin{corollary}
Let $A$ be an elliptic tensor and $p>1$. Assume that
\begin{equation}
\label{eq: p-ell sufficient cond}
|1-2/p|<\frac{\lambda(A)}{\Lambda(A)}.
\end{equation}
Then $A$ is $p$-elliptic.
\end{corollary}

\begin{proof}
Let us first show that for any elliptic tensor $A$ we have
\begin{equation}
\label{eq: mu bilateral est}
\frac{\lambda(A)}{\Lambda(A)}\leq\mu^\pm(A)\leq1.
\end{equation}
When $m=1$, these inequalities were already established in \cite[(1.11)]{CD-DivForm}. Thus let us now show that they hold when $m>1$.

From \eqref{eq: Xi^pm} and \eqref{eq: J_0 isometry} it follows that for almost every $x\in\Omega$ we have
$$
\left |
\Xi_A^\pm(x,\xi)
\right |
\leq \Lambda(A)|\xi|^2.
$$
Together with the ellipticity assumption on $A$, this implies the first inequality in \eqref{eq: mu bilateral est}.

As for the second one, choose $\omega=(1,0,\hdots,0)$ and $\xi$ with $\xi_h=(0,*,\hdots,*)$ for all $h$; such a choice of nonzero vectors is possible since $m>1$. Then $\Phi_\omega \xi=0$, hence $\cJ_{0,\omega}\xi=\xi$, hence $\mu^+(A)\leq1$. Similarly, by choosing $\xi,\omega$ such that $\Phi_\omega\xi=\xi$, for example, by choosing $\xi_j=\omega$ for any $j=1,\hdots,n$, we eventually see that $\mu^-(A)\leq1$.

Now the assumption \eqref{eq: p-ell sufficient cond} and the first inequality from \eqref{eq: mu bilateral est} merge into
$$
|1-2/p|<\mu^{\pm}(A),
$$
which by Theorem \ref{t: decoupling} means that $A$ is $p$-elliptic.
\end{proof}

\begin{remark}
We saw in Example \ref{ex: Theta & m=1} that it is possible for some $\xi$ to have $\Theta_A^-(\xi)>0$ and
$$
\frac{\hskip -51.5pt \Re\sk{A(x)\xi}{\xi}}{\Re\sk{A(x)\xi}{\xi}-2\Theta^-_A(\xi)}>1.
$$
This is however not in contradiction with \eqref{eq: mu bilateral est}, since in order to get $\mu^+(A)$ we are supposed to take the {\it infimum} of the quotients above over all $\xi$.
\end{remark}

\begin{remark}
From \eqref{eq: Delta_p expanded} and \eqref{eq: J_{0,omega}} we see that it makes sense to define, for any bounded tensor $A$, 
\begin{equation*}
\Delta_\infty(A)
=
\underset{x\in\Omega}{{\rm ess}\inf}
\min_{
\xi\in\Sigma^m_n}
\min_{
\omega\in\Sigma^m}
\Re\sk{A(x)\xi}{2\Phi_\omega\xi}_{\C^m_n}.
\end{equation*}
It is not difficult to see that with this definition we have
$$
\Delta_\infty(A)=\lim_{p\rightarrow\infty}\Delta_p(A).
$$
Indeed, rewrite the middle expression in \eqref{eq: Delta_p expanded} as 
\begin{equation*}
\Delta_p(A)
=
\underset{x\in\Omega}{{\rm ess}\inf}
\min_{
\xi\in\Sigma^m_n}
\min_{
\omega\in\Sigma^m}
\left[
\Re\sk{A(x)\xi}{2\Phi_\omega\xi}_{\C^m_n}
+\frac2p\,\Re\sk{A(x)\xi}{\cJ_{0,\omega}\xi}_{\C^m_n}
\right]
\end{equation*}
and take into account the uniform bound $|\sk{A(x)\xi}{\cJ_{0,\omega}\xi}|\leq\Lambda(A)$, cf. \eqref{eq: A bddness} and \eqref{eq: J_0 isometry}.

If $A$ is also elliptic, then Theorem \ref{t: decoupling} and \eqref{eq: mu bilateral est} imply that there exists a unique $p_0=p_0(A)\in(2,\infty]$ such that $\Delta_{p_0}(A)=0$; this notation is from \cite{EHRT21}. By the same theorem, one has $\mu^+(A)=1-2/p_0$.
If $p_0<\infty$, we see that the middle term in \eqref{eq: vector decoupling} becomes a ``difference quotient'' for $p>2$:
$$
-p_0\,\frac{\lambda_p(A)-\lambda_{p_0}(A)}{p-p_0}\,.
$$
This makes it natural to ask whether (or when) the limit of the above expression as $p\rightarrow p_0$ exists and, if it does, what it is equal to.

More generally, we saw in Theorem \ref{t: lambda is Lipschitz} that $p\mapsto \lambda_p(A)$ is Lipschitz continuous on $(1,\infty)$. 
By Rademacher's theorem, it is then also differentiable almost everywhere. 
In general, however, it is not differentiable {\it everywhere} on that set. An example is provided already by taking $m=1$ and $A=I$, as in that case $\lambda_p(I)=p-1$ if $1\leq p\leq 2$ and $\lambda_p(I)=1$ if $p\geq 2$.  
However, we do not know whether, or for which tensors, we may still have differentiability away from $p=2$.

The same differentiability questions linger for $p\mapsto \Delta_p(A)$ with $p\in(1,\infty)$, for this map is also Lipschitz continous (Corollary \ref{eq: Delta is Lipschitz}) and the same example of $m=1$, $A=I$ shows that it is not differentiable (precisely) at $p=2$.
\end{remark}

\section{Other properties of $\lambda_p(A)$}

\subsection{$p$-ellipticity of tensors implies $p$-ellipticity of diagonal block entries}

\begin{proposition}
\label{p: lambda_p on diagonal}
Assume the notation from Section \ref{eq: matrices, tensors}. 
Suppose $A$ is a bounded tensor on $\Omega$ acting on $\C_n^m$, and $p>1$. Then: 
\begin{itemize}
\item
for every $\alpha\in\{1,\hdots,m\}$ we have 
\begin{equation}
\label{p: lambda_p diag up}
\lambda_p(A)\leq\lambda_p(A^{\alpha,\alpha}); 
\end{equation}
thus if $A$ is a $p$-elliptic tensor, then $A^{\alpha,\alpha}$ is $p$-elliptic as a complex $\C^{n\times n}$-valued function on $\Omega$;

\item
for every $j\in\{1,\hdots,n\}$ we have 
\begin{equation}
\label{p: lambda_p diag down}
\lambda_p(A)\leq\lambda_p(A_{j,j}); 
\end{equation}
thus if $A$ is a $p$-elliptic tensor, then $A_{j,j}$ is 
$p$-elliptic as a complex $\C^{m\times m}$-valued function on $\Omega$.
\end{itemize}

When $m=1$, the latter means that the scalar functions $a_{jj}$ must map into the open sector 
$S_{\arccos|1-2/p|}$, cf. notation from Section \ref{sec}.
\end{proposition}

\begin{proof}
This part is rather straightforward. 
Take
$$
\aligned
\xi & = (0,\hdots,0,\xi^\alpha\, ,0,\hdots,0)\in\C_n^m\\
\omega & = (0,\hdots,0,\omega^\alpha,0,\hdots,0)\in\Sigma^m.
\endaligned
$$
Here the nonzero terms ($\xi^\alpha\in\C^n$, $\omega^\alpha\in\C$) lie at the $\alpha$-th spot.
Then, by recalling \eqref{eq: quasi inner product - up summation} and the fact that, owing to the particular choice of $\omega$, one has $|\omega^\alpha|=1$,
$$
\aligned
\Re\sk{A\xi}{\cJ_{p,\omega}\xi}
& = 
\Re\sk{A\xi}{\xi+(p-2)\Phi_{\omega}\xi}_{\C_n^m}\\
& =
\sum_{\beta=1}^m\Re\Sk{(A\xi)^{\beta}}{\xi^\beta+(p-2)(\Phi_\omega\xi)^\beta}_{\C_n}\\
& =
\sum_{\beta,\gamma=1}^m\Re\Sk{A^{\beta,\gamma}\xi^\gamma}{\xi^\beta+(p-2)\omega^\beta\Re[\xi,\omega]}\\
& =
\Re\Sk{A^{\alpha,\alpha}\xi^\alpha}{\xi^\alpha+(p-2)\omega^\alpha\Re(\overline{\omega}^\alpha
\xi^\alpha)}\\
& =
\Re\Sk{A^{\alpha,\alpha}(\overline{\omega}^\alpha
\xi^\alpha)}{(\overline{\omega}^\alpha
\xi^\alpha)+(p-2)\Re(\overline{\omega}^\alpha
\xi^\alpha)}.
\endaligned
$$
Write $\eta:=\overline{\omega}^\alpha\xi^\alpha\in\C^n$. 
Thus we proved
\begin{equation}
\label{eq: xi rotation}
\Re\Sk{A^{\alpha,\alpha}\eta}{\eta+(p-2)\Re\eta}
= 
\Re\sk{A\xi}{\cJ_{p,\omega}\xi}
\geq
\lambda_p(A)|\xi|^2
=
\lambda_p(A)|\eta|^2\,.
\end{equation}
Clearly, any $\eta\in\C^n$ can be represented in the above manner, which implies \eqref{p: lambda_p diag up}.

\medskip
Now let us prove the second assertion.
This time we choose
$$
\xi=(0,\hdots,0,\xi_j,0,\hdots,0)\in\C^m\times\hdots\times\C^m
$$
with $\xi_j$ at the $j$-th spot, giving
$$
\aligned
\Re\sk{A\xi}{\cJ_{p,\omega}\xi}
& =
\sum_{k=1}^n\Re\Sk{(A\xi)_k}{\xi_k+(p-2)(\Phi_\omega\xi)_k}_{\C^m}\\
& =
\sum_{k,l=1}^n
\Re\Sk{A_{k,l}\xi_l}{\xi_k+(p-2)\Re\sk{\xi_k}{\omega}_{\C^m}
\omega}_{\C^m}\\
&=
\Re\Sk{A_{j,j}\xi_j}{\xi_j+(p-2)\Re\sk{\xi_j}{\omega}_{\C^m}\omega}_{\C^m}.
\endaligned
$$
We want to find $\omega\in\Sigma^m$ such that
$
\Re\sk{\xi_j}{\omega}_{\C^m}\omega
=
\Re \xi_j.
$
If $\Re \xi_j=0$, any real $\omega$ suffices. 
Otherwise, one checks that the natural candidate, $\omega=\Re \xi_j/|\Re \xi_j|$, is actually adequate.

We thus proved that for every $\sigma\in\C^m\backslash\{0\}$ we have
$$
\Re\Sk{A_{j,j}\sigma}{\sigma+(p-2)\Re\sigma}_{\C^m}
=
\Re\sk{A\xi}{\cJ_{p,\omega}\xi},
$$
where 
$$
\xi:=(0,\hdots,0,\sigma,0,\hdots,0)\in\C^m\times\hdots\times\C^m
$$
with $\sigma$ in the $j$-th place and
$$
\omega :=
\frac{\Re \sigma}{|\Re \sigma|}\in\Sigma^m\,.
$$
Since $|\xi|=|\sigma|$, 
one concludes as in \eqref{eq: xi rotation} that \eqref{p: lambda_p diag down} holds.

The very last assertion follows from a short calculation.
\end{proof}

\begin{corollary}
For any bounded tensor $A$ we have $\lambda_1(A)\leq0$.
\end{corollary}

\begin{proof}
By Proposition \ref{p: lambda_p on diagonal}, it suffices to prove this for {\it matrices} $A$, that is, when $m=1$. In that case we have, for any $\xi\in\R_n$ (that is, when the components of $\xi$ are {\it real}) with $|\xi|=1$, 
\begin{alignat*}{2}
\tag*{\qedhere}
\lambda_1(A)
\leq
\Re\sk{A\xi}{\xi-\Re\xi}
=
\Re\sk{A\xi}{i\Im\xi}
=0.
\end{alignat*}
\end{proof}

\subsection{The $p$-ellipticity of direct sums of elliptic matrices}

\label{s: dir sums}

Here we address a converse to Proposition \ref{p: lambda_p on diagonal}.

Let $A_1,\hdots,A_m\in\C^{n\times n}$. 
Recall the notation from Section \ref{eq: matrices, tensors}. 
Denote
$$
\left(A_1\oplus\hdots\oplus A_m\right)^{\text{up}}
:=
\begin{bmatrix}
A_1 & & \\
& \ddots & \\
& & A_m
\end{bmatrix}.
$$
That is, for any $\xi\in\C_n^m$ as in \eqref{{eq: xi components}} we define
$$
(A_1\oplus\hdots\oplus A_m)(\xi^1,\hdots,\xi^m)=(A_1\xi^1,\hdots,A_m\xi^m).
$$
In particular, for $A\in\C^{n\times n}$ we denote
$$
A^{\oplus m}
:=\underbrace{A\oplus\hdots\oplus A}_{m}.
$$
We consider the following question:

\medskip
\begin{center}
\framebox{
If $A_1,\hdots,A_m$ are $p$-elliptic matrices, is then $A_1\oplus\hdots\oplus A_m$ a $p$-elliptic tensor?
}
\end{center}

First we prove that when $A$ is a direct sum of $m$ elliptic matrices, we may take {\it real} testing vectors $\omega$.

\begin{lemma}
\label{l: complex vs. real omega}
Let $M^1,\hdots,M^m\in\C^{n\times n}$ be elliptic matrices. Define the complex tensor $A=M^1\oplus\hdots\oplus M^m$, acting on $\C_n^m$. Then we may in \eqref{eq: lambda_p} replace $\omega\in\Sigma^m$ by real unit vectors $\omega\in S^{m-1}$ with nonnegative components. That is, for $p>1$ we have
\begin{equation*}
\lambda_p(A)=
\min_{\xi\in\Sigma^m_{n}}
\min_{\omega\in S_+^{m-1}}
\Re\sk{A\xi}{\xi+(p-2)\Phi_\omega\xi}_{\C^m_n}.
\end{equation*}
\end{lemma}

\begin{proof}
We start by our reasoning, based on \eqref{eq: lambda_p}, \eqref{eq: Phi up/down} and \eqref{eq: quasi inner product - up summation}. Take $\xi\in\C_n^m$. Owing to the tensor structure of $A$, we have
\begin{equation*}
\aligned
 & \sk{A\xi}{\xi+(p-2)\Phi_\omega\xi} \\
& = 
\Sk{	(M^1\xi^1, \hdots, M^m\xi^m) }
{
(
\xi^1+(p-2)\omega^1\Re[\xi,\omega],	 \hdots,	 \xi^m+(p-2)\omega^m\Re[\xi,\omega]
)
}\\
&=\sum_{\alpha=1}^m\Sk{M^\alpha\xi^\alpha}{\xi^\alpha+(p-2)\omega^\alpha\Re[\xi,\omega]}\\
&=\sum_{\alpha=1}^m
\left(
\Big\langle{M^\alpha\xi^\alpha},{\xi^\alpha}\Big\rangle
+(p-2)
\Big\langle{M^\alpha(\overline \omega^\alpha\xi^\alpha)},{\Re\sum_{\beta=1}^m\overline{\omega}^\beta\xi^\beta}\Big\rangle
\right).
\endaligned
\end{equation*}
By taking $\sigma^\alpha=\xi^\alpha$ if $\omega^\alpha=0$ and otherwise
\begin{equation}
\label{eq: sigma}
\sigma^\alpha
=\frac{\overline\omega^\alpha}{|\omega^\alpha|}\xi^\alpha,
\end{equation}
so that $\overline \omega^\alpha\xi^\alpha=|\omega^\alpha|\sigma^\alpha$,
we get
\begin{equation}
\label{eq: change of var}
\aligned
& \sk{A\xi}{\xi+(p-2)\Phi_\omega\xi}
=\sum_{\alpha=1}^m
\Big\langle M^\alpha\sigma^\alpha, 
\sigma^\alpha
+(p-2)
|\omega^\alpha|
\Re\sum_{\beta=1}^m
|\omega^\beta|
\sigma^\beta
\Big\rangle.
\endaligned
\end{equation}
We proved that
\begin{equation}
\label{eq: quad form sigma}
\Re \sk{A\xi}{\cJ_{p,\omega}\xi}
=
\Re \sk{A\sigma}{\cJ_{p,\wt\omega}\sigma},
\end{equation}
where $\wt\omega=\left(|\omega^1|,\hdots,|\omega^m|\right)\in S_+^{m-1}$ and $\sigma=(\sigma^1,\hdots,\sigma^m)$.

Since $|\sigma^\alpha|=|\xi^\alpha|$, this finishes the proof.
\end{proof}

Note that we proved \eqref{eq: quad form sigma} even without $\Re$.

\begin{notation}
The $p$-ellipticity of the matrix $M\in\C^{n\times n}$ means that, uniformly in $\zeta\in\C_n$, we have
$$
\Re\sk{M\zeta}{\zeta+(p-2)\Re\zeta}\geqsim|\zeta|^2.
$$
With this in mind, introduce the {\it ad hoc} notation 
$$
\cK_r(\zeta)=\cK_{M,r}(\zeta):=\Re\sk{M\zeta}{\cJ_r\zeta}=\Re\sk{M\zeta}{\zeta+(r-2)\Re\zeta}.
$$
Thus $\cK_r$ is the (real) quadratic form whose positivity determines $r$-ellipticity.
\end{notation}

\begin{theorem}
\label{t: dir sum J_p formula}
Let $A^1,\hdots, A^m\in\cA(\Omega)$ and $p>1$. 
Assuming the notation \eqref{eq: sigma}, for every $\xi\in\C_n^m$ we have
\begin{equation*}
\aligned
\Re & \sk{(A^1\oplus\hdots\oplus A^m)\xi}{\cJ_{p,\omega}\xi} \\
&
=\Re\Sk{\sum_{j=1}^mA^j|\omega_{j}|\sigma_{j}}{\cJ_p\sum_{k=1}^m|\omega_{k}|\sigma_{k}}
+
\sum_{j< k}\Re\sk{A^j|\omega_k|\sigma_{j}-A^k|\omega_j|\sigma_k}{|\omega_k|\sigma_{j}-|\omega_j|\sigma_k}.
\endaligned
\end{equation*}
In particular, for $A^1=\hdots=A^m=:A$ this gives
\begin{equation*}
\label{Semjonov - Beloshitski Ispanskaya syuita 5. del}
\Re \sk{A^{\oplus m}\xi}{\cJ_{p,\omega}\xi} \\
=
\cK_{A,p}\bigg(\sum_{j=1}^m|\omega_j|\sigma_j\bigg)
+
\sum_{j<k}\cK_{A,2}\left(|\omega_j|\sigma_k-|\omega_k|\sigma_j\right). 
\end{equation*}
\end{theorem}

\begin{proof}
Let us first illuminate the idea behind devising the formul\ae\ in this theorem. Its starting point was the study of the special case 
\begin{center}
\framebox{$m=2$ and $A^1=A^2=:A$}. 
\end{center}

\smallskip
Take $\xi=(\xi_1,\xi_2)\in\C^n\times\C^n$ and $\omega=(\omega^1,\omega^2)\in\Sigma^2$. 
Write, for some 
$\f\in[0,\pi/2]$ and $\phi_1,\phi_2\in[0,2\pi)$, 
$$
\aligned
\omega^1&=e^{i\phi_1}\cos\f\\
\omega^2&=e^{i\phi_2}\sin\f
\endaligned
$$
and 
$$
\aligned
x=\sigma^1 & := e^{-i\phi_1}\xi_1\\
y=\sigma^2 & := e^{-i\phi_2}\xi_2.
\endaligned
$$
In particular,
$$
\aligned
|\omega^1|&=\cos\f\\
|\omega^2|&=\sin\f.
\endaligned
$$
From \eqref{eq: quad form sigma} we now obtain 
\begin{equation}
\label{eq: dir sum m=2}
\aligned
\Re & \sk{A^{\oplus 2}\xi}{\xi+(p-2)\Phi_\omega\xi} \\
&
\hskip 30 pt 
=\Re\sk{A x}{x}\sin^2\!\f +\Re\sk{Ay}{y}\cos^2\!\f \, \\
&\hskip 40 pt +\Re\sk{A x}{x+(p-2)\Re x}\cos^2\!\f+\Re\sk{Ay}{y+(p-2)\Re y}\sin^2\!\f \, \\
&\hskip 50pt +(p-2)\left(\sk{\Re(A x)}{\Re y}+\sk{\Re(Ay)}{\Re x}\right)\sin\f\cos\f.
\endaligned
\end{equation}

We feel that the expression in \eqref{eq: dir sum m=2} should relate to $\cK_{A,p}$, applied to some linear combination of $x,y$. So as the first attempt we set to calculate, for any pair of $a,b\in\R$ and $x,y\in\C^n$,
\begin{equation}
\label{eq: K formula m=2}
\aligned
\cK_{A,p}(ax+by) 
& = \hskip 5pt 
a^2  \Re\sk{Ax}{x+(p-2)\Re x}
\\ &  \hskip 8pt 
+ b^2  \Re\sk{Ay}{y+(p-2)\Re y}
\\
 & \hskip 32pt +(p-2)ab\left(\Re\sk{Ax}{\Re y}+\Re\sk{Ay}{\Re x}\right)\\
  & \hskip 67.6pt +ab\left(\Re\sk{Ax}{y}+\Re\sk{Ay}{x}\right).
\endaligned
\end{equation}

By comparing this again to \eqref{eq: dir sum m=2}, we observe that the bottom row in \eqref{eq: K formula m=2} is the most ``problematic'', as no such expression appears in 
 \eqref{eq: dir sum m=2}. 
We may try to neutralize it by adding a term of the form $\cK_{A,2}(\hdots)$ as a correcting factor, such that the coefficients in front of the ``problematic terms'' will cancel out.
Since with $p=2$ all the expressions involving $p-2$ will disappear, this will not affect the expression much; on the contrary, it may contribute to taking care of the first row.

That said, we try with 
$
\cK_{A,p}(ax+by)+\cK_{A,2}(cx+dy), 
$
where $a,b,c,d$ are yet to be determined. We get
$$
\aligned
\cK_{A,p}(ax+by)& +\cK_{A,2}(cx+dy)\\
& = \hskip 5pt c^2  \Re\sk{Ax}{x}+ d^2  \Re\sk{Ay}{y}\\
& \hskip 10pt + a^2  \Re\sk{Ax}{x+(p-2)\Re x}\\
&  \hskip 10pt + b^2  \Re\sk{Ay}{y+(p-2)\Re y}\\
 & \hskip 15pt +(p-2)ab\left(\Re\sk{Ax}{\Re y}+\Re\sk{Ay}{\Re x}\right)\\
  & \hskip 20pt +(ab+cd)\left(\Re\sk{Ax}{y}+\Re\sk{Ay}{x}\right).
\endaligned
$$
This expression will be exactly equal to \eqref{eq: dir sum m=2} if we can find $a,b,c,d$ such that
$$
\aligned
a^2&=\cos^2\f\\
b^2&=\sin^2\f\\
c^2&=\sin^2\f\\
d^2&=\cos^2\f\\
ab&=\sin\f\cos\f\\
ab+cd & =0.
\endaligned
$$
Clearly, this can be achieved 
by choosing 
$$
\aligned
a & = \cos\f\\
b & = \sin\f\\
c & = \pm\sin\f\\
d & = \mp\cos\f.
\endaligned
$$

Thus, with the notation above,
\begin{equation}
\label{eq: p-ell quad form via K m=2}
\Re  \sk{A^{\oplus 2}\xi}{\xi+(p-2)\Phi_\omega\xi}
=
\cK_{A,p}(x\cos\f+y\sin\f) +\cK_{A,2}(x\sin\f-y\cos\f).
\end{equation}
In other words,
\begin{equation*}
\Re  \sk{A^{\oplus 2}\xi}{\cJ_{p,\omega}\xi}
=
\cK_{A,p}(|\omega_1|\sigma_1+|\omega_2|\sigma_2) +\cK_{A,2}(|\omega_2|\sigma_1-|\omega_1|\sigma_2).
\end{equation*}
This settles the special case specified at the beginning of the proof and represents a starting point for the general case that we consider next.

\medskip
Now we have to address the general case, that is, when 

\begin{center}
\framebox{$m\in\N$ and $A^1,\hdots,A^m\in\cA(\Omega)$ are arbitrary}. 
\end{center}

\medskip
\noindent
We recall the notation \eqref{eq: sigma}. That is, we define $\sigma^j$ by $\overline\omega^j\xi^j=|\omega^j|\sigma^j$.
Consequently,
$
[\xi,\omega]=\sum_\alpha|\omega^\alpha|\sigma^\alpha.
$
Recall the definition of the operator $\cJ_p$ from Section \ref{s: aux ops} which gives, for $\sigma\in\C_n$, $(p-2)\Re\sigma=\cJ_p\sigma-\sigma$.
Now continue the calculation initiated in \eqref{eq: change of var}, as follows:
\begin{equation*}
\aligned
\Re & \sk{(A^1\oplus\hdots\oplus A^m)\xi}{\xi+(p-2)\Phi_\omega\xi} \\
&=\sum_{j=1}^m\Re\sk{A^j\sigma_j}{\sigma_j}\sum_k|\omega_k|^2
-\sum_{j,k=1}^m |\omega_j||\omega_k|\Re\sk{A^j\sigma_j}{\sigma_k}\\
&\hskip 40pt
+\sum_{j,k=1}^m |\omega_j||\omega_k|\Re\sk{A^j\sigma_j}{\cJ_p\sigma_k}\\
&=\sum_{j\ne k}\left(\Re\sk{A^j\sigma_{j}}{\sigma_{j}}|\omega_k|^2
-|\omega_{j}||\omega_k|\Re\sk{A^j\sigma_{j}}{\sigma_{k}}\right)\\
&\hskip 40pt
+\sum_{j,k=1}^{m}|\omega_{j}||\omega_k|\Re\sk{A^j\sigma_{j}}{\cJ_p\sigma_{k}}\\
&=\sum_{j< k}\left(\Re\sk{A^j\sigma_{j}}{\sigma_{j}}|\omega_k|^2+\Re\sk{A^k\sigma_{k}}{\sigma_{k}}|\omega_j|^2\right.\\
&\hskip 60pt \left.-|\omega_{j}||\omega_k|\Re\sk{A^j\sigma_{j}}{\sigma_{k}}-|\omega_{j}||\omega_k|\Re\sk{A^k\sigma_{k}}{\sigma_{j}}\right)\\
&\hskip 40pt
+\sum_{j,k=1}^m\Re\sk{A^j|\omega_{j}|\sigma_{j}}{\cJ_p|\omega_{k}|\sigma_{k}}\\
&
=\Re\Sk{\sum_{j=1}^mA^j|\omega_{j}|\sigma_{j}}{\cJ_p\sum_{k=1}^m|\omega_{k}|\sigma_{k}}
+
\sum_{j< k}\Re\sk{A^j|\omega_k|\sigma_{j}-A^k|\omega_j|\sigma_k}{|\omega_k|\sigma_{j}-|\omega_j|\sigma_k}.
\endaligned
\end{equation*}
This finishes the proof.
\end{proof}

\begin{corollary}
\label{c: A^m vs. A p-ell}
$A^{\oplus m}$ is a $p$-elliptic tensor if and only if $A$ is a $p$-elliptic matrix.
\end{corollary}

\begin{proof}

($\Rightarrow$)
This is Proposition \ref{p: lambda_p on diagonal}.

($\Leftarrow$)
By $p$-ellipticity of $A$, we get $\cK_{A,p}(\zeta)\geq\lambda_p(A)|\zeta|^2$. 
We combine this with Theorem \ref{t: dir sum J_p formula}.

When $m=2$, Theorem \ref{t: dir sum J_p formula} expressed as in \eqref{eq: p-ell quad form via K m=2} gives
$$
\aligned
\cK_{A,p}(x & \cos\f+y\sin\f) +\cK_{A,2}(x\sin\f-y\cos\f)\\
& \geqsim
|x\cos\f+y\sin\f|^2+|x\sin\f-y\cos\f|^2\\
& = |x|^2+|y|^2\\
& =|\xi|^2,
\endaligned
$$
which proves 
$
\Re  \sk{A^{\oplus 2}\xi}{\xi+(p-2)\Phi_\omega\xi}
\geqsim 
|\xi|^2,
$
as desired.

For general $m$, similarly,
\begin{alignat*}{2}
\cK_p\bigg(\sum_{j=1}^m|\omega_j|\sigma_j\bigg)
+
\sum_{j<k}& \cK_2\left(|\omega_j|\sigma_k-|\omega_k|\sigma_j\right)\\
& \geqsim
\bigg|\sum_{j=1}^m|\omega_j|\sigma_j\bigg|^2
+
\sum_{j<k}\big||\omega_j|\sigma_k-|\omega_k|\sigma_j\big|^2\\
&=|\xi|^2.
\end{alignat*}

Let us give an alternative proof of the  implication ($\Leftarrow$). 
It is not restrictive to assume that $C(x)=C$ is constant. As noted in Remark~\ref{r: commute}, $\mathrm{Id}_{\C^m}\otimes C$ commutes with every $\cU_\omega$, $\omega\in\Sigma^m$. Therefore, by \eqref{eq: Deltap bis}, we have 
\begin{equation*}
\Delta_p(\mathrm{Id}_{\C^m}\otimes C)=
\min_{\xi\in\Sigma_n^m 
}
\Re\sk{(\mathrm{Id}_{\C^m}\otimes C)\xi}{(\mathrm{Id}+(1-2/p)\cB)\xi}_{\C^m_n}.
\end{equation*}
Let $\xi\in\C^m_n$. By \eqref{eq: K in up}, we have
\begin{alignat*}{2}
\Re&\sk{(\mathrm{Id}_{\C^m}\otimes C)\xi}{(\mathrm{Id}+(1-2/p)\cB)\xi}_{\C^m_n}=\sum^m_{\alpha=1}\Re \sk{C\xi^\alpha}{\xi^\alpha+(1-2/p)(\cB\xi)^\alpha}_{\C_n}\\
&=\Re \sk{C\xi^1}{\xi^1+(1-2/p)\overline{\xi^1}}_{\C_n}+\sum^m_{\alpha=2}\Re \sk{C\xi^\alpha}{\xi^\alpha-(1-2/p)\xi^\alpha}_{\C_n}\\
&=\Re \sk{C\xi^1}{\xi^1+(1-2/p)\overline{\xi^1}}_{\C_n}+\frac2p\sum^m_{\alpha=2}\Re \sk{C\xi^\alpha}{\xi^\alpha}_{\C_n}\\
&\geq \Delta_p(C)|\xi^1|^2+\frac2p\lambda_2(C)\sum^m_{\alpha=2}|\xi^\alpha|^2\\
&\geq \frac2p\min\{\lambda_p(C),\lambda_2(C)\}|\xi|^2.
\tag*{\qedhere}
\end{alignat*}
\end{proof}

\subsection{The $p$-ellipticity under inversion}
\label{s: inversion}

Let $A$ be an elliptic tensor on $\Omega$. Fix $x\in\Omega$ and write $A=A(x)$. This can be perceived as a linear map $A:\C_n^m\rightarrow\C_n^m$. From the ellipticity condition we get for almost every $x\in\Omega$, with 
$\lambda=\lambda(A)$ and $\Lambda=\Lambda(A)$,
$\lambda|\xi|^2\leq\Re\sk{A\xi}{\xi}\leq|A\xi||\xi|$, 
thus
$|A\xi|\geq\lambda|\xi|$.
Consequently, $A=A(x)$ is invertible for a.e. $x\in\Omega$ and 
$$
\nor{A^{-1}}\leq\frac{1}{\lambda}\,.
$$

\begin{theorem}
\label{t: lambda_p of 1/A}
Let $A\in\cA^m(\Omega)$ and $p>1$. 
If $A$ is $p$-elliptic, then so is $A^{-1}$. In this case,
\begin{equation*}
\lambda_p(A^{-1})\geq\frac{\lambda_p(A)}{\Lambda(A)^2}.
\end{equation*}
\end{theorem}

\begin{proof}
We saw before that, for a.e. $x\in\Omega$, the map $A(x):\C_n^m\rightarrow\C_n^m$ is invertible. Fix such $x$ and denote $A^{-1}=A(x)^{-1}$. 
Since $A$ is $p$-elliptic, we have, by \eqref{eq: p-ell},
\begin{equation*}
\Re\sk{A\xi}{\xi}+(p-2)\Re\sk{A\xi}{\Phi_\omega\xi}\geq\lambda_p(A)|\xi|^2
\end{equation*}
for all $\xi\in\C_n^m$ and $\omega\in\Sigma^m$. Write $\eta=A\xi$. Then the above inequality reads
\begin{equation*}
\Re\sk{\eta}{A^{-1}\eta}+(p-2)\Re\sk{\eta}{\Phi_\omega(A^{-1}\eta)}\geq\lambda_p(A)|A^{-1}\eta|^2.
\end{equation*}
By Proposition \ref{p: Phi bilinear form} and the simple fact that $\Re\sk{\mathbf{u}}{\mathbf{v}}=\Re\sk{\mathbf{v}}{\mathbf{u}}$, this becomes
\begin{equation*}
\Re\sk{A^{-1}\eta}{\eta}+(p-2)\Re\sk{A^{-1}\eta}{\Phi_\omega\eta}\geq\lambda_p(A)|A^{-1}\eta|^2.
\end{equation*}
The definition of $\Lambda$ implies $|A\xi|\leq\Lambda(A)|\xi|$, thus $|A^{-1}\eta|\geq\Lambda(A)^{-1}|\eta|$.
Since $A$ is $p$-elliptic, meaning that $\lambda_p(A)>0$,
the above estimates combine into
$$
\Re\sk{A^{-1}(x)\eta}{\cJ_{p,\omega}\eta}\geq\frac{\lambda_p(A)}{\Lambda(A)^2}|\eta|^2
$$
for every $\eta\in\C_n^m$ and every $\omega\in\Sigma^m$. 
This holds for almost every $x\in\Omega$, 
therefore the theorem is proved.
\end{proof}

\subsection{The $p$-ellipticity with respect to the conjugate exponent}

Let $A$ be a tensor on $\Omega$. Define the tensor $A^*$ as the adjoint with respect to the pairing in $\C_n^m$. That is, we request, for any $\xi,\eta\in\C_n^m$,
$$
\sk{A\xi}{\eta}_{\C_n^m}=\sk{\xi}{A^*\eta}_{\C_n^m}.
$$
A quick calculation gives that this is satisfied if
$$
(A^*)_{jk}^{\alpha\beta}:=\overline{a_{kj}^{\beta\alpha}}.
$$
Equivalently, $A^*_{j,k}:=(A^*)_{j,k}=(A_{k,j})^*$ or $(A^*)^{\alpha,\beta}=(A^{\beta,\alpha})^*$.

\begin{proposition}
\label{p: p-ell of A vs. A*}
Let $1/p+1/q=1$. Then $A$ is $p$-elliptic if and only if $A^*$ is $q$-elliptic. 
That is, $\cA_p^m(\Omega)^*=\cA_q^m(\Omega)$.
\end{proposition}

By Corollary \ref{c: DiLiPi} we have that this was already proven in \cite[Theorem 1.1 {\it (iv)}]{DiLiPi20}. That statement, in turn, followed an earlier statement in the same vein proven in \cite[Corollary 5.17]{CD-DivForm}.

\begin{proof}
This is a consequence of Corollaries \ref{c: JoJ} and \ref{c: norm of J}. 

Recall the notation from Section \ref{s: auxops}.
Then
$$
\aligned
\Re\sk{A\xi}{\cJ_{p,\omega}\xi}
& =
\sum_{j,k=1}^n
\Re\Sk{A_{j,k}\xi_k}{\iota_{p,\omega}\xi_j}_{\C^m}\\
&=
\sum_{j,k=1}^n
\Re\Sk{\xi_k}{(A_{j,k})^*(\iota_{p,\omega}\xi_j)}_{\C^m}
.
\endaligned
$$
Write $\eta_j:=\iota_{p,\omega}\xi_j$. Then, by Corollary \ref{c: JoJ}, $\xi_j=\iota_{q,\omega}\eta_j$. 
Thus we obtained
$$
\Re\sk{A\xi}{\cJ_{p,\omega}\xi}
=
\sum_{j,k=1}^n
\Re\Sk{\iota_{q,\omega}\eta_k}{A_{k,j}^*\eta_j}_{\C^m}
=
\Re\sk{A^*\eta}{\cJ_{q,\omega}\eta}.
$$
Since $A$ is $p$-elliptic and, by Corollary \ref{c: norm of J}, $|\eta|=|\cJ_{p,\omega}\xi|\sim|\xi|$, we conclude that $A^*$ is $q$-elliptic.
\end{proof}

\section{Invariance under conjugation of $p$}

Recall 
\eqref{eq: Delta vs. lambda}
and the notation from Section \ref{s: aux ops}.
Let $q$ be the conjugate exponent of $p$, in other words, a number given by $1/p+1/q=1$.
When \framebox{$m=1$}, we have 
\begin{equation}
\label{eq: invariance via Delta}
\Delta_p(A)=\Delta_q(A). 
\end{equation}

Consequently, still for $m=1$,
$p$-ellipticity and $q$-ellipticity are equivalent, that is, 
\begin{equation}
\label{eq: invariance via cA}
\cA_q^m(\Omega)=\cA_p^m(\Omega). 
\end{equation}
The proof of \eqref{eq: invariance via Delta}
(in the scalar case $m=1$) 
goes through observing that
\begin{equation}
\label{eq: Ip Iq scalar}
\cI_p(i\xi)=i\cI_q\xi. 
\end{equation}
(Equivalently, $\cI_q=\cM_i\cI_p \cM_i^{-1}$, with $\cM_i$ as in Remark \ref{r: I,J at m=1}.) Indeed, with $\eta=i\xi$ we then have, by \eqref{eq: Ip Iq scalar},
$$
\Re\sk{A\eta}{\cI_p\eta}_{\C^{n}}
=
\Re\sk{iA\xi}{i\cI_q\xi}_{\C^{n}}
=
\Re\sk{A\xi}{\cI_q\xi}_{\C^{n}}.
$$
Since the unit sphere in $\C^n$ is invariant under multiplication by $i$, by recalling \eqref{eq: Delta_p} this gives \eqref{eq: invariance via Delta}.

\medskip
In this section we want to investigate what happens for general \framebox{$m\in\N$}. That is, we want to understand whether
$p$-ellipticity and $q$-ellipticity are still equivalent.
We start with a vector-case variant of \eqref{eq: Ip Iq scalar}.

\begin{proposition}
Assume that $p>1$, $1/p+1/q=1$, $X\in\C_n^m$, $\omega\in\Sigma^m$. Then 
\begin{eqnarray}
& \cI_{p,\omega} (iX) \hskip -7pt & =i\cI_{p,i\omega} X\nonumber\\
\label{eq: Ip Iq vector}
&\cI_{p,i\omega} X  \hskip -16pt  & =\cI_{q,\omega} X
+
2(1-2/q)(X-[X,\omega]\omega).
\end{eqnarray}
\end{proposition}

\begin{proof}
The first identity is a special case of Proposition \ref{p: pi - unitary similarity} for $R=\cM_i$, cf. Remark  \ref{r: I,J at m=1}.
So let us turn to the second one. 
We will repeatedly use \eqref{eq: Ip Iq scalar}. By Lemma \ref{l: I,J vect vs. I,J scal}
we have
\begin{alignat*}{2}
\cI_{p,i\omega} X
&=\cI_{p}([X,i\omega])i\omega+\frac2p(X-[X,i\omega]i\omega)\\
&=\cI_{p}(-i[X,\omega])i\omega+\frac2p(X-[X,\omega]\omega)\\
&=\cI_{q}([X,\omega])\omega+\frac2p(X-[X,\omega]\omega)\\
&=\cI_{q,\omega} X
+\left(\frac2p-\frac2q\right)
(X-[X,\omega]\omega).
\tag*{\qedhere}
\end{alignat*}
\end{proof}

As a curiosity we note that \eqref{eq: Ip Iq vector} reads in a more symmetric way as
$$
\cI_{p,i\omega} X
+(1-2/p)(X-[X,i\omega]i\omega)
=
\cI_{q,\omega} X
+(1-2/q)(X-[X,\omega]\omega).
$$
This is the vector analogue of \eqref{eq: Ip Iq scalar}. 

\medskip
If $m=1$ then $X-[X,\omega]\omega=0$. 
Consequently,
\eqref{eq: Ip Iq vector} then reads 
$\cI_{p,i\omega} X =\cI_{q,\omega} X$,
which reduces to just \eqref{eq: Ip Iq scalar}, as expected.

If $m\ne1$, however, the term $X-[X,\omega]\omega$ is nonzero, moreover, replacing $\omega$ by $i\omega$ there does not affect it.

\subsection{Conjugation invariance for $m=1$ or $n=1$}
Suppose that $1/p+1/q=1$ and $1<q<2<p<\infty$. 
Let $A$ be an elliptic tensor. For the purpose of this section we may assume that it is constant. 
We want to relate $\Delta_p(A)$ to $\Delta_q(A)$. 
By \eqref{eq: Delta via Xi} we get
\begin{equation}
\label{eq: DpXpm}
\aligned
\Delta_p(A)
&=\min_{|\xi|=1}(\Re\sk{A\xi}{\xi}-(1-2/p)\Xi^+_A(\xi))\\
\Delta_q(A)
&=\min_{|\xi|=1}(\Re\sk{A\xi}{\xi}-(1-2/q)\Xi^-_A(\xi))\\
&=\min_{|\xi|=1}\big(\Re\sk{A\xi}{\xi}-(1-2/p)(-\Xi^-_A(\xi))\big).
\endaligned
\end{equation}
So we want to understand whether there is any relation between $\Xi^+_A(\xi)$ and $\Xi^-_A(\xi)$.
To be more precise, in fact we want to know the {\it sign} of $\Xi^+_A(\xi)+\Xi^-_A(\xi)$.
Even more concretely, if we could prove that
\begin{equation}
\label{eq: X++X-}
\Xi^+_A(\xi)+\Xi^-_A(\xi)
=0,
\end{equation}
this would obviously imply \eqref{eq: invariance via Delta}.

\medskip
We emphasize, however, that \eqref{eq: X++X-} looks like a ``pointwise'' variant of 
\eqref{eq: invariance via Delta}, that is, an identity before even taking the minima in \eqref{eq: DpXpm}.
So it seems like a considerably stronger statement than \eqref{eq: invariance via Delta}.

\begin{remark}
From \eqref{eq: Xi via Theta} we get that \eqref{eq: X++X-} is equivalent to
\begin{equation}
\label{eq: T++T-}
\Theta_A^+(\xi)+\Theta_A^-(\xi)=\Re\sk{A\xi}{\xi}.
\end{equation}
\end{remark}

By \eqref{eq: Xi^pm}, the identities \eqref{eq: X++X-} and \eqref{eq: T++T-} would follow if we could prove that 

\begin{equation}
\label{eq: J,-J}
\aligned
&for\ every\ \xi\in\C_n^m\ and\ \omega\in\Sigma^m\\ 
&there\ exists\ \omega'\in\Sigma^m\ such\ that\ 
\cJ_{0,\omega'}\xi=-\cJ_{0,\omega}\xi. 
\endaligned
\end{equation}

In order to study \eqref{eq: J,-J}, we need the following result.

\begin{proposition}
\label{p: surjective}
Take $y\in\Sigma^m=\mn{\xi\in\C^m}{|\xi|=1}$. Then the map
$$
\begin{array}{rccl}
f=f_y: &  \Sigma^m & \longrightarrow & \Sigma^m\\
          &  \zeta    &   \longmapsto    & \cJ_{0,\zeta}y=y-2\Re\sk{y}{\zeta}\zeta
\end{array}
$$
is surjective. 
\end{proposition}

\begin{proof}
If $m=1$ then $f_y(\zeta)=-\bar y\zeta^2$, thus $f_y$ is clearly surjective. 

Now let $m>1$.
By \eqref{eq: J_0 isometry}, applied with $n=1$, the map indeed maps into $\Sigma^m$.
Take $\varpi\in\Sigma^m$. We want to find $\zeta\in\Sigma^m$ such that $\varpi=f(\zeta)=\cJ_{0,\zeta}y$, that is,
\begin{equation}
\label{eq: varpi id}
2\Re\sk{y}{\zeta}\zeta
=
y-\varpi.
\end{equation}
If $y=\varpi$, we may take for $\zeta$ anything from $y^\perp\cap\Sigma^m$ (this is possible, as $m>1$). 

Now assume that $y\ne\varpi$. On the left-hand side of \eqref{eq: varpi id} we have a real multiple of $\zeta$, while on the right-hand side we have (a real multiple of) $y-\varpi$. 
Moreover, we are supposed to have $\zeta\in\Sigma^m$. Thus we must try with
\begin{equation}
\label{eq: zeta def}
\zeta=\pm\frac{y-\varpi}{|y-\varpi|}.
\end{equation}
Since $2(1-\Re\sk{y}{\varpi})=|y-\varpi|^2$,
this indeed solves \eqref{eq: varpi id}.
\end{proof}

\begin{theorem}
\label{t: invar m=n=1}
When $m=1$ or $n=1$, the statement \eqref{eq: J,-J} holds. 
Consequently, for any elliptic tensor $A$, the identities \eqref{eq: X++X-} and \eqref{eq: invariance via Delta} hold, as well. 
\end{theorem}

\begin{proof}
Take $\xi,\omega$ as in the assumption of \eqref{eq: J,-J}. We need to find $\omega'\in\Sigma^m$ such that for all $h=1,\hdots,n$ we have
\begin{equation*}
(\cJ_{0,\omega'}\xi)_h=(-\cJ_{0,\omega}\xi)_h,
\end{equation*}
that is, in the notation of Section \ref{s: J_{p,omega}},
\begin{equation}
\label{eq: iota, -iota}
\iota_{0,\omega'}(\xi_h)=-\iota_{0,\omega}(\xi_h).
\end{equation}
If $\xi_h=0$, this is trivially satisfied for any $\omega,\omega'$. So we may assume that $\xi_h\ne0$.

In terms of Proposition \ref{p: surjective}, the identity \eqref{eq: iota, -iota} reads 
$f_{\xi_h/|\xi_h|}(\omega')=-f_{\xi_h/|\xi_h|}(\omega)$.
For a {\it fixed index} $h$, we can solve this problem by using Proposition \ref{p: surjective}. Moreover,
\eqref{eq: zeta def} tells us that, unless $\xi_h$ is a real multiple of $\omega$, the solution is given by 
$$
\omega'=\pm\frac{\xi_h+\iota_{0,\omega}\xi_h}{|\hdots|}\,.
$$
For a given $z\in\C^m$ we calculate 
$
z+\iota_{0,\omega}z 
=
2(z-\Re\sk{z}{\omega}\omega)
=
2\iota_{1,\omega}(z)\,.
$
Thus 
\begin{equation}
\label{eq: omega' def}
\omega'
=\pm\frac{\xi_h-\Re\sk{\xi_h}{\omega}\omega}{|\hdots|}
=\pm\frac{\iota_{1,\omega}(\xi_h)}{|\iota_{1,\omega}(\xi_h)|}
\,.
\end{equation}
The right-hand side above in principle changes with $h$, while the left-hand side should be independent of $h$. 
Hence \eqref{eq: omega' def} surely defines a solution only if $n=1$.
Given that the solution \eqref{eq: omega' def} lies in the (real) linear span of $\xi_h$ and $\omega$, when $n>1$ we cannot expect it to always hold simultaneously for all $h$. 

Apart from the case $n=1$, the formula \eqref{eq: omega' def} also defines a solution when $m=1$. Indeed, then for $z\in\C_n$ and $\omega\in\Sigma^1$ we obtain
$$
\aligned
z-\Re\sk{z}{\omega}\omega
& = z-\frac{(\overline\omega z+\omega\overline z)\omega}{2}\\
& = \frac{z-\omega^2\overline z}{2}\\
& = \frac{(\overline\omega z-\omega\overline z)\omega}{2}\\
& = \Im(\overline\omega z)i\omega.
\endaligned
$$
Note that $z=\xi_h$ is only contained in the {\it real} term $\Im(\overline\omega z)$, which disappears after normalizing, cf. \eqref{eq: omega' def}. Thus for $m=1$ we obtain $\omega'=\pm i\omega$, which is of course independent of $h$ and thus admissible as a solution of \eqref{eq: J,-J}.
This is consistent with the discussion at the beginning of this section, cf. \eqref{eq: Ip Iq scalar}.
\end{proof}

\begin{remark}
As mentioned before, the case $m=1$ was already established in \cite{CD-DivForm}; the real novelty here is the case $n=1$.
\end{remark}

\subsubsection{Summary of the known results on conjugation invariance of $\Delta_p$}

\begin{theorem}
Let $p>1$ and $A$ be a $p$-elliptic tensor. Suppose that either:
\begin{itemize}
\item
$m=1$
\item
\ $n=1$
\item
$A$ is unitary.
\end{itemize}
Then $A$ is $q$-elliptic, where $1/p+1/q=1$.
\end{theorem}

\begin{proof}
If $m=1$ or $n=1$, use Theorem \ref{t: invar m=n=1}.

Now assume that $A$ is unitary, meaning that $A(x)$ is unitary for almost every $x\in\Omega$.
Since $A$ is $p$-elliptic, it is elliptic (Corollary \ref{t: A_p inclusions}), hence it is also invertible (Section \ref{s: inversion}). Denote $B=A^{-1}$. Since $A$ is $p$-elliptic, so is $B$ (Theorem \ref{t: lambda_p of 1/A}). By Proposition \ref{p: p-ell of A vs. A*}, $B^*$ is $q$-elliptic. Since $A$ is unitary, $B^*=A$. Thus we proved that $A$ is $q$-elliptic.
\end{proof}

It turns out that, in general, the invariance \eqref{eq: invariance via cA} fails when $m>1$. We will give a proof of this fact in a forthcoming paper.

\section{Examples}

Here we present a few simple real elliptic tensors whose $p$-ellipticity coefficients $\Delta_p$ we are in a position to calculate explicitly. Doing so, we show that in the vector case ($m>1$), real elliptic tensors may not necessarily be $p$-elliptic for all $p>1$, unlike in the scalar case.

\subsection{The $p$-ellipticity for direct sums of positive multiples of identity}

Here we show that a direct sum of {\it different} $p$-elliptic matrices may not be $p$-elliptic. 
In fact, this already fails with the simplest possible example that we present next: real multiples of the identity. 
Another example we have are direct sums of rotation matrices; see Theorem \ref{t: Delta_p of U_theta,k}. 

So let $a,b>0$. Then $aI,bI$ are $p$-elliptic for any $p>1$, cf. \cite[p.3179]{CD-DivForm}. In contrast, we claim that $(aI)\oplus(bI)$ is not $p$-elliptic for all $p>1$, unless $a=b$. 

If $p>1$, denote 
\begin{equation*}
\label{eq: p*}
p^*=\max\{p,q\},
\hskip 30pt
\text{ where }\frac1p+\frac1q=1.
\end{equation*}

\begin{theorem}
\label{t: aI+bI}
Let $a,b>0$ and define the tensor $T=T_{a,b}=(aI)\oplus(bI)$ acting on $\C_n^2$. Then for every $p>1$ we have
$$
\Delta_p(T)=
\left\{
\renewcommand{\arraystretch}{2.4}
\begin{array}{ccl}
{\displaystyle \frac{2}{p^*}\,\min\{a,b\}} & ; & \text{if } {\displaystyle |1-2/p|\leq\frac{2\min\{a,b\}}{a+b}}\\
 {\displaystyle\frac{a+b}4\left(\frac{4ab}{(a+b)^2}-(1-2/p)^2\right)} 
 & ; &  \text{if } {\displaystyle |1-2/p|>\frac{2\min\{a,b\}}{a+b}}.
\end{array}
\renewcommand{\arraystretch}{1}
\right.
$$
Consequently, $(aI)\oplus(bI)$ is $p$-elliptic if and only if $|1-2/p|<2\sqrt{ab}/(a+b)$.
\end{theorem}

Note the contrast between the general case and the special case of $a=b$, when we have simply 
$$
\Delta_p(T_{a,a})=\frac{2}{p^*}\,a.
$$

\begin{proof}
We may assume that \framebox{$a<b$}.
From Lemma \ref{l: complex vs. real omega} we have 
$$
\lambda_p(T)=\min_{\xi\in\Sigma_n^2}\min_{\f\in[0,2\pi)}
\Re \sk{T\xi}{\xi+(p-2)\Phi_\omega\xi}_{\C_n^2}, 
$$
where $\omega=(\cos\f,\sin\f)$.

If $p=2$, the minimum of this expression is $\min\{a,b\}$ which we assume to be $a$.

In the continuation we work with $p>2$; for $1<p<2$ we repeat the proof in the $p>2$ case with the required changes; these are obvious, so we leave that case out.

Write $\xi=(x,y)$. 
Picking up from \eqref{eq: dir sum m=2}, we get, after an obvious adaptation to the case $A^1\ne A^2$, 
\begin{equation*}
\aligned
\Re & \sk{T\xi}{\xi+(p-2)\Phi_\omega\xi} \\
&= 
a|x|^2+b|y|^2\\
&\hskip 10 pt +
(p-2)
\left[
a|\Re x|^2\cos^2\f+b|\Re y|^2\sin^2\f+(a+b)\sk{\Re x}{\Re y}\sin\f\cos\f
\right]\\
&= 
a|x|^2+b|y|^2+
\frac{p-2}{2}
\left(
a|\Re x|^2+b|\Re y|^2
\right)\\
& 
\hskip 10pt
+
\frac{p-2}{2}
\left[
\left(a|\Re x|^2-b|\Re y|^2\right)\cos2\f+(a+b)\sk{\Re x}{\Re y}\sin2\f
\right].
\endaligned
\end{equation*}
Minimizing the above expression over $\f$ gives that $\lambda_p(T)$ is the minimum, taken over all $\xi=(x,y)\in\Sigma_n^2$, of
\begin{equation*}
\aligned
a|x|^2+b|y|^2
+\frac{p-2}{2}
\Big(
a|\Re x|^2+b|\Re y|^2
-\sqrt{\left(a|\Re x|^2-b|\Re y|^2\right)^2
+(a+b)^2\sk{\Re x}{\Re y}^2}
\Big). 
\endaligned
\end{equation*}
The above expression depends only on the {\it modulus} of $\Re x,\Re y$ -- except in the inner product. The latter, however, is the largest when the two vectors are aligned, and then it equals to the product of $|\Re x|,|\Re y|$. So we may write $|\Re x|=u$, $|\Re y|=v$ and $\lambda_p(T)$ is then the minimum over all $(x,y)\in\Sigma_n^2$ of
\begin{equation*}
\aligned
a|x|^2+ & b|y|^2
+\frac{p-2}{2}
\Big(
au^2+bv^2
-\sqrt{\left(au^2-bv^2\right)^2
+(a+b)^2u^2v^2}
\Big)\\
& =\
a|\Im x|^2+b|\Im y|^2
+\frac{p}{2}
\left(
au^2+bv^2
-(1-2/p)\sqrt{(a^2u^2+b^2v^2)(u^2+v^2)}
\right).
\endaligned
\end{equation*}
Since $|x|^2+|y|^2=1$, we may write
$$
\aligned
u=|\Re x|&=\cos\alpha\cos s\\
v=|\Re y|&=\cos\alpha\sin s\\
|\Im x|&=\sin\alpha\cos t\\
|\Im y|&=\sin\alpha\sin t.
\endaligned
$$
Then we get to minimize, over $\alpha,s,t$, the expression
$$
\aligned
& (a\cos^2 t+b\sin^2t)\sin^2\alpha\\
& \hskip 40pt
+\frac{p}{2}\,\cos^2\alpha
\left(
a\cos^2s+b\sin^2s
-(1-2/p)\sqrt{a^2\cos^2s+b^2\sin^2s}
\right).
\endaligned
$$
The dependence on $t$ only appears within the first set of parentheses, where we have convex combinations of $a,b$, thus
the minimum is $\min\{a,b\}=a$. Hence we now have to minimize
$$
a\sin^2\alpha
+\frac{p}{2}\cos^2\alpha
\left(
a\cos^2s+b\sin^2s
-(1-2/p)\sqrt{a^2\cos^2s+b^2\sin^2s}
\right).
$$
Again note that a linear combination of $\cos^2\alpha$ and $\sin^2\alpha$ is a convex combination of the coefficients, thus its minimum is the minimum of the coefficients. 
This means that $\lambda_p(T)$ is the minimum between $\min\{a,b\}=a$ and 
$$
\min_{s\in[0,2\pi)}
\frac{p}{2}
\left(
a\cos^2s+b\sin^2s
-(1-2/p)\sqrt{a^2\cos^2s+b^2\sin^2s}
\right).
$$
Write $\cos^2s=1-\tau$, $\sin^2s=\tau$ for $\tau\in[0,1]$. Then the above expression is 
$$
f(\tau)
:=
\frac{p}{2}
\left(
a+(b-a)\tau
-(1-2/p)\sqrt{a^2+(b^2-a^2)\tau}
\right).
$$
It is at this point elementary to verify that the function $f$, extended to all $\tau\in\R$ such that $a^2+(b^2-a^2)\tau>0$, attains its minimum (again, provided that $a<b$) at
$$
\tau=\tau_0
=\frac14 \cdot \frac{b+a}{b-a} 
	\left(
 			(1-2/p)^2-\Big(\frac{2a}{a+b}\Big)^2
	\right)
$$
and there it amounts to 
$$
f(\tau_0)
=\frac{p(a+b)}8\cdot
\left(
\frac{4ab}{(a+b)^2}-(1-2/p)^2
\right).
$$
We consider two cases:
\begin{enumerate}[1.)]
\item
First, \framebox{$\tau_0\in(0,1)$}. That is,  
$$
\frac{2a}{a+b}<1-2/p<\frac{2b}{a+b}.
$$
Since we always have $1-2/p<1<2b/(a+b)$, the second condition is empty. Thus we only require
$$
\frac{2a}{a+b}<1-2/p.
$$
In this case, $f(\tau_0)$ is the minimum of $f$ on $[0,1]$.

\item
Second, \framebox{$\tau_0\not\in(0,1)$}. Then $f$ is monotonous on $[0,1]$, therefore its minimum there equals the minimum of $f(0),f(1)$. Now notice that $f(0)=a$, $f(1)=b$.

\end{enumerate}
Finally, keep in mind that in order to obtain $\Delta_p(T)$, we need to multiply by $2/p$.
\end{proof}

\subsection{The $p$-ellipticity of rotation matrices and related tensors}

The second set of examples of ours is built around euclidean rotation matrices. We start by revisiting their relation to $p$-ellipticity.

\subsubsection{Scalar rotation matrices}

Take $\vt\in[-\pi,\pi)$ and let $R_\vt\in\R^{2\times2}$ be the standard  $2\times2$ rotation matrix by the angle $\vt$:
\begin{equation}
\label{eq: R_theta}
R_\vt=
\begin{bmatrix}[r]
\cos\vt & -\sin\vt\\
\sin\vt & \cos\vt
\end{bmatrix}.
\end{equation}
Clearly, $R_\vt^*=R_{-\vt}$. For $k\in\N$ let 
\begin{equation}
\label{eq: R_theta,k}
R_{\vt,k}=R_\vt\otimes I_k
=
\begin{bmatrix}[r]
\cos\vt & -\sin\vt\\
\sin\vt & \cos\vt
\end{bmatrix}
\otimes
I_k
=
\begin{bmatrix}[r]
(\cos\vt) I_k
& 
(-\sin\vt) I_k 
\\
(\sin\vt) I_k
& 
(\cos\vt) I_k
\end{bmatrix}
\in\R^{(2k)\times(2k)}
\end{equation}
be the corresponding matrix acting on $\C^{2k}$. 
Explicitly, for $(x_1,x_2)\in\C^k\times\C^k$ we have
$$
R_{\vt,k}(x_1,x_2)=(x_1\cos\vt-x_2\sin\vt,x_1\sin\vt+x_2\cos\vt)\in\C^k\times\C^k.
$$

\begin{proposition}
\label{p: Delta_p of R_theta,k}
Let $\vt\in[-\pi,\pi)$ and $k\in\N$. 
Then for every $p>1$, 
$$
\Delta_p(R_\vt\otimes I_k)=\min\left\{\frac2p\cos\vt,\frac2q\cos\vt\right\},
$$
where $1/p+1/q=1$.

In particular, for any $p>1$, $R_\vt\otimes I_k$ is $p$-elliptic if and only if $\vt\in(-\pi/2,\pi/2)$.
\end{proposition}

\begin{proof}
For any $\xi=(u,v)\in\C^k\times\C^k$,
$$
\aligned
& \hskip -40pt
\Re\sk{(R_\vt\otimes I_k)\xi}{\xi+(p-2)\Re\xi}_{\C^{2k}}\\
& =
\Re\Sk{\begin{bmatrix}[r]
\cos\vt & -\sin\vt\\
\sin\vt & \cos\vt
\end{bmatrix}
\begin{bmatrix}
u\\
v
\end{bmatrix}
}{
\begin{bmatrix}
u\\
v
\end{bmatrix}
+(p-2)
\begin{bmatrix}
\Re u\\
\Re v
\end{bmatrix}
}_{\C^{2k}}\\
& =
\left[
\left(
|u|^2+|v|^2
\right)
+(p-2)
\left(
|\Re u|^2+|\Re v|^2
\right)
\right]
\cos\vt.
\endaligned
$$
We have
$0\leq|\Re u|^2+|\Re v|^2\leq|\xi|^2$, 
with both bounds being sharp. 
Thus for $p\geq2$ and $\cos\vt\geq0$ we get 
$$
\lambda_p(R_\vt\otimes I_k)
=
\min_{\xi\in\Sigma^{2k}}
\Re\sk{(R_\vt\otimes I_k)\xi}{\xi+(p-2)\Re\xi}_{\C^{2k}}
=
\cos\vt.
$$
The other combinations of $p,\vt$ are virtually identical.
\end{proof}

\subsubsection{Direct sum of rotations in the upper representation from \eqref{eq: A up-down}}

We saw in Proposition \ref{p: Delta_p of R_theta,k} that
$
\Delta_p(R_{\vt,k})=\left(1-|1-2/p|\right)\cos\vt
$
if $\vt\in(-\pi/2,\pi/2)$. 
Thus, in particular, such $R_{\vt,k}$ are $p$-elliptic for all $p>1$ (and $k\in\N$).

Passing to $m>1$ creates a stark contrast with the situation at $m=1$, in the sense that a direct sum of two rotations is a real elliptic tensor on $\C_n^2$ that may not be $p$-elliptic for all $p>1$. 
The example we present is for $m=2$, but surely it can be extended to arbitrary $m>1$. 
Compare with Theorem \ref{t: aI+bI}, where direct sums of positive multiples of identity were studied. 

If $R\in\C^{2k}$ is a matrix, let $R\oplus R^*$ be the tensor on $\C^{2k}_2$, defined as in Section \ref{s: dir sums}. 
That is,  for $(\xi^1,\xi^2)\in\C^{2k}\times\C^{2k}$ we set
$$
(R\oplus R^*)(\xi^1,\xi^2)=(R\xi^1,R^*\xi^2).
$$
Yet in other words,
$$
(R\oplus R^*)^{\text{up}}
:=
\begin{bmatrix}
R & \\
& R^*
\end{bmatrix}.
$$

\begin{theorem}
\label{t: Delta_p of U_theta,k}
Let $\vt\in(-\pi/2,\pi/2)$.
Denote $U_{\vt,k}:=R_{\vt,k}\oplus R_{\vt,k}^*$, 
where  $R_{\vt,k}$ is as in \eqref{eq: R_theta,k}. 
Then for $p>1$ we have
$$
\Delta_p(U_{\vt,k})=\cos\vt-|1-2/p|.
$$
\end{theorem}

\begin{proof}
We present the proof for $p>2$; the remaining case is virtually identical. 
As usual, we will use \eqref{eq: lambda_p} and \eqref{eq: Delta vs. lambda}. 
Take a unit vector 
$\xi=(x,y)\in\C^{2k}\times\C^{2k}=\C_{2k}^2$ and $\omega\in\Sigma^2$. By Lemma \ref{l: complex vs. real omega}, we can restrict ourselves to {\it real} vectors $\omega$. Thus write $\omega=(\cos\f,\sin\f)$.
Then, in order to calculate $\lambda_p(R\oplus R^*)$, we have to minimize over all such $\xi,\omega$ the expression,  cf. \eqref{eq: dir sum m=2},
\begin{equation*}
\aligned
\Re & \sk{(R\oplus R^*)\xi}{\xi+(p-2)\Phi_\omega\xi} \\
&= 
\Re\sk{Rx}{x}+\Re\sk{R^*y}{y}
+
(p-2)
\left[
\Re\sk{Rx}{\Re x}\cos^2\f+\Re\sk{R^*y}{\Re y}\sin^2\f\right.\\
&\left. \hskip 180pt
+(\Re\sk{Rx}{\Re y}+\Re\sk{R^*y}{\Re x})\sin\f\cos\f
\right].
\endaligned
\end{equation*}
We have 
$$
\aligned
\Re\sk{Rx}{x}& =|x|^2\cos\vt\\
\Re(Rx)& =R(\Re x)
\endaligned
$$
and thus
$$
\aligned
\Re\sk{Rx}{\Re x}& =\sk{\Re(Rx)}{\Re x}=|\Re x|^2\cos\vt\\
\Re\sk{Rx}{\Re y}+\Re\sk{R^*y}{\Re x}&=2\sk{R(\Re x)}{\Re y}.
\endaligned
$$
On top of using the above identities, we introduce the matrix 
$$
B=B_{\vt,k}=\frac{R_{\vt,k}}{\cos\vt},
$$
pass to double angles and finally recall that $|x|^2+|y|^2=1$. In this vein we obtain
\begin{equation*}
\aligned
\Re & \sk{(B\oplus B^*)\xi}{\xi+(p-2)\Phi_\omega\xi} \\
&= 
1+
\frac{p-2}{2}
\left[
(|\Re x|^2+|\Re y|^2)+(|\Re x|^2-|\Re y|^2)\cos 2\f+2\sk{B(\Re x)}{\Re y}\sin2\f
\right].
\endaligned
\end{equation*}
Write $X=\Re x, Y=\Re y$. 
After minimizing in $\f$ and remembering that $p>2$ we get 
\begin{equation*}
\aligned
\min_{\f\in[0,2\pi)}\Re & \sk{(B\oplus B^*)\xi}{\xi+(p-2)\Phi_{(\cos\f,\sin\f)}\xi} \\
&=
1+
\frac{p-2}{2}
\left[
(|X|^2+|Y|^2)-\sqrt{(|X|^2-|Y|^2)^2+4\sk{BX}{Y}^2}
\right].
\endaligned
\end{equation*}
The inner product can attain its largest value (hence, owing to the minus sign in front of the square root, the entire expression its minimum) when $Y=\lambda BX$ for some $\lambda>0$, when it equals $|BX||Y|=|X||Y|/\cos\theta$.
Hence 
$
\lambda_p(B\oplus B^*)
$
is the minimum of 
$$
1+
\frac{p-2}{2}
\left[
(|X|^2+|Y|^2)-\sqrt{(|X|^2-|Y|^2)^2+4|X|^2|Y|^2/\cos^2\vt}
\right]
$$
over all $\xi=(x,y)\in\C^{2k}\times\C^{2k}$ with norm 1, where $X=\Re x, Y=\Re y$.
We have
$$
(|X|^2-|Y|^2)^2+4|X|^2|Y|^2/\cos^2\vt
=
(|X|^2+|Y|^2)^2+4|X|^2|Y|^2\tan^2\vt.
$$
Consequently, the expression in the square brackets is always nonpositive, thus we say that we may take $x,y$ to be real. Thus we showed that 
$
\lambda_p(B\oplus B^*)
$
is the minumum of 
$$
1+
\frac{p-2}{2}
\left(
1-\sqrt{1+4|x|^2|y|^2\tan^2\vt}
\right)
$$
over all $\xi=(x,y)\in\C^{2k}\times\C^{2k}$ with $|x|^2+|y|^2=1$.
Clearly, the minimum in question is 
$$
1+
\frac{p-2}{2}
\left(
1-\sqrt{1+\tan^2\vt}
\right)
=
1+
\frac{p-2}{2}
\left(
1-1/\cos\vt
\right).
$$
This is $\lambda_p(B\oplus B^*)$. In order to obtain $\Delta_p(R\oplus R^*)$, we just need to multiply by $(2/p)\cos\vt$.
\end{proof}

\subsubsection{Representation of $\Phi_\omega$ for real $\omega$ in the case $m=2$}\ 

\medskip
Recall that the matrices $R_{\f,n}$ were introduced in \eqref{eq: R_theta,k}.
Denote also
$$
\aligned
\hskip 5.6pt
\cT=\cT_n& =
\left[
\begin{array}{rr}
1 & \ \ \ \ \ \
\\
  & 
\end{array}
\right]
\otimes I_n;
\endaligned
$$
this is $T_{1,0}$ in the notation of Theorem \ref{t: aI+bI}. 
These matrices act on $\C^n\times \C^n$.
For simplicity we will write just $R_\f$ and $\cT$, omitting the indication of the tensorization with the identity.

\begin{claim}
\label{cl: real vs. complex omega}
Fix $n\in\N$ and $\vt\in[-\pi,\pi)$. Then for every $\xi\in\Sigma_n^2$ and $\omega\in\Sigma^2$
there exist $\sigma\in\Sigma_n^2$ and $\omega_0\in S^1$ (that is, a \underline{real} $\omega_0\in\Sigma^2$) such that 
\begin{equation}
\label{eq: real omega}
\Re\Sk{R_{\vt,n}\xi}{\Phi_\omega\xi}_{\C_n^2}
=
\Re\Sk{R_{\vt,n}\sigma}{\Phi_{\omega_0}\sigma}_{\C_n^2}.
\end{equation}
We can take $\omega_0=(1,0)$.
\end{claim}

\begin{proof}
By using the definition of $R_{\vt,n}$ in the down-notation, we get
$$
\aligned
\Re\Sk{R_{\vt,n}\xi}{\Phi_\omega\xi}_{\C_n^2}
& = 
\sum_{j=1}^n
\Re\Sk{R_{\vt}\xi_j}{\omega}_{\C^2}
\Re\Sk{\xi_j}{\omega}_{\C^2}\\
& = 
\sum_{j=1}^n
\Sk{U_{\vt}\wt\xi_j}{\wt\omega}_{\R^4}
\Sk{\wt\xi_j}{\wt\omega}_{\R^4},
\endaligned
$$
where
$$
U_\vt=
\begin{bmatrix}[r]
R_\vt & \\
 & R_\vt
\end{bmatrix}\in\R^{4\times 4}.
$$

Since $\sk{U_\vt\sigma}{\sigma}$ is the same ($\cos\vt$) for all $\sigma\in\R^4$, $|\sigma|=1$, there exists an orthogonal map $\tau:\R^4\rightarrow\R^4$ such that for $\wt\omega_0:=(1,0,0,0)\in\R^4$ we have
$$
\aligned
\wt\omega_0 & \longmapsto \wt\omega\\
 U_\vt^T \wt\omega_0 & \longmapsto U_\vt^T \wt\omega.
\endaligned
$$
Then for $\wt\sigma_j:=\tau^T\wt\xi_j$ we have
$$
\aligned
\Sk{U_{\vt}\wt\xi_j}{\wt\omega} & = \Sk{U_{\vt}\wt\sigma_j}{\wt\omega_0}\\
\Sk{\wt\xi_j}{\wt\omega} & = \Sk{\wt\sigma_j}{\wt\omega_0}.
\endaligned
$$

We thus showed that for $\sigma:=(\sigma_1,\hdots,\sigma_n)$ we get \eqref{eq: real omega}.
It remains to observe that now $\omega_0=(1,0)$ is {\it real} and 
$
|\sigma_j|
=
|\wt\sigma_j|
=
|\wt\xi_j|
=
|\xi_j|
$
by the orthogonality of $\tau$ (and hence $\sigma\in\Sigma_n^2$ if $\xi\in\Sigma_n^2$), which proves the claim.
\end{proof}

\begin{lemma}
\label{l: Phi for real omega}
Let $m=2$ and $\omega\in\Sigma^2$ be real, that is, $\omega=(\cos\f,\sin\f)$ for some $\f\in[0,2\pi)$.
Then for $\xi=(\xi^1,\xi^2)\in\C_n^2$ we have
\begin{align}
\Phi_\omega\xi & =(R_\f\circ \cT\circ R_{-\f})(\Re\xi). \notag
\end{align}
Here $R_\f$ should be understood as $R_{\f,n}$ and likewise $R_{-\f}$.
\end{lemma}

\begin{proof}
We have 
$
(\Phi_\omega\xi)^{\text{up}}
=
\left(
	\omega^1\Re[\xi,\omega],
	\omega^2\Re[\xi,\omega]
\right)
$;
here $[\xi,\omega]=\overline\omega^1\xi^1+\overline\omega^2\xi^2$.
Write $\xi^1=x,\xi^2=y$, thus $x,y\in\C^n$. Since $\omega^1=\cos\f,\omega^2=\sin\f$ are real, we get
$$
\Re[\xi,\omega]=\Re(x\cos\f+y\sin\f)=(\Re x)\cos\f+(\Re y)\sin\f\in\R^n.
$$
Consequently,
$$
\aligned
(\Phi_\omega\xi)^{\text{up}}
& =
\left(
\cos\f((\Re x)\cos\f+(\Re y)\sin\f),\sin\f((\Re x)\cos\f+(\Re y)\sin\f)
\right)\\
& = 
\left[
\begin{array}{rr}
\cos\f \cdot \cos\f & \cos\f \cdot \sin\f 
\\
\sin\f \cdot \cos\f  & \sin\f \cdot \sin\f 
\end{array}
\right]
\begin{bmatrix}
\Re x\\
\Re y
\end{bmatrix},
\endaligned
$$
the matrix understood as being composed of and acting on blocks (that is, strictly speaking, we should have tensorized with $I_n$). 
Finally, from \eqref{eq: ab+cd formula} we conclude that
\begin{alignat*}{2}
\begin{bmatrix}[r]
\cos\f \cdot \cos\f & \cos\f \cdot \sin\f 
\\
\sin\f \cdot \cos\f  & \sin\f \cdot \sin\f 
\end{bmatrix}
& = R_\f \cT R_{-\f}.
\tag*{\qedhere}
\end{alignat*}
\end{proof}

\subsubsection{Direct sum of rotations in the lower representation from \eqref{eq: A up-down}}

Here we present, again in the case $m=2$, another example devised with rotation matrices as building blocks.
Its main difference with respect to the previous two examples, featuring in Theorems \ref{t: aI+bI} and \ref{t: Delta_p of U_theta,k}, is that in the upper representation it is {\it not} a direct sum of elliptic matrices.

\begin{theorem}
\label{t: Delta of V_theta}
Choose $\vt\in(-\pi/2,\pi/2)$.  Following the notation of Section \ref{eq: matrices, tensors},
define the tensor $V_\vt=V_{\vt,n}$ as acting on $\C_n^2$ by 
\begin{equation}
\label{t: V_theta def}
V_\vt^{\rm{up}}
:=
\left[
\begin{array}{cc}
V^{1,1} & 
V^{1,2}\\
V^{2,1} & 
V^{2,2}
\end{array}
\right]
=
\begin{bmatrix}[r]
(\cos\vt) I_n 
& 
(-\sin\vt) I_n 
\\
(\sin\vt) I_n 
& 
(\cos\vt) I_n
\end{bmatrix}
=
R_{\vt,n}
\end{equation}
or, equivalently,
$$
(V_\vt)_{\rm{down}}
:=
\begin{bmatrix}
R_\vt & 
& \\
& \ddots & \\
& 
& R_\vt
\end{bmatrix},
$$
with $n$ copies of $R_\vt$ on the diagonal.

Then
$$
\Delta_p(V_{\vt})=\cos\vt-|1-2/p|.
$$
\end{theorem}

\begin{proof}
We will again use \eqref{eq: lambda_p} and \eqref{eq: Delta vs. lambda}. 
Building upon the representation \eqref{t: V_theta def}, we take $\xi\in\Sigma_n^2$, $\omega\in\Sigma^2$  and calculate the expression in \eqref{eq: lambda_p}, where for simplicity we write just $R_\vt$ instead of $R_{\vt,n}=R_\vt\otimes I_n$. 
We get
$$
\aligned
\Re\sk{R_\theta\xi}{\xi+(p-2)\Phi_\omega\xi}
=
|\xi|^2\cos\vt+(p-2)\Re\sk{R_\theta\xi}{\Phi_\omega\xi}.
\endaligned
$$
By Claim \ref{cl: real vs. complex omega}, we may assume that $\omega$ is in fact real.
Now Lemma \ref{l: Phi for real omega} gives
$$
\aligned
\Re\sk{R_\theta\xi}{\Phi_\omega\xi}
& =\sk{\Re(R_\theta\xi)}{(R_\f \cT R_{-\f})(\Re\xi)}_{\R_n^2}\\
& =\sk{R_{-\f} R_\theta(\Re\xi)}{(\cT R_{-\f})(\Re\xi)}_{\R_n^2}\\
& =\sk{ R_\theta(R_{-\f}(\Re\xi))}{\cT (R_{-\f}(\Re\xi))}_{\R_n^2}\\
& =\sk{\cT R_\theta\eta}{\eta}_{\R_n^2}\\
& =\sk{M_\theta\eta}{\eta}_{\R_n^2},
\endaligned
$$
where $\eta=R_{-\f}(\Re\xi)\in\R_n^2$ and
$$
M_\theta
=\frac12\left(\cT R_\theta+R_{-\theta}\cT\right)
=
\begin{bmatrix}
\cos\vt & -(\sin\vt)/2\\
-(\sin\vt)/2 & 
\end{bmatrix}
\otimes I_n.
$$
Observe that $|\eta|\leq|\xi|$ always, and $|\eta|=|\xi|$ for real $\xi$. 
We have
$$
\det(M_\theta-\lambda I_{2n})
=
\left(
\lambda^2-\lambda\cos\vt-\frac{\sin^2\vt}{4}
\right)^n.
$$
Consequently,
the eigenvalues of $M_\theta$ are $(\cos\vt\pm1)/2$, the smaller of the two being $(\cos\vt-1)/2$; when $n=1$, the eigenvectors corresponding to this eigenvalue are multiples of $(\sin\vt/2,\cos\vt/2)$. 

Suppose that $p\geq2$. 
By Proposition \ref{p: real symmetric}, we get
$$
\lambda_p(V_\theta)=
\cos\vt+(p-2)(\cos\vt-1)/2,
$$
hence multiplication by $2/p$ returns
$$
\Delta_p(V_\vt)
=
\cos\vt-(1-2/p).
$$
For $1<p<2$ the proof works just as well, under minor obvious modifications.
\end{proof}

\subsubsection{Comparison between $U_\vt$ and $V_\vt$}

For the sake of clarity, let us here present a comparison between the two rotation-type tensors, $U_{\vt,k}$ and $V_{\vt,k}$, that we studied in 
Theorems \ref{t: Delta_p of U_theta,k} and \ref{t: Delta of V_theta}, respectively.

Fix $\vt\in(-\pi/2,\pi/2)$ and $k\in\N$. 
By resorting to the notation \eqref{eq: R_theta,k} and recalling the figures from Remark \ref{r: up-to-down colors}, we may write 
$$
\aligned
(U_{\vt,k})^{\text{up}}
& =R_{\vt,k}\oplus R_{\vt,k}^*\\
(V_{\vt,2k})_{\text{down}}
& =R_{\vt,k}\oplus R_{\vt,k}.
\endaligned
$$
In other words and more explicitly,
$$
(U_{\vt,k})^{\text{up}}
=
\begin{bmatrix}[r]
\cos\vt & -\sin\vt & & \\
\sin\vt & \cos\vt & & \\
& & \cos\vt & \sin\vt \\ 
& & -\sin\vt & \cos\vt 
\end{bmatrix}
\otimes
I_k
$$
\hskip 10pt
\text{ and }
\hskip 10pt
$$
(V_{\vt,2k})^{\text{up}}
=
\begin{bmatrix}[r]
\cos\vt & & -\sin\vt & \\
& \cos\vt &  & -\sin\vt \\
\sin\vt & & \cos\vt & \\
& \sin\vt &  & \cos\vt
\end{bmatrix}
\otimes
I_k.
$$

\part{Roots}

\label{part: Roots}

In this section of the paper we aim to explain ideas and reasoning that led us to our definition of $p$-ellipticity for systems.

\section{Roots of $p$-ellipticity for systems: formal time-derivative of the $L^p$ norm of the semigroup}
 
\label{s: rootsyst}

We start by presenting an argument -- not fully rigorous -- involving the contractivity of the semigroups and which suggests that $p$-ellipticity is a natural condition for {\it scalar} elliptic operators in divergence form and with complex coefficients.
Then we will try to explain how mimicking this reasoning in the vector case paved way for our definition of $p$-ellipticity for {\it systems}. 
Finally, we complement these intuitive arguments with a rigorous statement on the contractivity of $\exp(-tL_A)$ on $L^p$ (Theorem \ref{t: contractivity}).

\medskip
We indicate the following key steps of that argument.

\begin{itemize}

\item
An operator semigroup $(T_t)_{t>0}$ on a Banach space $X$ is contractive if and only if for every $u\in X$, 
the norms $\nor{T_tu}_X$ are decreasing in $t$. When $X=L^p(\Omega)$, this means that 
\begin{equation}
\label{eq: Tt Lp contr}
-\frac{d}{dt}\int_{\Omega}\left|T_tu\right|^p\geq0,
\end{equation}
provided that the integrals are differentiable with respect to $t$.

\item
In the special case when the semigroup $(T_t)_{t>0}$ is generated by the divergence-form operator $L_A$, we may (formally, at least) differentiate under the integral sign, use the definition of an $L_A$-generated semigroup, integrate by parts and reshape the left-hand side in \eqref{eq: Tt Lp contr} as 
\begin{equation*}
p\,\Re\int_\Omega
\Sk{A\nabla  g}{\nabla\big(|  g|^{p-2}  g\big)},
\end{equation*}
where $g=g_t=e^{-tL_A}u$.

\item
We recognize $\nabla\big(|  g|^{p-2}  g\big)$, up to a constant factor, as the {\it Hessian matrix of the power function $F_p$}.
This is not too surprising since $F_p$ is a crucial part of the definition of the $L^p$ norm.

\item
Finally, we express the Hessian matrix of $F_p$ through a formula involving $\cI_p$ \cite[Remark 2.4]{CD-Potentials}, which thus clearly leads to $p$-ellipticity. 

\end{itemize}

In \cite[Theorem 1.3]{CD-DivForm} we showed that $p$-ellipticity indeed bears close ties with the $L^p$ contractivity of the semigroup (in particular, it implies it). In this sense, we put the above -- partially informal -- considerations onto firm ground.

\subsection{Informal calculations of the heat flow's derivative in the vector case}

Now we want to pursue the same path in the vector case.
As in the scalar case, the ensuing calculation is rather formal and serves as a motivation for introducing the $p$-ellipticity condition, whose consequences we then prove rigorously.

\subsection*{Meta-lemma}
{\it
Let $A$ be $p$-elliptic and $L=L_A$. Write
$$
\vec g=\vec g_t(x):=\left(e^{-tL}u\right)(x)\in\C^{m}.
$$
Assuming that the integrals below are well defined and that we may exchange the derivative and the integral,
we have
\begin{equation}
\label{eq: Tt Lp norm derivative}
-\frac{d}{dt}\int_{\Omega}\left|e^{-tL}u\right|^p
=
p\,\Re\int_\Omega
\Sk{A\nabla\vec g}{\nabla\big(|\vec g|^{p-2}\vec g\big)}_{\C^m_n}.
\end{equation}
}

\begin{proof}[{\bf Explanation.}]
Suppose that the following terms exist and the identities hold:
$$
\aligned
-\frac{d}{dt}\int_{\Omega}\left|e^{-tL}u\right|^p
& =
-\frac{d}{dt}\int_{\Omega}F_p
(\vec g)\\
& =
-\int_{\Omega}\frac{d}{dt}F_p
(\vec g).
\endaligned
$$

Let us elaborate on the integrand. Introduce the notation 
$$
\vec g=(g_1,\hdots,g_m)
$$
and
\begin{equation*}
g_j=\alpha_j+i\beta_j\,
\end{equation*}
with $\alpha_j,\beta_j\in\R$. 
We calculate
$$
\aligned
\frac{d}{dt}F_p(\vec g)
& =\frac{d}{dt}\widetilde F_p\underbrace{\big(\alpha_1(t),\hdots,\alpha_m(t),\beta_1(t),\hdots,\beta_m(t)\big)}_{=\widetilde{\vec g}\,=:\sigma(t)}\\
& =\sum_{j=1}^m
\left[\frac{\pd\widetilde F_p}{dx_j}\big(\sigma(t)\big)\cdot\alpha_j'(t)+\frac{\pd\widetilde F_p}{dy_j}\big(\sigma(t)\big)\cdot\beta_j'(t)\right]\\
& =\sum_{j=1}^m
\Re
\left[
\Big(\frac{\pd\widetilde F_p}{dx_j}+i\frac{\pd\widetilde F_p}{dy_j}\Big)\cdot(\alpha_j'-i\beta_j')
\right]\\
& =\sum_{j=1}^m
\Re
\sk{2\pd_{\bar z_j}\widetilde F_p}{g_j'}_\C\\
&
=2\Re
\Sk{ \big(\pd_{\bar z}\widetilde F_p\big) (\vec g) }{\frac{d\vec g}{dt}}_{\C^{m}}\\
&=-2\Re
\Sk{ \big(\pd_{\bar z}\widetilde F_p\big) (\vec g) }{L\vec g}_{\C^{m}}.
\endaligned
$$ 
Hence from \eqref{eq: L def} we get
\begin{equation*}
\aligned
-\frac{d}{dt}\int_{\Omega}\left|e^{-tL}u\right|^p
& =2\Re\int_{\Omega}
\Sk{ \big(\pd_{\bar z}\widetilde F_p\big) (\vec g) }{L\vec g}_{\C^{m}}\\
& =2\Re\int_\Omega
\Sk{A\nabla\vec g}{\nabla\big[\big(\pd_{\bar z}\widetilde F_p\big)(\vec g)\big]}_{\C^m_n}.
\endaligned
\end{equation*}
In order to conclude, observe that 
$\pd_{\bar z} \widetilde F_p(z) 
=
(p/2) |z|^{p-2}z$ for $z\in\C^m$.
\end{proof}

The ``meta-identity" \eqref{eq: Tt Lp norm derivative} now suggests that the term $\nabla\left(|f|^{p-2}f\right)$ deserves to be looked into. This led us to Theorem \ref{t: nabla |f|^{p-2}f} which is the centrepiece of the next section.

\section{Roots of $\cJ_{p,\omega}$: Calculation of $\nabla\left(|f|^{p-2}f\right)$}
\label{s: calculation}

We will see in Section \ref{s: Nittka theorem} that it is important to study the positivity of the expression
$$
\Re\a(f,|f|^{p-2}f)=\Re\Sk{A\nabla f}{\nabla\left(|f|^{p-2}f\right)}_{L^2(\Omega;\C_n^m)},
$$
where $\a$ is as in Section \ref{s: form def}. 
To this end, one needs to express in more tangible terms the rightmost term that appears above. In this section we aim to fulfill this task. This will eventually lead to the crucial identity \eqref{eq: dissip term formula}. 
We summarize our findings in the following identity.
It explains the significance of $\cJ_{p,\omega}$ introduced in
\eqref{eq: J_pvect} and featured in \eqref{eq: lambda_p via J_p}, 
and the reason behind its introduction. Scalar equivalents of \eqref{eq: nabla |f|^{p-2}f} appeared in \cite[(7.5)]{CD-DivForm} and \cite[Lemma 2.3]{CD-Potentials}.

\begin{theorem}
\label{t: nabla |f|^{p-2}f}
Suppose that $p > 1$ and that $f$ and $|f|^{p-2}f$ belong to $H^1_{loc}(\Omega;\C^m)$. Then
\begin{equation}
\label{eq: nabla |f|^{p-2}f}
\nabla\left(|f|^{p-2}f\right)
=
|f|^{p-2}
\cJ_{p,f/|f|}\nabla f.
\end{equation}
Consequently,
$$
\left|
\nabla\left(|f|^{p-2}f\right)
\right|
\sim
|f|^{p-2}
|\nabla f|.
$$
\end{theorem}

The rest of this section is devoted to gradually proving this statement. 

\subsection{Power functions}
Take $p>1$. Introduce the {\it power function} $F_p$, defined on complex vectors of any (finite) length $k$, by
$$
\begin{array}{cccc}
F_p : & \C^{k} & \longrightarrow & [0,\infty)\\
& \xi  & \longmapsto & |\xi|^p\,.
\end{array}
$$
The proof of Theorem \ref{t: nabla |f|^{p-2}f} will rest on recognizing the expression $|f|^{p-2}f$ as the $\bar\pd$-derivative of the power function. Namely, 
$
\pd_{\bar z} \widetilde F_p(z) 
=
(p/2) |z|^{p-2}z,
$
thus 
\begin{equation}
\label{eq: |f|^{p-2}f via F_p}
|f|^{p-2}f
=
\frac 2p
\big(\pd_{\bar z}\widetilde F_p\big)(f)
\,.
\end{equation}

\subsection{Chain rule}

We need to calculate the expressions of the type
$\pd_{x_j}\left[\left(\pd_{\bar z_k}F_p\right) (\vec g)\right]$.
For this purpose we recall the following theorem, whose proof was given in \cite[Theorem B.3]{CD-Potentials}.
Recall the notation from Section \ref{s: complex vs. real}.

\begin{theorem}
\label{t: chain rule scalar}
Let $\Omega\subset\R^n$ be open. Take functions $f=u+iv\in H^1_{loc}(\Omega)$ and $\Phi:\C\rightarrow\C$,
such that 
$\wt\Phi\in C^1(\R^2)$ and
\begin{center}
$\Phi\circ f=\wt\Phi(u,v)$ and $|\nabla \wt\Phi|(u,v)$ belong to $L_{loc}^2(\Omega)$\,.
\end{center}
Then $\Phi\circ f$ is weakly differentiable and we have 
\begin{equation*}
\nabla(\Phi\circ f) 
=\big(\pd_z\wt\Phi\big)(u,v)\cdot\nabla f
+\big(\pd_{\bar z}\wt\Phi\big)(u,v)\cdot\nabla\bar f.
\end{equation*}
\end{theorem} 

Here we apply it to vector-valued case.
We will use the notation from Section \ref{s: complex vs. real}.
Let $G=G(u_1,\hdots,u_m,v_1,\hdots,v_m)\in C^2(\R^m\times\R^m)$ be a given complex function.
For $l,k\in\{1,\hdots,m\}$ denote 
\begin{equation*}
\label{Mirko - Mlad se momak}
d_{lk}^2G:=
\begin{bmatrix}
G_{u_lu_k} & G_{v_lu_k} \\
G_{u_lv_k} & G_{v_lv_k} 
\end{bmatrix}.
\end{equation*}
If now $F:\C^m\rightarrow\C$ is such that $\wt F= F\circ \cV^{-1}\in C^2(\R^m\times\R^m)$, define
$$
D_{lk}^2F
:=
d_{lk}^2(F\circ \cV^{-1})\circ \cV.
$$

\begin{proposition}
\label{p: chain rule vector}
Assume that:
\begin{itemize}

\item
$\Omega\subset\R^n$ is an open set and $\vec g=(g_1,\hdots,g_m)\in H^1_{loc}(\Omega;\C^m)$. 

\noindent
Write $\vec g=\vec\alpha+i\vec\beta$ for $\vec\alpha,\vec\beta:\Omega\rightarrow\R^m$;
thus $\cV(\vec g)=(\vec\alpha,\vec\beta)$.

\medskip
\item
$F:\C^m\rightarrow\R$ is such that
$\wt F=\wt F(u_1,\hdots,u_m,v_1,\hdots,v_m)
\in C^2(\R^m\times\R^m)$.

\noindent
Introduce $\Psi=(\Psi_1,\hdots,\Psi_m):\C^m\rightarrow\C^m$ by $\widetilde\Psi=\pd_{\bar z}\wt F$.
Assume that
\begin{center}
$\Psi\circ \vec g=\wt\Psi(\vec\alpha,\vec\beta)$ and $|\nabla \wt\Psi|(\vec\alpha,\vec\beta)$ belong to $L_{loc}^2$\,.
\end{center}

\end{itemize}

Then $\Psi\circ\vec g$ is weakly differentiable and for $j\in\{1,\hdots,n\}$, $k\in\{1,\hdots,m\}$ we have 
$$
\cV\left(
\pd_{x_j}\left[\big(\pd_{\bar z_k}\wt F\big)\circ\cV(\vec g)\right]
\right)
=
\frac12
\sum_{l=1}^m
(D_{lk}^2F)(\vec g)\cdot\widetilde{\pd_{x_j}g_l}\,.
$$
\end{proposition}

\begin{proof}
We are differentiating the composition
$$
\Omega\xrightarrow{\vec g=(g_1,\hdots,g_m)}\C^{m}\xrightarrow{\Psi_k}\C\,.
$$
Let $j,k$ be as in the formulation of the theorem. Write $g_l=\alpha_l+i\beta_l$. 
From the definition of $\Psi$ we get $\Psi_k:\C^m\rightarrow\C$ and
$$
\wt\Psi_k=\pd_{\bar z_k}\wt F:\R^m\times\R^m\rightarrow\C\,.
$$
Observe that $\big(\pd_{\bar z_k}\wt F\big)\circ\cV(\vec g)=(\Psi_k\circ\cV^{-1})\circ(\cV\circ \vec g)=\Psi_k\circ\vec g$.
Applying a direct generalization of Theorem \ref{t: chain rule scalar} -- to mappings to and from $\C^m$ rather than $\C$ -- yields
\begin{alignat*}{2}
\frac{\pd}{\pd x_j}\left[\Psi_k\circ\vec g\right]
& =
\sum_{l=1}^m\left[\left(\pd_{z_l}\pd_{\bar z_k}\wt F\right)(\vec\alpha,\vec\beta)\cdot\frac{\pd g_l}{\pd x_j}+\left(\pd_{\bar z_l}\pd_{\bar z_k}\wt F\right)(\vec\alpha,\vec\beta)\cdot\frac{\pd \overline g_l}{\pd x_j}\right]\\
\label{Nikon D610}
& =
\sum_{l=1}^m\bigg[
\frac14(\pd_{u_l}-i\pd_{v_l})(\pd_{u_k}+i\pd_{v_k})\wt F(\cV(\vec g))\cdot\big(\pd_{x_j}\alpha_l+i\,\pd_{x_j}\beta_l\Big)\\
& 
\hskip 27pt
+\frac14(\pd_{u_l}+i\pd_{v_l})(\pd_{u_k}+i\pd_{v_k})\wt F(\cV(\vec g))\cdot\big(\pd_{x_j}\alpha_l-i\,\pd_{x_j}\beta_l\big)\bigg]\\
& =
\frac12
\sum_{l=1}^m
\bigg[
\left(
\pd_{x_j}\alpha_l\cdot \wt F_{u_lu_k} + \pd_{x_j}\beta_l\cdot \wt F_{v_lu_k}
\right)\\
& 
\hskip 33.3pt
+i
\left(
\pd_{x_j}\alpha_l\cdot \wt F_{u_lv_k} + \pd_{x_j}\beta_l\cdot \wt F_{v_lv_k}
\right)
\bigg](\cV(\vec g)).
\end{alignat*}
Consequently,
\begin{alignat*}{2}
\cV\left(
\frac{\pd}{\pd x_j}\left[\Psi_k\circ\vec g\right]
\right)
& =
\frac12
\sum_{l=1}^m
\begin{bmatrix}
\wt F_{u_lu_k}(\cV(\vec g)) & \wt F_{v_lu_k}(\cV(\vec g)) \\
\wt F_{u_lv_k}(\cV(\vec g)) & \wt F_{v_lv_k}(\cV(\vec g)) 
\end{bmatrix}
\begin{bmatrix}
\pd_{x_j} \alpha_l\\
\pd_{x_j} \beta_l
\end{bmatrix}.
\tag*{\qedhere}
\end{alignat*}
\end{proof}

\subsection{Hessian of power functions $F_p$}

In view of the identity \eqref{eq: |f|^{p-2}f via F_p}, our foremost interest in Proposition \ref{p: chain rule vector} lies in applying it to $F=F_p$.
This calls for an evaluation of the matrices $D_{lk}^2F_p$, which is the content of the next result.
It generalizes \cite[(5.5)]{CD-DivForm}; see also \cite[Remark 2.4]{CD-Potentials}.

\begin{lemma}
\label{l: proto I_p}
For any $\xi=(\xi_1,\hdots,\xi_m)\in\C^m\backslash\{0\}$ 
and $l,k\in\{1,\hdots,m\}$ such that $\xi_k,\xi_l\ne0$
we have
$$
(D_{lk}^2F_p)(\xi)
=
p|\xi|^{p-2}\cdot
R_{\arg\xi_k}
\left(
\delta_{kl}I+
(p-2)\frac{|\xi_k|}{|\xi|}\cdot\frac{|\xi_l|}{|\xi|}
\begin{bmatrix}
1 & & \\
 & &
\end{bmatrix}
\right)
R_{\arg \xi_l}^{-1}.
$$
Here $\delta_{kl}$ is the Kronecker delta, $I$ is the $2\times2$ identity matrix, while $R_\f$ is the standard rotation matrix in $\R^2$ by the angle $\f\in\R$, as specified in \eqref{eq: R_theta}.
Also, $(d_{lk}^2F_p)(0)=0$.
\end{lemma}

\begin{proof}
Let $\sigma=\wt\xi=(u_1,\hdots,u_m,v_1, \hdots, v_m)$. Recall that
$$
F_p(\xi)=\wt F_p(\wt \xi)=\left(u_1^2+v_1^2+\hdots+ u_m^2+v_m^2\right)^{p/2}.
$$
Then 
$$
\aligned
\big(\pd_{u_k}\wt F_p\big)(\sigma) &= pu_k|\sigma|^{p-2}
\\
\big(\pd_{v_k}\wt F_p\big)(\sigma)  &= pv_k|\sigma|^{p-2}.
\endaligned
$$
Consequently,
$$
\aligned
(d_{lk}^2\wt F_p)(\sigma)
&=
\begin{bmatrix}
\pd^2_{u_lu_k}\wt F_p & \pd^2_{v_lu_k}\wt F_p \\
\pd^2_{u_lv_k}\wt F_p & \pd^2_{v_lv_k}\wt F_p 
\end{bmatrix}\\
&=
\begin{bmatrix}
p|\sigma|^{p-2}\delta_{kl}+pu_k(p-2)u_l|\sigma|^{p-4} &  pu_k(p-2)v_l|\sigma|^{p-4} \\
pv_k(p-2)u_l|\sigma|^{p-4}   & p|\sigma|^{p-2}\delta_{kl} + pv_k(p-2)v_l|\sigma|^{p-4}
\end{bmatrix}\\
&=
p|\sigma|^{p-2}
\begin{bmatrix}
\delta_{kl} +(p-2)u_ku_l/|\sigma|^{2} &  (p-2)u_kv_l/|\sigma|^{2} \\
(p-2)u_lv_k/|\sigma|^{2}   & \delta_{kl} +(p-2)v_lv_k/|\sigma|^{2}
\end{bmatrix}.
\endaligned
$$
Therefore, with 
$$
\aligned
\xi_k
&=u_k+iv_k\\
&=|\xi_k|(\cos\f_k+i\sin\f_k),
\endaligned
$$
we obtained 
$$
\aligned
\frac{(d_{lk}^2\wt F_p)(\sigma)
}{p|\xi|^{p-2}}
& =
\delta_{kl}I_2+(p-2)\frac{|\xi_l|}{|\xi|}\cdot\frac{|\xi_k|}{|\xi|}
\begin{bmatrix}
u_l/|\xi_l|\cdot u_k/|\xi_k| &  v_l/|\xi_l|\cdot u_k/|\xi_k|\\
u_l/|\xi_l|\cdot v_k/|\xi_k|  & v_l/|\xi_l|\cdot v_k/|\xi_k|
\end{bmatrix}\\
&=
\delta_{kl}I_2+(p-2)\frac{|\xi_l|}{|\xi|}\cdot\frac{|\xi_k|}{|\xi|}
\begin{bmatrix}
\cos\f_k\cos\f_l & \cos\f_k\sin\f_l\\
\sin\f_k\cos\f_l & \sin\f_k\sin\f_l
\end{bmatrix}.
\endaligned
$$
We observe that
\begin{equation}
\label{eq: ab+cd formula}
\begin{bmatrix}
ac & ad\\
bc & bd
\end{bmatrix}
=
\begin{bmatrix}
a  \\
b 
\end{bmatrix}
\begin{array}{cc}
\pmb{\left[\vphantom{\frac{1}{2}}
\right.}
\, c & d\, 
\pmb{\left.\vphantom{\frac{1}{2}}\right]}\\
&
\end{array}
=
\begin{bmatrix}[r]
a & -b\\
b & a
\end{bmatrix}
\begin{bmatrix}
1 & & \\
 & &
\end{bmatrix}
\begin{bmatrix}
c & d\\
-d & c
\end{bmatrix},
\end{equation}
which means that
$$
\aligned
\frac{(d_{lk}^2\wt F_p)(\sigma)}{p|\xi|^{p-2}}
&=
\delta_{kl}I_2+
(p-2)\frac{|\xi_k|}{|\xi|}\cdot\frac{|\xi_l|}{|\xi|}
R_{\arg\xi_k}\Pi_1R_{\arg\xi_l}^{-1}\\
&=
R_{\arg\xi_k}
\left(
	\delta_{kl}I+
	(p-2)\frac{|\xi_k|}{|\xi|}\cdot\frac{|\xi_l|}{|\xi|}
	\Pi_1
\right)
R_{\arg\xi_l}^{-1}.
\endaligned
$$
Here 
$$
\Pi_1:=
\begin{bmatrix}
1 & &   \\
   & &
\end{bmatrix}
$$
represents the operator taking the real part of a complex number. Alternatively, it is the projection on the first variable ($x$) in $\C$. Finally, recall that $(d_{lk}^2\wt F_p)(\sigma)=(D_{lk}^2F_p)(\xi)$.
\end{proof}

\subsection{Inversion formul\ae}
\label{s: Inversion formulae}

Let $p>1$.
Define, for $f\in\C^m\backslash\{0\}$ and $F\in\C^m_n$, 
\begin{equation*}
\aligned
\tau_p(f) & = |f|^{p-2}f\\
\Lambda_p(f,F) 
& = |f|^{p-2}\cJ_{p,f/|f|}F.
\endaligned
\end{equation*}
In this notation, formula \eqref{eq: nabla |f|^{p-2}f} reads $\nabla\left(\tau_p(f)\right)=\Lambda_p(f,\nabla f)$.

We may easily verify the following pair of {\it inversion formul\ae}.
We remark that when $m=1$, they were already observed in \cite[p.103]{CD-Potentials}.

\begin{corollary}
Let $1/p+1/q=1$. Then for any $f\in\C^m\backslash\{0\}$ and $F\in\C^m_n$ we have
\begin{align}
f  & = \tau_q\big(\tau_p(f)\big).  
   \label{eq: tp}\\
F  & = \Lambda_q\big(\tau_p(f),\Lambda_p(f,F) \big).
   \label{eq: Lp}
\end{align}
\end{corollary}
\begin{proof}
Indeed, the first identity is straightforward, while for the second one we 
apply Proposition \ref{p: Phi properties} {\it \ref{p: Phi prop 3}} and calculate
\begin{alignat*}{2}
\Lambda_q\big(\tau_p(f),\Lambda_p(f,F) \big)
& =\left|\tau_p(f)\right|^{q-2}\cJ_{q,\tau_p(f)/|\tau_p(f)|}\Big(|f|^{p-2}\cJ_{p,f/|f|}F\Big)\\
& =|f|^{(p-1)(q-2)+(p-2)}\cJ_{q,f/|f|}(\cJ_{p,f/|f|}F)\\
& =F.
\tag*{\qedhere}
\end{alignat*}
\end{proof}

\subsection{Proof of Theorem \ref{t: nabla |f|^{p-2}f}}

First assume that \framebox{$p\geq2$}, as then the function $F_p$ is of class $C^2$, thus fitting the assumptions of Proposition \ref{p: chain rule vector}.

Choose $\xi=(\xi_1,\hdots,\xi_m)\in\C^{m}$ and indices $k,l\in\{1,\hdots,m\}$. Define the ($\R$-linear) operator 
$\cE_{p,\xi,k,l}:\C\rightarrow\C$ by
$$
\cE_{p,\xi,k,l}(z):=\delta_{kl}z+(p-2)\frac{|\xi_k|}{|\xi|}\cdot\frac{|\xi_l|}{|\xi|}\,\Re z
\hskip 30pt\forall z\in\C.
$$
It serves the purpose of restating Lemma \ref{l: proto I_p}, as we have
$$
\delta_{kl}I+
(p-2)\frac{|\xi_k|}{|\xi|}\cdot\frac{|\xi_l|}{|\xi|}
\begin{bmatrix}
1 & & \\
 & &
\end{bmatrix}
=
\cV \circ \cE_{p,\xi,k,l} \circ \cV^{-1}.
$$
Thus the action on $\R^2$ of the matrix on the left is represented by the action on $\C$ of the operator $\cE_{p,\xi,k,l}$.
Together with Proposition \ref{p: chain rule vector} and Lemma \ref{l: proto I_p} 
and the identity $R_\f(\cV(z))=\cV(e^{i\f} z)$, we now obtain
$$
\pd_{x_j}\left[\big(\pd_{\bar z_k}\wt F\big)\circ\cV(\vec g)\right]
=
\frac  p2|\vec g|^{p-2}
\sum_{l=1}^m
e^{i\arg g_k}
\cE_{p,\vec g,k,l}
\big(e^{-i\arg g_l}
(\pd_{x_j}g_l)\big)\,.
$$
Finally, by observing that
$$
\frac{g_k\bar g_l}{|g_k||g_l|}\delta_{kl}=\delta_{kl}
$$
we get
\[
\aligned
\sum_{l=1}^m
e^{i\arg g_k}
\cE_{p,\vec g,k,l}
\big(e^{-i\arg g_l}
(\pd_{x_j}g_l)\big)
& =
\sum_{l=1}^m
\frac{g_k}{|g_k|}
\left(\delta_{kl}I+(p-2)\frac{|g_k|}{|\vec g|}\cdot\frac{|g_l|}{|\vec g|}\,\Re \right)
\Big(\frac{\bar g_l}{|g_l|}
(\pd_{x_j}g_l)\Big)\\
& =
\pd_{x_j}g_k+(p-2)\frac{g_k}{|\vec g|}
\sum_{l=1}^m
\Re
\Big(\frac{\bar g_l}{|\vec g|}
(\pd_{x_j}g_l)\Big).
\endaligned
\]

The case \framebox{$1<p<2$} is now settled by repeating verbatim (except for replacing $\C$ by $\C^m$ at several places) the trick from \cite[p.103]{CD-Potentials}. It relies on the inversion formul\ae\ from Section \ref{s: Inversion formulae}.  
For the sake of the reader's convenience we repeat it here. 

Assume that $1/p+1/q=1$.
Write $g=\tau_p(f)$. Because of \eqref{eq: tp}, the assumption 
$f,\tau_p(f)\in H^1_{loc}(\Omega)$
now reads 
$g, \tau_q(g)\in H^1_{loc}(\Omega)$. 
Since $q>2$, we may apply the part proven so far, which gives
that $f=\tau_q(g)$ is weakly differentiable, moreover,  
$
\nabla f = \nabla\left(\tau_q(g)\right)=\Lambda_q(g,\nabla g).
$
Now \eqref{eq: Lp} returns 
$$
\nabla\left(\tau_p(f)\right)
= \nabla g 
= \Lambda_p\big(\tau_q(g),\Lambda_q(g,\nabla g) \big) 
= \Lambda_p\big(f,\nabla f \big),
$$
as desired. 
These reasonings are surely valid on 
$\{f,g\ne0\}=\{f\ne0\}=\{g\ne0\}$.
On the other hand, $\nabla g=0$ on $\{f,g=0\}$ by a known result.

The last statement now follows from Corollary \ref{c: norm of J}.
\hfill \qed

\part{Applications of $p$-ellipticity for systems}

\label{part: Applications}

The syntagm ``applications of $p$-ellipticity" naturally refers to theorems where $p$-ellipticity features prominently, as a basic (usually -- but not always -- sufficient) condition. 
In the scalar case, such theorems were in recent years proven in a number of papers by various authors, see Section \ref{s: scalar p-ell}.

As we have seen throughout this paper, the vector case notably differs from the scalar one because the algebraic condition is considerably more involved, owing to the ``noncommutative" action of the sphere $\Sigma^m$. This said, it is understandable that any fundamentally new type of an estimate shall first be tested in the scalar case. 

Hence, the main challenge addressed in this setting is to carry the ``scalar" theorems to the vector case while trying to make use in a natural way of the generalization of $p$-ellipticity for systems introduced and explored in this paper. This is what we mean by the notion of ``applications of $p$-ellipticity" for systems. 

In order to keep the size of this paper manageable, we shall limit our attention to one example only -- the semigroup contractivity on $L^p$ (Theorem \ref{t: contractivity}). 
Other examples, most notably bilinear embedding (whose pursuit in the scalar complex case gave rise to our introduction of $p$-ellipticity in \cite{CD-DivForm}) shall be studied elsewhere.

\section{$L^p$ contractivity of semigroups}

\subsection{Semigroup properties on $L^p$}

The purpose of this section is to prove Theorem \ref{t: contractivity}. It generalizes the implication ${\mathit(a)}\Rightarrow{\mathit(b)}$ of \cite[Theorem 1.3]{CD-DivForm}, in the sense of proving it in the context of elliptic systems. See also \cite[Proposition 1]{CD-Mixed}, where the implication mentioned above was proven in special cases (a)-(c) from Section \ref{s: bdry conds}, but still in the scalar setting.  

The proof of Theorem \ref{t: contractivity} is based on the one from \cite{CD-Potentials}, which itself built upon \cite{CD-DivForm}. In all of those cases, we crucially rely on a criterion by Nittka (Theorem \ref{t: Nittka theorem}).
Likewise, we highlight the role played here by Theorem \ref{t: nabla |f|^{p-2}f}.

\begin{theorem}
\label{t: contractivity}
If $p\in(1,\infty)$, $\Omega\subset\R^n$ is an open set and $A$ is $p$-elliptic tensor acting on $\C_n^m$, then $\big(e^{-tL_A}\big)_{t>0}$ extends from $(L^2\cap L^p)(\Omega;\C^m)$ to a contractive semigroup on $L^p(\Omega;\C^m)$.
\end{theorem}

\medskip
Let us cite the following result by Arendt, ter Elst and Sauter \cite{AtES23}: 

\begin{theorem}\cite[Theorem 3.1]{AtES23}
Let $\Omega\subset\R^n$ be open, $H$ a separable Hilbert space and 
$A=[a_{jk}]$ an elliptic matrix such that $a_{jk}:\Omega\rightarrow\cL(H)$. If
\begin{equation}
\label{eq: AtES}
2|1-2/p|+|1-2/p|^2\leq\frac{\lambda(A)}{\Lambda(A)},
\tag{AtES}
\end{equation}
then, assuming Neumann boundary conditions, $e^{-tL_A}$ extends to a contraction semigroup on $L^p(\Omega,H)$.
\end{theorem}

The context of their theorem is more general than ours, in the sense that it includes our setting as the special case $H=\C^m$. 
In that case, however, by \eqref{eq: vector decoupling}, the range of contractivity given by Theorem \ref{t: contractivity} improves the range provided by \eqref{eq: AtES}. 
Indeed, by the \eqref{eq: mu bilateral est}, the \eqref{eq: AtES} condition implies, for $p\ne2$,
$$
|1-2/p|
<
2|1-2/p|+|1-2/p|^2\leq\frac{\lambda(A)}{\Lambda(A)}
\leq\mu^\pm(A),
$$
which is, by Theorem \ref{t: decoupling}, exactly $p$-ellipticity of $A$.

Recently, Böhnlein and Egert studied $L^p$ boundedness of the said semigroup. 
They formulate their main theorem in terms of the ``distance function" (a measure of $L^\infty$ distance from the identity to the set of all positive multiples of $A$); see \cite{BE24}.

\bigskip
The rest of this section is devoted to the proof of Theorem \ref{t: contractivity}.

\subsection{Nittka's invariance criterion for semigroups}
\label{s: Nittka theorem}

In his excellent paper \cite{Nittka}, Nittka specialized Ouhabaz's invariance theorem \cite[Theorem 2.1]{Ouh1996}, see also \cite[Theorem 2.2]{O}, to the case of $\cH=L^2(\Omega,\mu)$ and $K=\mn{u\in L^2\cap L^p}{\nor{u}_p\leq 1}$, as a result obtaining a remarkably simple characterization of the $L^p$ contractivity of semigroups generated by sesquilinear forms. 
In order to state Nittka's result in a concise manner, we need another definition.

\subsection*{$L^{p}$-dissipativity of forms}
The notion of {\it $L^{p}$-dissipativity} of sesquilinear forms was introduced by Cialdea and Maz'ya in \cite[Definition~1]{CM} for the case of forms defined on $C_c^1(\Omega)$ and associated with complex matrices. 
We extended that notion in \cite[Definition 7.1]{CD-DivForm}. Here we give another extension, by treating forms acting on vector-valued $L^2$ spaces. 

\begin{definition}
Let $\Omega$ be a measure space, $m$ a natural number, $\gotb$ a sesquilinear form defined on the domain $\cD(\gotb)\subset L^2(\Omega;\C^m)$ and $1<p<\infty$. Denote 
$$
\cD_p(\gotb):=\mn{u\in\cD(\gotb)}{|u|^{p-2}u\in\cD(\gotb)}. 
$$
We say that $\gotb$ is {\it $L^{p}$-dissipative} if 
$$
\Re\gotb\left(u,|u|^{p-2}u\right)\geq0
\hskip40pt
\forall\,u\in\cD_p(\gotb).
$$
Let $q$ be given by $1/p+1/q=1$. By \eqref{eq: tp}, we have $u\in\cD_p(\gotb)$ if and only if $\tau_p(u)=|u|^{p-2}u\in\cD_q(\gotb)$.
\end{definition}

The following theorem is due to Nittka \cite[Theorem 4.1]{Nittka}. 
It was stated for scalar ($\C$-valued) functions, but the proof there works just as well for the vector ($\C^m$-valued) ones. See also \cite[Theorem 1.1]{AtES23} for a generalization to Hilbert-space-valued functions. 
We also remark that Nittka formulated his result for sectorial forms, but it seems that sectoriality is not needed for our version of his result, since it is not needed for Ouhabaz's criterion \cite[Theorem 2.2]{O} on which Nittka's own criterion is based. Of course, the forms we are dealing with in this paper are all sectorial anyway.

\begin{theorem}
\label{t: Nittka theorem} 
Let $(\Omega,\mu)$ be a measure space. 
Suppose that the sesquilinear form $\a$ on $L^2=L^2(\Omega,\mu;\C^m)$ is densely defined, accretive, 
continuous and closed. Let $\cL$ be the operator associated with $\a$ in the sense of \cite[Section 1.2.3]{O}.

Take $p\in(1,\infty)$ and define $B^p:=\mn{u\in L^2\cap L^p}{\nor{u}_p\leq 1}$. Let $\bP_{B^p}$ be the orthogonal projection $L^2\rightarrow B^p$. 
Then the following assertions are equivalent: 
\begin{itemize}
\item
$\nor{\exp(-t\cL)f}_p\leq\nor{f}_p$ for all $f\in L^2\cap L^p$ and all $t\geq0$;
\item
$\cD(\a)$ is invariant under $\bP_{B^p}$ and 
$\text{\gotxii a}$ is {\it $L^p$-dissipative}.
\end{itemize}
\end{theorem}

As a consequence of Theorem \ref{t: nabla |f|^{p-2}f} we get, with $p>1$ and $\a$ as in Section \ref{s: the form},
\begin{equation}
\label{eq: dissip term formula}
\Re\a\left(f,|f|^{p-2}f\right)
=
\int_\Omega
| f |^{p-2}\,
\Re\sk{A\nabla f }{\cJ_{p, f /| f |}\nabla f }_{\C_n^m}
\hskip 30pt
\forall f\in\cD_p(\a).
\end{equation}

Therefore if this left-hand side of \eqref{eq: dissip term formula} is to be positive, as required by Nittka's theorem in order to obtain $L^p$ contractivity, it is {\it natural} to impose the {\it scalar} condition $\lambda_p(A)>0$, with $\lambda_p(A)$ as in \eqref{eq: lambda_p via J_p}.

\subsection{Proof of Theorem \ref{t: contractivity}}

This is a consequence of Nittka's theorem (Theorem \ref{t: Nittka theorem}). 
In order to apply it, we must check the following: 
\begin{enumerate}[i)]
\item
\label{Nittka cond 1}
$\cD(\a)$ is invariant under $\bP_{B^p}$;
\item
\label{Nittka cond 2}
$\a$ is $L^p$-dissipative. 
\end{enumerate}
We start with \ref{Nittka cond 1}. 
Let 
$$
{-\mathbf \Delta}:=-
\underbrace{\Delta\oplus\hdots\oplus\Delta}_{m}
$$ 
denote the operator associated with the form $\a_{{\mathbf I},\oV}$, that is, the form generated by the identity tensor (${\mathbf I}^{\alpha,\beta}:=\delta_{\alpha,\beta}I_n$). 
By the basic assumption \eqref{eq: invariance} and \cite[Theorem 4.31 {\it 2)}]{O}, the semigroup 
$$
\left(
e^{-t({-\mathbf \Delta})}
=
e^{-t(-\Delta)}
\oplus\hdots\oplus
e^{-t(-\Delta)}
\right)_{t>0}
$$
is contractive on $L^\infty(\Omega)$, since so is each of its direct summands. 
By duality, it is then also contractive on $L^1$. 
Thus, by interpolation with the $L^2$-estimates, the semigroup in question is also contractive on $L^p(\Omega)$ for all $p\in[1,\infty]$. Hence Nittka's Theorem \ref{t: Nittka theorem} gives that $\cD(\a_{{\mathbf I},\oV})$ is invariant under $\bP_{B^p}$. Now use that $\cD(\a_{{\mathbf I},\oV})=\cD(\a)$.

The statement \ref{Nittka cond 2} follows from \eqref{eq: dissip term formula} and \eqref{eq: lambda_p via J_p}. Indeed, if $A\in\cA_{p}^m(\Omega)$, then
\begin{alignat*}{2}
\Re\a\left(f,|f|^{p-2}f\right)
\geq
\lambda_p(A)
\int_\Omega
| f |^{p-2}|\nabla f|^2
\geq0
\hskip 30pt
\forall f\in\cD_p(\a).
\tag*{\qedsymbol}
\end{alignat*}

\section*{Statement on use of artificial intelligence}

All the main results, proofs, and research directions in this paper were developed by the authors without the use of large language models (LLMs). GPT-6.1 Sol and GPT-6 Astra (OpenAI) were used on the completed manuscript for proofreading, language polishing, and identifying textual or notational inconsistencies and passages lacking detail. All suggested corrections were reviewed and verified by the authors, who take full responsibility for the paper.

\section*{Acknowledgements}

A. Carbonaro was partially supported by the ``National Group for Mathematical Analysis, Probability and their Applications'' (GNAMPA-INdAM).

O. Dragi\v{c}evi\'c was partially supported by the Slovenian Research Agency, ARIS (research grants J1-1690, 
J1-3005, N1-0237 and research program P1-0291).

\bibliography{biblio_total}{}

    \bibliographystyle{alpha}

\end{document}